\UseRawInputEncoding
\documentclass[a4paper,oneside,11pt]{article}
\usepackage{amsmath}
\usepackage{amssymb}
\usepackage{titlesec}
\usepackage{graphicx}
\usepackage{fancyhdr}
\usepackage{latexsym}
\usepackage{mathrsfs}
\usepackage{booktabs}
\usepackage{mathrsfs}
\usepackage{mathtools}
\usepackage{amsthm}
\usepackage{wasysym}
\usepackage{cite}
\usepackage{color}
\usepackage{geometry}
\numberwithin{equation}{section}

\newtheorem{theorem}{Theorem}[section]
\newtheorem{lemma}{Lemma}[section]
\newtheorem{remark}{Remark}[section]

\newcommand{\n}{\rho}

\renewcommand{\div}{ {\rm div }  }

\renewcommand{\r}{\mathbb{R}}

\newcommand{\bt}{\begin{theorem}}
\newcommand{\bl}{\begin{lemma}}
\newcommand{\el}{\end{lemma}}
\newcommand{\et}{\end{theorem}}
\newcommand{\ga}{\gamma}
\newcommand{\thatsall}{\hfill$\Box$}
\newcommand{\curl}{{\rm curl} }

\newcommand{\ve}{\varepsilon}
\newcommand{\la}{\label}

\newcommand{\bn}{\begin{eqnarray}}
\newcommand{\en}{\end{eqnarray}}
\newcommand{\bnn}{\begin{eqnarray*}}
\newcommand{\enn}{\end{eqnarray*}}

\newcommand{\bnnn}{\begin{eqnarray*}}
\newcommand{\ennn}{\end{eqnarray*}}

\newcommand{\ba}{\begin{aligned}}
\newcommand{\ea}{\end{aligned}}
\newcommand{\be}{\begin{equation}}
\newcommand{\ee}{\end{equation}}
\def\O{{\Omega }}

\def\norm[#1]#2{\|#2\|_{#1}}

\def\la{\label}

\def\na{\nabla}

\makeatletter      
\@addtoreset{equation}{section}
\makeatother       
\date{}
\title{ Global Classical Solutions to  Compressible Navier-Stokes Systems with Vacuum under Non-Slip Boundary Conditions}
\author{Xinyu F{\small AN}$^{a}$,  Jing L{\small I}$^{b,c}$ \thanks{Email addresses:  fanxinyu17@mails.ucas.edu.cn (X. Y. Fan), ajingli@gmail.com  (J. Li)  }  \\  {\normalsize a.  Hua Loo-Keng Center for Mathematical Sciences, AMSS,}\\
{\normalsize Chinese Academy of Sciences, Beijing 100190, P. R. China;}\\
{\normalsize b. Department of Mathematics, }\\ {\normalsize  Nanchang University, Nanchang 330031, P. R. China;} \\ {\normalsize c. Institute of Applied Mathematics, AMSS,} \\ {\normalsize \&   Hua Loo-Keng Key Laboratory of Mathematics,}\\
{\normalsize  Chinese Academy of Sciences,    Beijing 100190,
P. R. China }}
\begin{document}
\maketitle
\begin{abstract}
This paper studies the global well-posedness of  classical solutions to the isentropic compressible Navier-Stokes equations in 3D domains $\Omega$ under non-slip boundary conditions. $\Omega$ will separate into the inner and boundary parts along a free surface: In the inner part,
the density is allowed to vanish and the gradient of it grows with an exponential rate when vacuum appears initially; while in the boundary part, no vacuum forms and the $L^2$-norm of higher order derivatives of the density remains uniformly bounded. We utilize the Lagrangian coordinates introduced by Christodoulou-Lindblad \cite{CL} to study such dichotomy.


\end{abstract}
\section{Introduction and the main result}
\quad We consider the barotropic compressible
Navier-Stokes equations in a simply connected bounded $C^{2,1}$ domain $\Omega\subset \r^3$,
\begin{equation}\label{11}
\begin{cases}
\rho_t+\mathrm{div}(\rho u)=0, \\
(\rho u)_t+\mathrm{div}(\rho u\otimes u)-\mu\Delta u-\lambda\nabla\mathrm{div} u+\nabla P=0,
\end{cases}
\end{equation}
where $\rho=\rho(x, t)$ and $u=(u_1(x, t), u_2(x, t), u_3(x, t))$ denote the density and the velocity field respectively. The pressure and the viscosity coefficients are given by 
\begin{equation}\label{CA101}
P(\rho)=\rho^\gamma,~~\mu>0,~~\lambda\geq\mu/3,
\end{equation}
for some constant $\gamma\geq 1$.
Moreover, the system is imposed the non-slip   boundary conditions:
\begin{equation}
u\big|_{\partial\Omega}=0,~~ \forall t\geq 0, \label{15}
\end{equation}
and the initial values are provided by 
\begin{equation}\label{CC102}
\rho(x, 0)=\rho_0(x),~~ u(x, 0)= u_0(x),~~ x\in\Omega.
\end{equation}
The main topic of this paper is the initial boundary value problem \eqref{11}--\eqref{CC102}.

There is huge number of literature about the well-posedness for the
multidimensional compressible Navier-Stokes systems with constant viscosity coefficients.
The history of the area may trace back to Nash \cite{Na} and Serrin \cite{se1},
who established the local existence and uniqueness  of classical solutions away from vacuum respectively.
An intensive treatment of the compressible Navier-Stokes equations started with the pioneering papers
by Kanel \cite{KAN}, Itaya \cite{nit1}, and Kazhikhov-Solonnikov \cite{Kaz2}. Matsumura-Nishida \cite{M1} first obtained the global existence of the 3D global smooth solutions, with properly small (in $H^3$-norm) initial values. Then 
Hoff \cite{H3, Hof2} studied the problem with the discontinuous initial data, by introducing a new type of a priori estimates based on the material derivatives $\dot{u}$.
The major breakthrough in the framework of weak solutions is due to Lions \cite{L1}, where he obtained the global existence of weak solutions,  merely under the condition of finite initial energy. For technical reasons, the method in \cite{L1} requires $\gamma\geq 9/5$ in dimension 3, while such index is further released to the critical case $\gamma>3/2$ by Feireisl-Novotn\'{y}-Petzeltov\'{a} \cite{fe1}. 
Recently, Huang-Li-Xin \cite{hlx1} and Li-Xin \cite{jx01} established the global well-posedness of classical solutions with the small initial energy, where large oscillations and vacuum are available. Bresch-Jabin \cite{BJ} obtained the weak solutions to the compressible system with the anisotropic viscous stress tensor and non-monotonic pressure law. Merle-Rapha\"{e}l-Rodnianski-Szeftel \cite{MRRS} studied the implosion of the compressible flow with smooth initial values and investigated the formation of  singularity in the finite time.

The works mentioned above mainly concern the global case $\mathbb{R}^n$, besides, an interesting problem is to consider the compressible Navier-Stokes equations in the bounded domains with proper boundary conditions. The study on the global smooth solutions subject to the Navier-slip boundary conditions is rather completed up to now, see \cite{Ho3, caili01, FLL} for examples, while the progress in Dirichlet  problems which provide another typical boundary condition is slow.
Salvi-Stra\v{s}kraba \cite{S2} obtained the global existence of strong solutions, provided that $\|u_0\|_{H^2}+\|\nabla\rho_0\|_{L^q}$ is initially small for some $q\in(3,6]$.
Matsumura-Nishida \cite{MN2} constructed the global classical solutions in the upper half plane and exterior domains, by assuming the initial values are close (in $H^3$-norm) to the non-vacuum equilibrium. Recently, under some stringent restriction upon the viscosity coefficients
\begin{equation}\label{116}
\mu/(\mu+\lambda)\leq\varepsilon_0,
\end{equation}
with $\varepsilon_0$  properly small, Perepelitsa \cite{PP1} established the global existence of the   non-vacuum weak solutions in the Hoff's regularity class. Huang-Su-Yan-Yu \cite{HY} obtained the symmetric smooth solutions with large initial values on solid balls to the Vaigant-Kazhikhov model.

In summary, for $\mu$ and $\lambda$ satisfying the physical restriction \eqref{CA101}, whether there exists a global smooth solution with initial vacuum to the system \eqref{11}--\eqref{CC102} remains open. This article will adopt some ideas in free boundary problems to give a positive answer. 

To begin with, we introduce the notations through out the paper: For $p\in[1,\infty]$, $k\in\mathbb{N}^+$, and $\alpha\in(0,1)$, the standard Lebesgue, Sobolev, and H\"{o}lder spaces on the domain $A$ are given by 
$L^p(A)$, $W^{k,p}(A)$, and $C^\alpha(A)$ respectively, in particular, $H^k(A)\triangleq W^{k,2}(A)$. 
The material derivative and the integral average of a function $f$ are defined by
$$\dot{f}\triangleq \partial_tf+u\cdot\nabla f,~~\bar{f}\triangleq \frac{1}{|\Omega|}\int_\Omega f\,dx.$$

Then, the regularity assumption on the domain $\Omega$ ensures that
we can select a constant $r_0$, such that the level set
\begin{equation}\label{CC104}
F_r\triangleq\{x\in\Omega|\,\mathrm{dist}(x,\partial\Omega)=r\}
\end{equation}
gives a family of 2 dimensional $C^{2,1}$ surfaces with $r\in[0,r_0]$. The closed domains lying between two surfaces $F_{i}$ and $F_{j}$ are denoted by
\begin{equation}\label{CC105}
\Omega^i_j\triangleq\big\{x\in\Omega\big|i \leq\mathrm{dist}(x,\partial\Omega)\leq j \big\},~~\Omega_{k}\triangleq\{x\in\Omega\big|0 \leq\mathrm{dist}(x,\partial\Omega)\leq k\big\},
\end{equation}
thus $\Omega_{k}$ is a domain near the boundary $\partial\Omega$. For convenience, let us assume  $r_0=10$, and the inner open domain $\Omega\setminus\Omega_{10}$ is not empty. 

We suppose that $\alpha<1/3$ and the initial values satisfy the regularity conditions:
\begin{equation}\label{CC103}
\begin{split}
&\qquad\quad 0\leq\rho_0\leq M_0~\mathrm{in}\ \Omega,~~ \rho_0\geq\underline{\rho}~\mathrm{in}\ \Omega_{10},\\
&\|\rho_0\|_{C^\alpha(\Omega_{10})}+\|\rho_0\|_{W^{1,4}(\Omega_{8})} +\|\rho_0\|_{H^2(\Omega_{6})}\leq M_0,\\
&\qquad\qquad\| u_0\|_{H^1(\Omega)}+\|u_0\|_{H^2(\Omega_{8})} 
\leq M_0,\\
\end{split}
\end{equation}
for some positive constants $M_0$ (not necessarily small) and $\underline\rho$ satisfying $2\underline{\rho}/3\leq\bar{\rho}_0\leq 3M_0/2$, which means $\rho_0$ is out of vacuum and gains higher regularities \textbf{near} the boundary. Moreover, the initial energy is defined by
\begin{equation}\label{C108}
\begin{split}
C_0
&\triangleq\int_{\Omega}\left(\frac{1}{2}\n_0|u_0|^2 +G(\rho_0)\right)dx,~~G(\rho)\triangleq \rho \int_{\bar{\rho}}^{\rho}\frac{P(s)-P(\bar{\rho})}{s^{2}} ds.\\
\end{split}
\end{equation}
\begin{remark}
By integrating \eqref{11}$_1$ over $\Omega\times[0,t]$, we apply \eqref{CC102} to infer that
\begin{equation*}
\int_\Omega\rho\,dx=\int_\Omega\rho_0\,dx.
\end{equation*}
Consequently, we will not distinguish $\bar{\rho}$ and $\bar{\rho}_0$ in the rest of the paper.
\end{remark}
In terms of these notations, we can state our main results as follows:
\begin{theorem}\label{T1}
For some $3<q<6$, suppose that  
\begin{equation}\label{1114}
\rho_0\in W^{2,q}(\Omega),\ u_0\in H^2(\Omega)\cap H^1_0(\Omega).
\end{equation}
If the initial value $(\rho_0,u_0)$ satisfies the regularity conditions \eqref{CC103} and the following compatibility condition holding for some $ g\in L^2(\Omega),$
\begin{equation}\label{1116}
-\mu\Delta u_0-(\mu+\lambda)\nabla\div u_0+\nabla P(\rho_0)=\rho_0^{1/2}g.
\end{equation}
Then there is a positive constant $\ve$ depending only on $\mu,\,\lambda,\,\ga,\,\underline{\rho},\,M_0,$ and $\Omega$, such that the initial boundary value problem \eqref{11}--\eqref{CC102} admits a unique classical solution $(\n,u)$
in $\Omega\times (0,\infty)$, satisfying for any $0<\tau<T<\infty$,
\begin{equation}\label{1118}
\begin{cases}
\rho\in C\big([0,T];W^{2,q}(\Omega)\big),\\
u\in C\big([0,T];H^2(\Omega)\big)\cap L^\infty\big(\tau,T; W^{3,q}(\Omega)\big),\\
u_t\in L^{\infty}\big(\tau,T; H^2(\Omega)\big)\cap H^1\big(\tau,T; H^1(\Omega)\big),\\
\sqrt{\n}u_t\in L^\infty\big(0,\infty;L^2(\Omega)\big),
\end{cases}
\end{equation}
provided that
\be \la{lj03} C_0\le\varepsilon.\ee  

Moreover,
$\forall\,(r,p)\in[1,\infty)\times[1,6)$, there are constants $\mathbf{C}$ and $\eta$ determined by $\mu,\,\lambda,\,\ga,\,\underline{\rho},\,M_0,\,\Omega,\,r,$ and $p$, such that $\forall t\geq 1$,
\begin{equation}\label{1151}
\|\rho(\cdot, t)-\bar{\rho}\|_{L^r(\Omega)}+\|u(\cdot, t)\|_{W^{1,p}(\Omega)}\leq \mathbf{C}\,e^{-\eta t}.
\end{equation}
\end{theorem}
As a direct consequence of Theorem \ref{T1} and \cite[Theorem 1.2]{lx}, we declare the blow-up behaviour of $\nabla\rho$ with an exponential rate when vacuum appears initially.
\begin{theorem}\label{T2}
Suppose that $(\n,u)$ is the global solution obtained in Theorem \ref{T1}. If $\n_0(x_0)=0,$  for some $x_0\in\Omega$, then $\forall r>3,$ there exist positive constants $\hat{\mathbf{C}}$ and $\hat{\eta}$ determined by $\mu,$ $\lambda,$ $\ga,$ $\underline{\rho},$ $M_0,$ $\Omega,$
and $r$, such that $\forall t\geq 1$,
\begin{equation*}
\|\nabla\rho (\cdot,t)\|_{L^r(\Omega)}\geq \hat{\mathbf{C}}\, e^{\hat{\eta} t}.
\end{equation*}
\end{theorem}

A few remarks are in order:

\begin{remark}
In view of the standard embedding
\begin{equation*}
\begin{split}
&L^\infty\big(\tau,T;H^1(\Omega)\big)\cap H^1\big(\tau,T;H^{-1}(\Omega)\big)\hookrightarrow C\big([\tau,T];L^q(\Omega)\big)~~\forall q\in[2,6),\\
\end{split}
\end{equation*}
we apply the Sobolev's inequality together with   \eqref{11}$_1$ and \eqref{1118} to deduce
\begin{equation*}
\begin{split}
&\rho,\nabla\rho \in C(\overline{\Omega}\times[0,T]),\\
\rho_t,u,&\na u,\na^2 u, u_t \in C(\overline{\Omega}\times [\tau,T]).
\end{split}
\end{equation*}
which means the solution obtained in Theorem \ref{T1} is a classical one for any positive time. Moreover, the initial condition \eqref{CC103}$_1$ only excludes vacuum in $\Omega_{10}$, thus this classical solution possibly contains vacuum in the inner domain $\Omega\setminus\Omega_{10}$.
\end{remark}

\begin{remark} 
Theorem \ref{T1} holds for all viscosity coefficients satisfying the physical conditions \eqref{CA101}, and the initial condition \eqref{lj03} only requires the smallness of $\|(\rho_{0}-\bar\rho_0,\sqrt{\rho_0}u_0)\|_{L^2(\Omega)}$. This result improves the theory due to Matsumura-Nishida \cite{MN2} and  Perepelitsa \cite{PP1} by dropping the additional assumption \eqref{116} on viscosity coefficients, and allowing the density with large oscillations and  vacuum.
\end{remark}

\begin{remark} 
Theorem \ref{T2} indicates that $\nabla\rho$ will grow unbounded in an exponential rate, which is in sharp contrast to the results of Matsumura-Nishida \cite{MN2}, where $\nabla\rho$ remains uniformly small.
Such phenomenon seems to occur only in bounded domains, and there is no analogous result for Cauchy problems up to now.
\end{remark}

We now make some comments on the analysis of the paper. In the case of Cauchy problems \cite{hlx1} and the Navier-slip boundary conditions \cite{caili01,Ho3}, the global well-posedness assertions rely heavily on the $L^1(0,T;L^\infty(\Omega))$ estimates of $F\triangleq(\mu+\lambda)\mathrm{div}u-(\rho^\gamma-\bar{\rho}^\gamma)$, which are obtained by solving a well-posed Poisson equation,
\begin{equation}\label{CC106}
\begin{cases}
\Delta F=\mathrm{div}(\rho\dot{u})~~&\mathrm{in}\ \Omega,\\
n\cdot\nabla F=\rho\dot{u}\cdot n+\nabla\times\mathrm{rot}u\cdot n~~&\mathrm{on}\ \partial\Omega.
\end{cases}
\end{equation}
However, the non-slip boundary conditions \eqref{15} fail to give proper boundary values of \eqref{CC106}, thus some new observations are in demand.
The essential ingredient of this paper is  the \textbf{local} $H^2$ energy estimates of $\rho$ and there are several key points.

\textit{1. The partition of the domain $\Omega$.}

We will divide $\Omega$ into the boundary domain $B$ and the inner domain $I=\Omega\setminus B$, such that vacuum only appears in $I$. Therefore the boarder line between these domains should be a free surface, in order to avoid the flow lines originated from the vacuum in $I$ touching $B$. Particularly, we will study this problem in the Lagrangian coordinates and geometric settings introduced by Christodoulou-Lindblad \cite{CL}, in which the moving domain $B$ and corresponding truncation functions all turn into fixed ones. We mention that the notations in this framework are intrinsically defined, which ensure us freely changing between the Lagrangian and Eulerian pictures and carrying out the calculations in the proper coordinates (see Remark \ref{R24} for detailed discussions).

\textit{2. The structure of the solution.}

The local energy estimates require integrations by parts over $B$, which also lead to some tough boundary terms on the interior free surface. To get over them, we will set up an intermediate layer $\Gamma$ strictly containing this free surface, in view of the trace theorem \cite[Theorem 1, Sec 5.5]{ELC}, the estimates in $\Gamma$ also yield proper controls on the free surface. We mention that $\Gamma$ is away from $\partial\Omega$ and intersects both $B$ and $I$. The position of $\Gamma$ can be described as $B\rightarrow \Gamma\rightarrow I$.

Inspired by \cite{PP1}, the key step is to  carry out the point-wise $C^\alpha$ estimates in $\Gamma$, which avoid integrations by parts (see Lemma \ref{lm51}) and drop the additional restriction on $\mu$ and $\lambda$ in \cite{PP1}.    
More inspirationally, note that the regularity of $\rho$ in the boundary domain $( H^2(B))$ is much higher than that in the inner domain $(  L^\infty(I))$, while $C^\alpha$ estimates exactly serve as a transition between $H^2$ and $L^\infty$. In fact, the regularity of $\rho$ is given by
\begin{equation*}
H^{2}(B)\rightarrow C^{\alpha}(\Gamma)\rightarrow L^{\infty}(I).
\end{equation*}
Totally speaking, the regularity of $\rho$  gradually decreases as the distance to $\partial\Omega$ increases, which is the main structural feature of the solution.  

\textit{3. The advantage of the interior estimates.}

Different from the classical free boundary problems, no boundary condition is imposed on the inner free surface. Our method is based on the interior elliptic estimates rather than the arguments depending on the regularity of $\partial\Gamma$ (like \cite{CL}).  
In fact, $\|\nabla^2u\|_{L^\infty(\partial\Gamma)}$ controls the second fundamental form on $\partial\Gamma$, however it possibly grows unbounded due to the blow up assertion in Theorem \ref{T2}. Thus a cusp on $\partial\Gamma$ may occur as $t\rightarrow\infty$, and the arguments in \cite{CL} are violated. In contrast, the interior estimates merely depend on the thickness of $\Gamma$ and $\mathrm{dist}(\Gamma,\partial\Omega)$ which are completely described by $\|u\|_{L^\infty(\partial\Gamma)}$, therefore the exponential decay of $\|\rho-\bar{\rho}\|_{L^2(\Omega)}$ along with the small initial energy yields uniform controls on these terms. Such advantage of interior estimates turns out to be crucial in the local energy estimates.

The rest of paper is organized as follows:
Some elementary properties and geometric tools are collected in Section \ref{S2}. Then we derive the basic energy estimates in Section \ref{sec5} and carry out the near boundary estimates in Section \ref{S4}. After that, the standard arguments provide higher order a priori estimates in Section \ref{S5} and finish the main theorems in Section \ref{S6}.

\section{Preliminaries}\label{S2}
\quad This section provides necessary tools for later arguments. To make notations and ideas clear, we always assume $\gamma=\mu=\lambda=1$ in the rest of the paper, and the methods given below are valid for general settings.  

\subsection{Some elementary properties}
\quad First of all, we quote the local well-posedness theory in\cite{hxd1},
which allows the initial vacuum.
\begin{lemma} \label{u21}
Under the conditions of Theorem \ref{T1},
there is a small time $T>0$ depending only on $\Omega,$ $g,$ $\|\rho_0\|_{W^{2,q}(\Omega)},$ and $\|u_0\|_{H^2(\Omega)}$, such that there exists a unique strong solution $(\rho, u)$ to the problem \eqref{11}-\eqref{CC102} on $\Omega\times(0, T]$
satisfying for any $\tau\in(0,T)$,
\begin{equation*} 
\begin{cases}
\rho\in C\big([0,T]; W^{2,q}(\Omega)\big), \\
u\in C\big([0,T];H^2(\Omega)\big)\cap L^\infty\big(\tau,T;W^{2,q}(\Omega)\big),\\   
u_t\in L^\infty\big(\tau,T;H^2(\Omega)\big),\ u_{tt}\in L^2\big(\tau,T; H^1(\Omega)\big),\\
\sqrt{\rho}u_t\in L^\infty\big(0,T;L^2(\Omega)\big).
\end{cases}
\end{equation*}
\end{lemma}

Next, we provide some useful div-rot type estimates which can be found in \cite{vww,CANEHS,caili01,Swz,bkm}.
\begin{lemma}
Let $(p,q)\in(1,\infty)\times(3,\infty)$ and $\Omega$ be a connected bounded $C^{2,1}$ domain in $\mathbb{R}^3$. For $v\in W_0^{1,p}(\Omega)$ and $u\in W_0^{1,p}(\Omega)\cap W^{2,q}(\Omega),$
there is a constant $C$ depending only on $p$, $q$, and $\Omega$, such that
\begin{equation}\label{24}
\begin{split}
&\|\nabla v\|_{L^p(\Omega)}\leq C\bigg(\|\mathrm{div}v\|_{L^p(\Omega)}+\|\mathrm{rot}v\|_{L^p(\Omega)}\bigg),\\
\|\nabla u\|_{L^\infty(\Omega)}\leq C\bigg(\|\mathrm{div}u&\|_{L^\infty(\Omega)} +\|\mathrm{rot}u\|_{L^\infty(\Omega)}\bigg)\,\log(e+\|\nabla^2u\|_{L^q(\Omega)})+C\|\nabla u\|_{L^2(\Omega)}+C. 
\end{split}
\end{equation}
\end{lemma}

The well-known Gagliardo-Nirenberg's inequality (see \cite{nir}) will be used frequently. 
\begin{lemma}
Suppose that $\Omega$ is a bounded Lipschitz domain in $\r^3$, and $u$ is a $C^\infty$ function on $\Omega$. Let $p,q\in[1,\infty]$, then following assertions are valid.

1). There are two constants $C_1$ and $C_2$ depending only on $p$, $q$, and $\Omega$, such that
\begin{equation}\label{GN}
\begin{split}
\|u\|_{L^r(\Omega)}&\leq C_1\,\|u\|_{L^p(\Omega)}^\theta\cdot\|\nabla u\|_{L^q(\Omega)}^{1-\theta}+C_2\,\|u\|_{L^2(\Omega)},\\
\|\nabla u\|_{L^{r'}(\Omega)}&\leq C_1\,\|u\|_{L^p(\Omega)}^{\theta'}\cdot\|\nabla^2 u\|_{L^q(\Omega)}^{1-\theta'}+C_2\,\|u\|_{L^2(\Omega)},\\
\end{split}
\end{equation}
where $1/r=\theta/p+(1-\theta)\cdot(1/q-1/3)$, and $1/r'=1/3+\theta'/p+(1-\theta')\cdot(1/q-2/3)$ with $\theta'\leq 1/2$. In particular, $C_2=0$ if $u|_{\partial\Omega}=0$ or $\bar{u}=0$.

2). Suppose that $\Gamma'\subset\Gamma\subset\Omega$ with $\mathrm{dist}(\Gamma',\partial\Gamma)>d$, then there is a constant $C$ depending on $d$, $p$, $q$, and $\Omega$, such that
\begin{equation}\label{GN2}
\|u\|_{L^r(\Gamma')}\leq C\big(\|u\|_{L^p(\Gamma)}^\theta\|\nabla u\|_{L^q(\Gamma)}^{1-\theta}+\|u\|_{L^2(\Gamma)}\big),
\end{equation}
where $1/r=\theta/p+(1-\theta)\cdot(1/q-1/3)$.
\end{lemma}

In view of \cite[Theorem III.3.1]{GPG}, we can introduce an ``inversion" of $\mathrm{div}(\cdot)$.
 
\begin{lemma}\la{lemmaz01}
Let $\hat{L}^p(\Omega)\triangleq \{f\in L^p(\O)|\,\bar f =0\}$, then there is a linear operator $\mathcal{B}$ well defined on $\hat{L}^p(\Omega)$  for all $p\in(1,\infty)$, such that the mapping
\begin{equation*}
\mathcal{B}: \hat{L}^p(\Omega)\rightarrow\big(W^{1,p}_0(\O)\big)^3,
\end{equation*} 
satisfies $\mathrm{div}\mathcal{B}(f)=f$, and there is a constant $C$ determined by $p$ and $\Omega$, 
\begin{equation}\label{CC220}
\|\mathcal{B}(f)\|_{W^{1,p}_0(\Omega)}\le C \|f\|_{L^p(\Omega)}.  
\end{equation}
Moreover, if $f=\mathrm{div}g$ with $g\cdot n|_{\partial\Omega}=0$, $\mathcal{B}(\cdot)$ is also well defined on $f$ and admits
\begin{equation*} 
\|\mathcal{B}(f)\|_{L^{p}(\Omega)}\leq C\|g\|_{L^p(\Omega)},~~\forall p\in(1,\infty). 
\end{equation*}
\end{lemma}

Then, we adopt the ideas in \cite{PP1} and assume that $w$ and $v$ are determined by
\begin{equation}\label{CC201}
\begin{split}
&\begin{cases}
\Delta w+\nabla\mathrm{div}w =\rho\dot{u}  &\mathrm{in}\ \Omega, \\
w=0 &\mathrm{on}\ \partial\Omega;
\end{cases}
\\ 
&\begin{cases} 
\Delta v+\nabla\mathrm{div}v=\nabla(\rho-\bar{\rho}) &\mathrm{in}\ \Omega, \\
v=0 &\mathrm{on}\ \partial\Omega.
\end{cases}
\end{split}
\end{equation}
The uniqueness assertion of the Lam\'{e} system under Dirichlet boundary conditions (See \cite[Theorem 10.1-10.2]{adn}) ensures the decomposition 
\begin{equation}\label{CC214}
u= w+v.
\end{equation}
The standard $C^\alpha$ and $L^p$ elliptic estimates (see  \cite[Theorems 10.1-10.2]{adn} and \cite[Theorems 9.11-9.13]{GT}) and the trace theorem  (\cite[Theorem 1, Sec 5.5]{ELC}) are available for \eqref{CC201}.
\begin{lemma}
Suppose that $w$ and $v$ are given by \eqref{CC201}, and $k\in\mathbb{N},p\in(1,\infty),\alpha\in(0,1)$, then we deduce the following global and local estimates.

1).$($Global estimates$)$ There is a constant $C$ depending on $k$, $p$, $\alpha$, and $\Omega$, such that
\begin{equation}\label{355}
\begin{split}
&\|w\|_{W^{k+2,p}(\Omega)}\leq C\|\rho\dot{u}\|_{W^{k,p}(\Omega)},\\
&\|v\|_{W^{k+1,p}(\Omega)}\leq C\|\rho-\bar{\rho}\|_{W^{k,p}(\Omega)},\\
&\|v\|_{C^{k+1,\alpha}(\Omega)}\leq C\|\rho-\bar{\rho}\|_{C^{k,\alpha}(\Omega)}.
\end{split}   
\end{equation}

2).$($Local estimates$)$ Suppose that $\Gamma'\subset\Gamma\subset\Omega$ with $\mathrm{dist}(\Gamma',\partial\Gamma)>d$, then there is a constant $C$ depending on $d$, $k$, $p$, $\alpha$, and $\Omega$, such that
\begin{equation}\label{C406}
\begin{split}
&\|w\|_{W^{k+2,p}(\Gamma')}\leq C\big(\|\rho\dot{u}\|_{W^{k,p}(\Gamma)}+\|w\|_{L^2(\Omega)}\big),\\
&\|v\|_{W^{k+1,p}(\Gamma')}\leq C\big(\|\rho-\bar{\rho}\|_{W^{k,p}(\Gamma)}+\|v\|_{L^2(\Omega)}\big),\\
&\|v\|_{C^{k+1,\alpha}(\Gamma')}\leq C\big(\|\rho-\bar{\rho}\|_{C^{k,\alpha}(\Gamma)}+\|v\|_{L^2(\Omega)}\big).
\end{split}
\end{equation}
In particular, $F_v\triangleq 2\,\mathrm{div}v-(\rho-\bar{\rho})$ and $\mathrm{rot}v$ satisfy that $\Delta F_v=\Delta\mathrm{rot}v=0$ in $\Omega$, therefore the interior estimates for harmonic functions yield that
\begin{equation}\label{CC215}
\|\nabla^k \mathrm{rot}v\|_{L^\infty(\Gamma')}+\|\nabla^k F_v\|_{L^\infty(\Gamma')}\leq C\|F_v\|_{L^2(\Gamma)}\leq C\|\rho-\bar{\rho}\|_{L^2(\Omega)}.
\end{equation}

3).$($Boundary estimates$)$ Suppose that $\partial\Omega\subset\Gamma'\subset\Gamma\subset\Omega$ with $\mathrm{dist}(\Gamma',\partial\Gamma\setminus\partial\Omega)>d$, then the estimates \eqref{C406} are valid as well. 
Moreover, for any $f\in H^{k+1}(\Gamma')$, there is a constant $C$ depending only on $d$, $k$, and $\Omega$, such that
\begin{equation}\label{C409}
\int_{\partial\Omega}|\nabla^k f|^2\,dS\leq C\,\|\nabla^{k} f\|_{H^1(\Gamma)}\cdot\|\nabla^k f\|_{L^2(\Gamma)}.
\end{equation}
\end{lemma}

The local estimates \eqref{C406}$_1$ can be slightly improved if we take $p=2$. Let us illustrate the following weighted local elliptic estimates.  
\begin{lemma}\label{lm42}
Suppose that 
$\partial\Omega\subset\Gamma'
\subset\Gamma\subset\Omega$ with $\mathrm{dist}(\partial\Gamma'\setminus\partial\Omega,\,\partial\Omega)>d$, $\mathrm{dist}(\partial\Gamma\setminus\partial\Omega,\,\Gamma')>d$. Let $\phi(x)\in C^\infty(\Omega)$ be a truncation function satisfying
\begin{equation*}
\begin{cases}
\phi=1~~\forall x\in\Gamma',\\
\phi=0~~\forall x\in\Omega\setminus\Gamma,\\
0\leq\phi\leq 1,\,|\nabla\phi|\leq C/d~~\forall x\in\Omega.
\end{cases}
\end{equation*} 
Then there is a constant $C$ depending on $\Omega$ and $d$, such that
\begin{equation}\label{C120}
\int_\Omega|\nabla^2 w|^2\cdot\phi^2\,dx\leq C\int_\Omega|\rho\dot{u}|^2\cdot\phi^2\,dx+C\,\|\nabla w\|_{L^2(\Omega)}^2.
\end{equation}
\end{lemma}
\begin{proof}
Let us first consider another truncation function $\zeta$ satisfying
\begin{equation}\label{CC701}
\begin{cases}
\zeta=0~~\mathrm{if}\ \mathrm{dist}(x,\partial\Omega)<d/4,\\
\zeta=1~~\mathrm{if}\ \mathrm{dist}(x,\partial\Omega)\geq d/2,\\
0\leq\zeta\leq 1,\,|\nabla\zeta|\leq C/d~~\forall x\in\Omega,
\end{cases}
\end{equation}
then we also define $\chi\triangleq\phi\cdot\zeta.$
Now multiplying \eqref{CC201}$_1$ by $\Delta w\cdot\chi^2$ and integrating over $\Omega$ implies 
\begin{equation}\label{C122}
\begin{split}
\int_\Omega\Delta w\cdot\Delta w\cdot\chi^2\,dx
+\int_\Omega\nabla\mathrm{div}w\cdot\Delta w\cdot\chi^2\,dx
&=\int_\Omega\rho\dot{u}\cdot\Delta w\cdot\chi^2\, dx.
\end{split}
\end{equation}

Observe that $\xi\in C_0^\infty(\Omega)$, thus no boundary term occurs during the integration by parts and we argue that
\begin{equation*}
\begin{split}
&\int_\Omega\Delta w\cdot\Delta w\cdot\chi^2\,dx\\
&=-\int_\Omega\nabla w\cdot\Delta(\nabla w)\,\chi^2 \,dx+\int_\Omega\nabla w\cdot\nabla^2 w\cdot\nabla\chi\cdot\chi\,dx\\
&=\int_\Omega|\nabla^2w|^2\chi^2\,dx
+\int_\Omega\nabla w\cdot\nabla^2 w\cdot\nabla\chi\cdot\chi\,dx\\
&\geq\frac{1}{2}\int_\Omega|\nabla^2w|^2\chi^2\,dx-C\,\int_\Omega|\nabla w|^2\,dx,
\end{split}
\end{equation*}
where we have taken advantage of the fact $|\nabla\chi|\leq C/d$ in the last line.
Similar arguments also lead to the fact that
\begin{equation*}
\begin{split}
&\int_\Omega\nabla\mathrm{div}w\cdot\Delta w\cdot\chi^2\,dx\\
&\geq\int_\Omega|\nabla\mathrm{div}w|^2\cdot\chi^2\,dx
-\frac{1}{10}\int_\Omega|\nabla^2w|^2\cdot\chi^2 \,dx-C\int_\Omega|\nabla w|^2\,dx.
\end{split}
\end{equation*}
Meanwhile, we directly declare that
\begin{equation*} 
\int_\Omega\rho\dot{u}\cdot\Delta w\cdot\chi^2 \,dx\leq C\int_\Omega|\rho\dot{u}|^2\cdot\chi^2\,dx
+\frac{1}{10}\int_\Omega|\nabla^2w|^2\cdot\chi^2\,dx.
\end{equation*}
Thus, substituting these estimates into \eqref{C122} gives that
\begin{equation}\label{C131}
\int_\Omega|\nabla^2 w|^2\cdot\chi^2\,dx\leq C\int_\Omega|\rho\dot{u}|^2\cdot\phi^2\,dx+C\|\nabla w\|_{L^2(\Omega)}^2.
\end{equation}

At last, we mention that 
\begin{equation*} 
\begin{split}
\int_\Omega|\nabla^2w|^2\cdot\phi^2\,dx
&=\int_\Omega|\nabla^2w|^2\cdot\chi^2\,dx+\int_\Omega|\nabla^2w|^2\cdot\phi^2(1-\zeta)^2\,dx.\\
\end{split}
\end{equation*}
According to the definition of $\phi$ and $\zeta$, we argue that
\begin{equation*}
\begin{split}
&\mathrm{dist}\big(\mathrm{supp}(\phi\cdot(1-\zeta)),\,\partial\Gamma'\setminus\partial\Omega\big)
>d/8,
\end{split}
\end{equation*}
which along with \eqref{C406}$_1$ applied to
$\mathrm{supp}\big(\phi\cdot(1-\zeta)\big) \subset \Gamma'$ gives
\begin{equation}\label{C130}
\int_\Omega|\nabla^2w|^2\cdot\phi^2\big(1-\zeta\big)^2\,dx\leq
C\int_\Omega|\rho \dot{u}|^2\cdot\phi^2\,dx+C\|\nabla w\|_{L^2(\Omega)}^2.
\end{equation}
Consequently, combining \eqref{C131}--\eqref{C130} yields \eqref{CC701} and the proof is therefore completed.
\end{proof}

\begin{remark} 
Note that no structural assumption on the weight $\phi$ is required to obtain \eqref{C120}. In particular, Lemma \ref{lm42} can be regarded to provide a bounded operator $\nabla^2\Delta^{-1}$ from $L^2(\Omega,\,\phi^2\,dx)$ to itself up to some lower order terms.
\end{remark}
\begin{remark}
We specially mention that the domain $\Gamma'$ is ``strictly" contained in $\Gamma$ and the interior estimates depend on the distance between $\Gamma'$ and $\Gamma$, rather than the regularity of $\partial\Gamma'$ and $\partial\Gamma$. Such fact plays important roles in our arguments below.
\end{remark}

Generally speaking, the regular part $w$ is rather easy to handle, since $\rho\dot{u}$ is properly controlled in Section \ref{sec5},
while the analysis on the singular part $v$ relies on some detailed local $C^\alpha$ estimates which is our next topic.
  
Suppose that $\Gamma'\subset\Gamma$, then we introduce the following modified H\"{o}lder type norms:
\begin{equation}\label{D51}
\|\rho\|_{C^\alpha(\Gamma,\Gamma')}
\triangleq\sup_{x\in\Gamma,\ y\in\overline{\Gamma'}}
\frac{|\rho(x)-\rho(y)|}{|x-y|^\alpha}+\sup_{x\in \Gamma}|\rho(x)|.
\end{equation}
We mention that
$C^\alpha(\Gamma)\subset C^\alpha(\Gamma,\Gamma')\subset C^\alpha(\Gamma').$
In particular, $\rho\in C^\alpha(\Gamma,\Gamma')$ is a usual $C^\alpha$ function on $\overline{\Gamma'}$, but is merely bounded on $\Gamma\setminus\overline{\Gamma'}$.
 
Let us consider the domains $\partial\Omega\subset\Gamma'\subset\Gamma\subset\Omega$ , which satisfy that $\mathrm{dist}(\Gamma',\partial\Gamma\setminus\partial\Omega)>d$ and $\mathrm{dist}\,(\partial\Gamma'\setminus\partial\Omega,\,\partial\Omega)>d$.
Then the key technical point of this paper is that the modified norm $\|\rho\|_{C^\alpha(\Gamma,\Gamma')}$ is sufficient to control $\|\nabla v\|_{L^\infty(\Gamma')}$, rather than the full norm $\|\rho\|_{C^\alpha(\Gamma)}$.

\begin{lemma}\label{le52}
Assume that $v$ is given by \eqref{CC201}$_2$, then for $x\in \Gamma$, $y\in \Gamma'$, and $\varepsilon\in(0,1/1000)$,
there is a constant $C$ depending only on $d$ and $\Omega$ such that 
\begin{equation*}
  \frac{|v(x)-v(y)|}{|x-y|}\leq C\,\varepsilon^{\alpha/2}\|\rho\|_{C^\alpha(\Gamma,\Gamma')}+
C\,\varepsilon^{-2}\|\rho-\bar{\rho}\|_{L^4(\Omega)}.
\end{equation*}
\end{lemma}
The rather long proof of Lemma \ref{le52} is postponed to the Appendix.

\subsection{Lagrangian coordinates and related Geometric issues} 
\quad  Let us  make a brief introduction to the Lagrangian coordinates in higher dimensions, and the framework below is due to Christodoulou-Lindblad \cite{CL}.
Suppose that $u(x,t)$ is a smooth vector field on $\Omega$, then $\forall(x,t)\in\Omega\times[0,\infty)$, there is a corresponding ODE called the flow line:
\begin{equation*} 
\begin{cases}
\frac{d}{d\tau}{X}(\tau;x,t)=u(X(\tau;x,t),\tau)~~\forall\tau\geq 0,\\
X(t;x,t)=x.
\end{cases}
\end{equation*}
Let $y=X(0;x,t)$, then the mapping
$(x,t)\mapsto(y,t)$ gives a smooth change of coordinates. Particularly, $(y,t)$ is called the Lagrangian coordinates, and $(x,t)$ is called the Eulerian coordinates. Following the classical text books \cite{Bo, GTM94}, we shortly discuss the geometric structure of the Lagrangian coordinates.
 
First of all, $\forall t\geq 0$ fixed, the flow map 
\begin{equation*}
X_t:\Omega\rightarrow\Omega,~~ 
y \mapsto x(y,t),
\end{equation*}
gives a diffeomorphism from $\Omega$ to itself.
Pulling back the Euclidean metric (\cite[Section 2]{CL}) in Eulerian coordinates via $X_t^*$ (\cite[Section 2.22]{GTM94}, \cite[Chapter V \S 1]{Bo}), we obtain a metric $g$ in Lagrangian coordinates with components $g_{ij}$ provided by 
\begin{equation*}
g_{ij}\triangleq\sum_{k=1}^3\frac{\partial x^k}{\partial y^i}\cdot\frac{\partial x^k}{\partial y^j},~~\forall i,j=1,2,3.
\end{equation*}
Moreover, the metric $g$ induces a Riemanian connection (\cite[Chapter VIII \S 4]{Bo}) with connection coefficients $\Gamma^{k}_{ij}$ (\cite[Section 2]{CL}) given by   
$$\Gamma_{ij}^k\triangleq\frac{1}{2}\,g^{kl}(\partial_{y^i}g_{lj} +\partial_{y^j}g_{li}-\partial_{y^l}g_{ij})=\frac{\partial y^k}{\partial x^l}\cdot\frac{\partial^2x^l}{\partial y^i\partial y^j},$$
where $g^{ij}$ is the matrix inversion of $g_{ij}$, meanwhile the repeated superscripts and subscripts have been summed over.

In Lagrangian coordinates, a tensor $A$ of type $(s,r)$ $\big($\cite[Section 2]{CL}, \cite[Chapter 33 \S 1]{Pos}$\big)$ can be written as:
\begin{equation*}
 A=A^{i_1i_2\cdots i_s}_{j_1j_2\cdots j_r}\,dy^{j_1}\otimes dy^{j_2}\cdots\otimes dy^{j_r}
 \otimes\frac{\partial}{\partial y^{i_1}}\otimes\frac{\partial}{\partial y^{i_2}}\cdots
 \otimes\frac{\partial}{\partial y^{i_s}}.
\end{equation*}
Generally, the inner product (\cite[Chapter V \S 2]{Bo}, \cite[Section 4.10]{GTM94}) between two tensors $A$, $B$ of the same type is defined by
\begin{equation}\label{C208}
\langle A,B\rangle\triangleq g_{i_1k_1}\cdots g_{i_sk_s} g^{j_1l_1}\cdots g^{j_rl_r}\cdot A^{i_1\cdots i_s}_{j_1\cdots j_r}\cdot B^{k_1\cdots k_s}_{l_1\cdots l_r},
\end{equation}
and $|A|\triangleq\langle A,A\rangle^{1/2}$ is the module of $A$. The inner product induces a bijection between tensors of type $(0,1)$ (called 1-forms) and tensors of type $(1,0)$ (called vector fields),
\begin{equation*}
\begin{split}
dy^i\mapsto g^{ij}\,\frac{\partial}{\partial y^j },~~~\frac{\partial}{\partial y^j}\mapsto g_{ij}\,dy^i,~~~\forall i,j=1,2,3.\\
\end{split}
\end{equation*}
Consequently, the 1-form $\omega=\omega_i\,dy^i$ is mapped to the vector field $\omega^j\,\partial/\partial y^j$  with $\omega^j=g^{ij}\omega_i$ and vice versa.

Suppose $A$ is a tensor of type $(s,r)$ and $X$ is a vector field, then the covariant derivative (\cite[Chapter 35 \S 6]{Pos}, \cite[Section 2]{CL}) of $A$ with respect to $X$ is denoted by $\nabla_X A$. In particular, $\nabla A$ gives a tensor of type $(s,r+1)$ with components  
\begin{equation}\label{C42}
\begin{split}
  A^{i_1i_2\cdots i_s}_{j_1j_2\cdots j_r,k}&=\frac{\partial}{\partial y^k}A^{i_1i_2\cdots i_s}_{j_1j_2\cdots j_r}-
  \Gamma_{j_1k}^l\cdot A^{i_1i_2\cdots i_s}_{lj_2\cdots j_r}-\cdots-
  \Gamma_{j_rk}^l\cdot A^{i_1i_2\cdots i_s}_{j_1j_2\cdots l}\\
  &\qquad+\Gamma_{lk}^{i_1}\cdot A^{li_2\cdots i_s}_{j_1j_2\cdots j_r}+\cdots+
  \Gamma_{lk}^{i_s}\cdot A^{i_1i_2\cdots l}_{j_1j_2\cdots j_r},
\end{split}
\end{equation}
where $A_{,k}$ corresponds to the covariant derivative with respect to $\partial/\partial y^k$.
We mention that, for a vector field $u=u^i\,\partial/\partial y^i$, the covariant derivatives   $u^i_{,j}=\partial_{y^j}u^i+\Gamma^{i}_{kj}\cdot u^k$ are usually \textbf{not} the space derivatives $\partial_{y^j}u^i$.

The higher order covariant derivatives $\nabla^k A=\nabla(\nabla^{k-1}A)$ are defined inductively. Some basic properties in Lagrangian coordinates are listed below.
\begin{lemma}\label{L27}
Let $A$ and $B$ be two tensors $($not necessarily of the same type$)$, and $\nabla$ be the covariant differentiation, then following assertions hold:

1). The order of covariant derivatives with respect to $\partial/\partial y^i$ and $\partial/\partial y^j$ can be changed freely,
\begin{equation}\label{CC206}
A_{,ij}=A_{,ji},~~\forall i,j=1,2,3.
\end{equation}

2). Let $X$ be a vector field, then $\nabla_X$ commutes with the contractions (\cite[Chapter 33 \S 1]{Pos}), and satisfies that
\begin{equation}\label{CC207}
\nabla_X(A\otimes B)=\nabla_XA\otimes B+A\otimes\nabla_XB.
\end{equation}
In particular, $(A\otimes B)_{,i}=A_{,i}\otimes B+A\otimes B_{,i}$.

\end{lemma}
\begin{proof}
The detailed proof of 
\eqref{CC206} and \eqref{CC207} can be found in \cite[Section 2 (2.9)]{CL} and \cite[Chapter 2 \S 1]{Pos} respectively. 
\end{proof}

\begin{remark}\label{R24}
As indicated by Christodoulou-Lindblad \cite{CL}, the covariant derivatives play more significant roles than usual space derivatives in Lagrangian coordinates, since the definition of the former  is independent of coordinates, which does great benefits to our  calculations.

We provide a typical example. Let $u=\tilde{u}^i\,\partial/\partial x^i=u^j\,\partial/\partial y^j$ be a vector field on $\Omega$ and $\nabla u$ be the covariant derivatives. In Eulerian coordinates, the covariant derivatives $\tilde{u}_{,j}$ are given by the space derivatives $\partial_{x^j}\tilde{u}$ $($\cite[Chapter VII \S 1]{Bo}, \cite[Section 2]{CL}$)$, thus $\nabla u$ is the usual gradient, and the module of $\nabla u$ is given by 
\begin{equation*} 
|\nabla u|^2=\sum_{i,j}|\partial_{x_i}\tilde{u}^j|^2=\delta^{ij}\delta_{kl}\cdot\partial_{x_i}\tilde{u}^k\,\partial_{x_j}\tilde{u}^l.
\end{equation*}
However if we take the change of coordinates, this term just transforms into
\begin{equation*} 
|\nabla u|^2=g^{ij}g_{kl}\cdot u^k_{,i}\,u^l_{,j},
\end{equation*}
which coincides with the definition $|\nabla u|=\langle\nabla u,\nabla u\rangle^{1/2}$ in Lagrangian coordinates as \eqref{C208}, rather than $\sum_{i,j}|\partial_{y^i}u^j|^2$. Such fact also ensures us \textbf{not} to distinguish the module $|\cdot|$ in different coordinates through out the paper.
\end{remark}

We utilize the idea in \cite{CL} to introduce the tangential metric induced by $g$ and the corresponding projection to the tangential plane.
Different form the free boundary problems, $\Omega$ is fixed due to the non-slip boundary conditions \eqref{15}, thus we start from the Eulerian coordinates $(x,t)$.

Precisely, let $\tilde{\mathcal{N}}=\tilde{N}^i(x)\,\partial/\partial x^i$ be the smooth global extension (given in \cite[(3.10)--(3.11)]{CL}) of the unit outer normal vector on $\partial\Omega$, which satisfies that
\begin{equation}\label{CC221}
\langle\tilde{\mathcal{N}}(x),\tilde{\mathcal{N}}(x)\rangle=\sum_{i}|\tilde{N}^i(x)|^2 =1~~\mathrm{when}\ \mathrm{dist}(x,\partial\Omega)\leq 8.
\end{equation}
Let $\tilde{N}=\tilde{N}_i(x)\,dx^i$ be the corresponding 1-form of $\tilde{\mathcal{N}}$, then the tangential metric $\widetilde{\gamma}(x)$ in Eulerian coordinates is given by
\begin{equation}\label{CC219}
\begin{split}
\widetilde{\gamma}=\,&\delta_{ij}\,dx^i\otimes dx^j-\tilde{N}\otimes\tilde{N},~~\mathrm{with}\ \widetilde\gamma_{ij}=\delta_{ij}-\tilde{N}_i\tilde{N}_j.
\end{split}
\end{equation}
By pulling back, we define $\gamma(y,t)=X_t^*(\widetilde{\gamma})$ as the tangential metric in Lagrangian coordinates. For $N=X_t^*(\widetilde{N})=N_i\,dy^i$, it holds that
\begin{equation}\label{C43}
\begin{split}
&\,\gamma=g-N\otimes N,~~\mathrm{with}\ \gamma_{ij}=g_{ij}-N_iN_j.
\end{split}
\end{equation}  
From another view point, \eqref{C43} is equivalent to
\begin{equation*}
g=\gamma+N\otimes N,~~\mathrm{with}\ g_{ij}=\gamma_{ij}+N_iN_j,
\end{equation*}
which provides a decomposition of the metric  $g$ into the tangential part $\gamma$ and the normal part $N\otimes N$.
Such case is similar with the upper half plane $\mathbb{R}^3_{+}$, where the global coordinates $(x_1$,$x_2$,$x_3)$ divide into the tangential part $(x_1,x_2)$, and the normal part $x_3$.

Then following the Definition 3.1 in \cite{CL}, the tangential projection of a tensor $A$ of type $(s,r)$ is given by a tensor $\Pi A$ of the same type:
\begin{equation*}
(\Pi A)^{i_1i_2\cdots i_s}_{j_1j_2\cdots j_r}
= \gamma^{i_1}_{k_1}\gamma^{i_2}_{k_2}\cdots\gamma^{i_s}_{k_s}\cdot\gamma_{j_1}^{l_1}\gamma_{j_2}^{l_2}
\cdots\gamma_{j_r}^{l_r}\cdot A^{k_1k_2\cdots k_s}_{l_1l_2\cdots l_r},~~\mathrm{with}\ \gamma^{i}_j=g^{ik}\gamma_{kj}.
\end{equation*}
In particular, for a vector field $X=X^i\,\partial/\partial y^i$, $\Pi X\big|_{\partial\Omega}$ is just the projection of $X$ to the tangential plane of $\partial\Omega$ in the usual sense,
\begin{equation*}
 \Pi X =X-\langle X, N\rangle\,N =(\gamma^i_j\cdot X^j)\,\frac{\partial}{\partial y^i}.
\end{equation*}
The tangential projection helps us establishing some useful boundary estimates on $\partial\Omega$.
\begin{lemma}
Suppose that $u$ is a tensor of type $(s,r)$ defined on $\Omega$ and $d\nu$ is the volume element $($Chapter V \S 7 of \cite{Bo}$)$, then the following assertions are valid.  

1) If $u=u^i\partial/\partial y^i$ is a vector field, Stokes' theorem (\cite[Section 4.9]{GTM94}) holds
\begin{equation}\label{CC216}
\int_\Omega \mathrm{div}u\,d\nu=\int_\Omega u^i_{,i}\,d\nu=\int_{\partial\Omega}N_i\,u^i\,d\nu.
\end{equation}
In particular, for some constant $C$ depending only on $\Omega$, we declare that
\begin{equation}\label{CC217}
\begin{split}
\bigg|\int_{\partial\Omega}\gamma^{ij}u_{i,j}\,d\nu\bigg|\leq C\int_\Omega|\nabla N|\cdot|\nabla u|\,d\nu.
\end{split}
\end{equation}

2) If $u\equiv 0$ on $\partial\Omega$, the tangential derivatives of $u$ on $\partial\Omega$ vanish: For any vector field $X$, 
\begin{equation}\label{CC218}
\begin{split}
\nabla_{\Pi X}u=0, ~~\mathrm{on}\ \partial\Omega.
\end{split}
\end{equation}
\end{lemma}
\begin{proof}
1). The detailed proof of Stokes' theorem can be found in \cite[Section 4.9]{GTM94}. To illustrate \eqref{CC217}, we observe that
$$\int_{\partial\Omega}\gamma^{ij}u_{i,j}\,d\nu=\int_{\partial\Omega}(u^i_{,i}-N^iN^j\,u_{i,j})\,d\nu.$$
In addition, \eqref{CC206} along with \eqref{CC216} implies that
\begin{equation*}
\begin{split}
\int_{\partial\Omega}N^iN^j\,u_{i,j}\,d\nu
&=\int_{\Omega}\big(N^j\,u^i_{,j}\big)_{,i}\,d\nu\\
&=\int_{\Omega} N^j_{,i}\,u^i_{,j}\,d\nu+\int_{\Omega} N^j\,u^i_{,ij}\,d\nu\\
&=\int_{\Omega} N^j_{,i}\,u^i_{,j}\,d\nu-\int_{\Omega} N^j_{,j}\,u^i_{,i}\,d\nu+\int_{\Omega} (N^j\,u^i_{,i})_{,j}\,d\nu\\
&=\int_{\Omega} N^j_{,i}\,u^i_{,j}\,d\nu-\int_{\Omega} N^j_{,j}\,u^i_{,i}\,d\nu+\int_{\partial\Omega} u^i_{,i}\,d\nu,\\
\end{split}
\end{equation*}
which yields \eqref{CC217} and finishes the first part.

2). Without loss of generality, let $u$ be a tensor of type $(0,2)$, then $\nabla_{\Pi(\cdot)} u$ gives a tensor of type $(0,3)$, which in view of \eqref{CC219} and \eqref{C43} can be written as 
\begin{equation}\label{CA201}
\nabla_{\Pi(\cdot)} u=\gamma^{s}_k\,u_{ij,s}\,dy^k\otimes dy^i\otimes dy^j=\widetilde{\gamma}^s_k \cdot\frac{\partial \tilde{u}_{ij}}{\partial x^s}\,dx^k\otimes dx^i\otimes dx^j~~\mathrm{with}\ dy^k=\frac{\partial y^k}{\partial x^j}\,dx^j.
\end{equation}
Note that $u\equiv0$ on $\partial\Omega$ implies that the components $\tilde{u}_{ij}\equiv0$ on $\partial\Omega$. Moreover we check that
\begin{equation}\label{CC222} 
\langle\widetilde{\gamma}^s_k\,\frac{\partial}{\partial x^s},\, N^s\frac{\partial}{\partial x^s}\rangle=0 ~~\mathrm{on}\ \partial\Omega,\ \forall k=1,2,3,
\end{equation}
which means that $\widetilde{\gamma}^s_k\,\partial/\partial {x^s}$ are tangential derivatives on $\partial\Omega$ and $\widetilde{\gamma}^s_k\,\partial\tilde{u}_{ij}/\partial {x^s}\equiv0$ on $\partial\Omega$. Thus we must have $\nabla_{\Pi(\cdot)} u\equiv0$ on $\partial\Omega$, which also finishes the proof. 
\end{proof}

Now let us calculate the derivatives of the metrics given above. We mention that the time derivative $D_t f$ in the Lagrangian coordinates corresponds to the material derivative $\dot{f}=(\partial_tf+u\cdot\nabla f)$ in the Eulerian coordinates.  
 
\begin{lemma}\label{411}
Suppose that $g$, $\gamma$, and $N$ are defined above. Let $\Gamma^k_{ij}$ be the connection coefficients and $d\nu$ be the volume element, then the following assertions are valid.

1). The time derivatives in the Lagrangian coordinates are given by
\begin{equation}\label{CC208}
\begin{split}
D_t\,g_{ij}=&u_{i,j}+u_{j,i},~~
D_t\,\Gamma^k_{ij}=u^k_{,ij},~~D_t\,d\nu=u^s_{,s}\,d\nu,\\
&D_t\gamma_{ij}=u^k\,\gamma_{ij,k}+\gamma_{ik}\,u^k_{,j}+\gamma_{jk}\,u^k_{,i}.\\
\end{split}
\end{equation}

2). For any $k\in\mathbb{N}$, there is a constant $C$ depending only on $k$ and $\Omega$, such that the covariant derivatives of these metrics satisfy  
\begin{equation}\label{C401}
\begin{split}
\nabla^k g=0,~~
|\nabla^k\gamma|+|\nabla^k N|\leq C.
\end{split}
\end{equation}
\end{lemma}
\begin{proof}
1). Suppose that $u=\tilde{u}^i\,\partial/\partial x^i=u^j\,\partial/\partial y^j$ with $\tilde{u}^i=u^j\cdot(\partial x^i/\partial y^j)$, then the arguments in \cite[Lemma 2.1]{CL} provide  \eqref{CC208}$_1$ and the fact
\begin{equation*} 
\begin{split}
D_t\frac{\partial x^i}{\partial y^j}=\frac{\partial\tilde{u}^i}{\partial x^k}\cdot\frac{\partial x^k}{\partial y^j}= u^k_{,j}\cdot\frac{\partial x^i}{\partial y^k}.\\
\end{split}
\end{equation*}
Moreover, note that
the definition of $X_t^*$ leads to
\begin{equation}\label{CC209}
\begin{split}
\gamma_{ij}(y,t)=\widetilde{\gamma}_{kl}\big(X_t(y)\big)\,\frac{\partial x^k}{\partial y^i}\cdot\frac{\partial x^l}{\partial y^j}.\\
\end{split}
\end{equation}
Thus let us take $D_t$ on \eqref{CC209} and infer that
\begin{equation*}
\begin{split}
D_t\gamma_{ij}&=(\tilde{u}^s\cdot\frac{\partial\widetilde{\gamma}_{kl}}{\partial x^s})\cdot\frac{\partial x^k}{\partial y^i}\cdot\frac{\partial x^l}{\partial y^j}+\widetilde{\gamma}_{kl}\cdot D_t\frac{\partial x^k}{\partial y^i}\cdot\frac{\partial x^l}{\partial y^j}+\widetilde{\gamma}_{kl}\cdot\frac{\partial x^k}{\partial y^i}\cdot D_t\frac{\partial x^l}{\partial y^j}\\
&=u^t\cdot\frac{\partial\widetilde{\gamma}_{kl}}{\partial x^s}\cdot\frac{\partial x^k}{\partial y^i}\cdot\frac{\partial x^l}{\partial y^j}\cdot\frac{\partial x^s}{\partial y^t}+
\widetilde{\gamma}_{kl}\cdot u^s_{,i}\cdot\frac{\partial x^k}{\partial y^s}\cdot\frac{\partial x^l}{\partial y^j}+
\widetilde{\gamma}_{kl}\cdot u^s_{,j}\cdot\frac{\partial x^l}{\partial y^s}\cdot\frac{\partial x^k}{\partial y^i}\\
&=\gamma_{ij,k}\,u^k+\gamma_{js}\,u^s_{,i}+\gamma_{is}\,u^s_{,j},
\end{split}
\end{equation*}
which gives \eqref{CC208}$_2$. We mention that the last line is due to the change of variables similar with \eqref{CA201}. For example, if we write $\nabla\gamma$ in both coordinates, it holds that
\begin{equation*}
\begin{split}
\nabla\gamma=\frac{\partial\widetilde{\gamma}_{ml}}{\partial x^s}\,dx^s\otimes dx^m\otimes dx^l= \gamma_{ij,k}\, dy^i\otimes dy^j\otimes dy^k,~~\mathrm{with}\ dx^s=\frac{\partial x^s}{\partial y_i}dy^i.
\end{split}
\end{equation*}  
Comparing the components on both sides provides the desired results. 

2). Remark \ref{R24} implies that we can check $|\nabla^k g|$ and $|\nabla^k\gamma|$ directly in the Eulerian coordinates. Observing that
\begin{equation}\label{CC210}
\begin{split}
\nabla^k&(\delta_{ij}\,dx^i\otimes dx^j)=0,~~|\nabla^k\widetilde{\gamma}|\leq C,~~\forall k\in\mathbb{N},
\end{split}
\end{equation}
for some constant $C$ depending only on $k$ and $\Omega$, since $\widetilde{\gamma}$ is fixed. Thus according to Remark \ref{R24}, the estimates \eqref{CC210}  are valid if we replace $(\delta_{ij}\,dx^i\otimes dx^j)$ and $\widetilde{\gamma}$ by $g$ and $\gamma$ respectively, which combining with \eqref{C43} gives \eqref{C401}. The proof is therefore completed.
\end{proof}

To end up this subsection, we carry out some typical calculations in Lagrangian coordinates. Recalling that we have assumed $\mu=\lambda=a=\gamma=1$, therefore it holds in $\Omega$ that
\begin{equation}\label{CC101}
\begin{cases}
D_t\rho =-\rho\,u^k_{,k},\\
\rho\,D_tu_i+\rho\,u^ku_{i,k}-g^{jk}u_{i,jk}-u^j_{,ji}+\rho_{,i}=0.
\end{cases}
\end{equation}

Let us first deal with the conservation of mass. Taking the covariant derivative with respect to $\partial/\partial y^i$ on \eqref{CC101}$_1$ and applying \eqref{CC207} lead to
\begin{equation}\label{C675}
D_t\rho_{,i}=-\rho_{,i}\,u^k_{,k}-\rho\,u^k_{,ki}.
\end{equation}
Note that $D_t(\rho_{,i})=(D_t\rho)_{,i}$, since $\rho_{,i}=\partial_{y^i}\rho$ is the usual space derivative. However, if we take one more order covariant derivatives on \eqref{C675} and apply \eqref{CC207} to infer that
\begin{equation}\label{C664}
(D_t\rho_{,i})_{,j}=-\rho_{,ij}u^k_{,k}-\rho_{,i} u^k_{,kj}-\rho_{,j} u^k_{,ki}-\rho u^k_{,kij},
\end{equation}
$D_t$ and $(\cdot)_{,j}$ no longer commute. In fact, the definition of covariant derivatives \eqref{C42} yields
\begin{equation*}
\begin{split}
(D_t\rho_{,i})_{,j}&=\frac{\partial}{\partial y^j}D_t\rho_{,i}-\Gamma^k_{ji}\cdot\big(D_t\rho_{,k}\big)\\
&=D_t\big(\frac{\partial}{\partial y^j}\rho_{,i}-\Gamma^k_{ji}\,\rho_{,k}\big)+\big(D_t\Gamma^k_{ji}\big)\cdot\rho_{,k}\\
&=D_t(\rho_{,ij})+\rho_{,k}\,u^k_{ji},
\end{split}
\end{equation*}
where we have used \eqref{CC208}$_1$ in the last equality. Thus \eqref{C664} actually transforms into
\begin{equation}\label{C665}
\begin{split}
D_t(\rho_{,ij})&=-\rho_{,k}\,u^k_{,ji}-\rho_{,ij}\, u^k_{,k}-\rho_{,i}\,u^k_{,kj}-\rho_{,j}\,u^k_{,ki}-\rho\,u^k_{,kij}\\
&=-\big(\rho_{,k}\,u^k_{,ji}+\rho_{,i}\,u^k_{,kj}+\rho_{,j}\,u^k_{,ki}\big)-\rho_{,ij}\,u^k_{,k}-\rho\,u^k_{,kij},
\end{split}
\end{equation}
which provides the time derivatives of $\rho_{,ij}$.

Then, let us consider the decomposition $u=w+v$ given by \eqref{CC201}, and the singular part $v$ solves the Lam\'{e} system \eqref{CC201}$_2$. In Lagrangian coordinates, it is given by
\begin{equation}\label{62}
g^{mn}v_{k,mn}+v^s_{,sk}=\rho_{,k}\,.
\end{equation}
According to Lemma \ref{L27} and \eqref{C401}, the above equation changes into
\begin{equation}\label{63}
g^{mn}(v_{k,m}-v_{m,k})_{,n}+\left(2 v^s_{,s}-\rho\right)_{,k}=0.\\
\end{equation}
Multiplying \eqref{63} by $g^{ik}$ and taking one more order covariant derivatives, we apply \eqref{CC207} again to derive that
\begin{equation*}
g^{ik}g^{mn}(v_{k,m}-v_{m,k})_{,ni}+g^{ik}\left(2v^s_{,s}-\rho\right)_{,ki}=0.\\
\end{equation*}
Because $g^{ik}g^{mn}(v_{k,m}-v_{m,k})_{,ni}=0$ due to the symmetry, it holds that
\begin{equation}\label{CC223} 
g^{ik}\left(2v^s_{,s}-\rho\right)_{,ki}=\Delta F_v=0,
\end{equation}
where $F_v\triangleq(2v^s_{,s}-\rho)$ and $\Delta$ is the Laplacian-Beltrami operator given by  \cite[Definition 6.1]{GTM94}. Thus $F_v$ is harmonic in $\Omega$ (such fact is quite direct in  Eulerian coordinates).
 
Similarly, according to \cite[Section 2.20, Definition 6.1]{GTM94}, we can rewrite \eqref{62} as
\begin{equation}\label{623}
\big(d*d*+*d*d\big)v+ d\,(u^k_{,k})=d\,\rho,
\end{equation}
where $d$ is the exterior differential operator, $*$ is the Hodge $*$-operator, and
$v$ is viewed as a 1-form respectively.
Operating the exterior differential operator on \eqref{623} and utilizing the fact $d\circ d=0$ 
(\cite[Section 2.20]{GTM94}) yield that
\begin{equation*} 
d*d*d\,v=0\Leftrightarrow \big(d*d*+*d*d\big)dv=\Delta(d v)=0.
\end{equation*}
Consequently, the 2-form $($\cite[Definition 2.15]{GTM94}$)$ $\mathrm{rot}v=dv$ is harmonic as well. Note that the components of $dv$ are given by
\begin{equation*}
\big(\mathrm{rot}v\big)_{ij}=\frac{1}{2}\big(v_{i,j}-v_{j,i}\big).
\end{equation*}

\subsection{The partition of $\Omega$} 
\quad 
We introduce some notations in Lagrangian coordinates $(y,t)$, analogous to \eqref{CC104} and \eqref{CC105},
\begin{equation}\label{CC213}
\begin{split}
F_r&\triangleq\{y\in\Omega|\,\mathrm{dist}_E(y,\partial\Omega)=r\},\\
\Omega_{k}&\triangleq\{y\in\Omega\big|0 \leq\mathrm{dist}_E(y,\partial\Omega)\leq k\big\},\\
\Omega^i_j&\triangleq\big\{y\in\Omega\big|i \leq\mathrm{dist}_E(y,\partial\Omega)\leq j \big\},\\
\end{split}
\end{equation} 
where $\mathrm{dist}_E(y,y')=(\sum_i|y_i-y'_i|^2)^{1/2}$ is the Euclidean distance in Lagrangian coordinates. Note that $\Omega_k$ and $\Omega^i_j$ are fixed in the Lagrangian coordinates, which in fact correspond to the moving domains $X_t(\Omega_k)$ and $X_t(\Omega^i_j)$ in the Eulerian coordinates respectively. Moreover, $X_t$ is the identity mapping at $t=0$.
\begin{remark}\label{R25}
We mention that $\mathrm{dist}_{E}(\cdot,\cdot)$ is \textbf{not} the geodesic distance $($see Chapter V \S 3 and Chapter VII \S 5 of \cite{Bo}$)$ in Lagrangian coordinates except for $t=0$. In fact, the geodesic distance coincides with the Euclidean distance in Eulerian coordinates which measures the  distance in the physical world. 
Through out this paper, the single word ``distance" and the symbol ``\ $\mathrm{dist}(\cdot,\cdot)$" refer to the geodesic distance. 
\end{remark}

We denote $B\triangleq\Omega_5$ as the boundary domain, $I=\Omega\setminus B $ as the inner domain, and $\Gamma\triangleq\Omega^1_9$ as the intermediate layer, thus $F_5$ separates the inner and boundary parts of $\Omega$.
To carry out local estimates in $B$, we   introduce some truncation functions in   Lagrangian coordinates and compute the derivatives of them. The module $|\cdot|$ below is given by \eqref{C208}.

\begin{lemma}\label{lm29}
Suppose that $f(y)$ is defined in Lagrangian coordinates and independent of $t$, then $\forall t\geq 0$, there is a generic constant $C$ such that the point-wise estimate of $|\nabla f|$ holds,
\begin{equation}\label{44}
\sup_{0\leq\tau\leq t}|\nabla f(y,\tau)|\leq C\,|\nabla f(y,0)|\cdot \exp\big(\int_0^t|\nabla u|(y,\tau)\,d\tau\big),~~\forall y\in\Omega.
\end{equation}
We mention that
$|\nabla f(y,0)|^2=\sum_{i}\big(\partial f/\partial y^i\big)^2(y).$
\end{lemma}
\begin{proof}
By taking $D_t$ on the equation $g^{ij}g_{jk}=\delta_{ik}$ and applying \eqref{CC208}$_1$, we argue that
\begin{equation}\label{CC211}
D_tg^{ij}=-g^{ik}u^j_{,k}-g^{jk}u^i_{,k},
\end{equation}
which together with the fact $D_tf=0$ implies that $\forall y\in\Omega$,
\begin{equation*}
\begin{split}
D_t|\nabla f|^2=D_t(g^{ij}\,f_{,i}f_{,j})=(D_tg^{ij})\,f_{,i}f_{,j}
=-2\,g^{ik}u^j_{,k}f_{,i}f_{,j}\leq C\,|\nabla u|\cdot|\nabla f|^2.
\end{split}
\end{equation*}
Thus Gr\"{o}nwall's inequality yields that $\forall t\geq 0$, there is a generic constant $C$ such that
$$\sup_{0\leq\tau\leq t}|\nabla f(y,\tau)|^2\leq C\,|\nabla f(y,0)|^2\cdot\exp\big(\int_0^t|\nabla u|(y,\tau)\,d\tau\big),$$
which gives \eqref{44} and finishes the proof.
\end{proof}

As a direct corollary of Lemma \ref{lm29}, we consider a typical truncation function $\psi(y)\in C^\infty(\Omega)$ which satisfies that
\begin{equation*}
\begin{cases}
\psi=1~~\forall y\in\Omega_1,\\
\psi=0~~\forall y\in\Omega\setminus\Omega_2,\\
0\leq\psi\leq 1,~\sum_i\big(\partial\psi/\partial y^i\big)^2\leq 1~~\forall y\in\Omega.
\end{cases}
\end{equation*}
Observe that $\partial\psi/\partial y^i=0$ in $\Omega\setminus\Omega_2$, thus \eqref{44} ensures that $\forall y\in\Omega$ and $\forall t\geq 0$,
\begin{equation}\label{CC212}
\sup_{0\leq\tau\leq t}|\nabla\psi(y,\tau)|\leq   C\cdot\exp\big(\int_0^t\|\nabla u\|_{L^\infty(\Omega_2)}\,d\tau\big).
\end{equation}
\begin{remark}
Note that $\psi(y)$ is fixed in the Lagrangian coordinates, however it corresponds to $\psi(X_t^{-1}(x))$ in Eulerian coordinates, which provides a truncation of the moving domain $X_t(\Omega_2)$.
\end{remark}

We end up this subsection by establishing a proper control of the (geodesic) distance in the Lagrangian coordinates.

\begin{lemma}\label{lm43}
Let $y_1,y_2\in\Omega$ be two points in Lagrangian coordinates, and $\mathrm{dist}\,(y_1,y_2)|_{t=0}=d_0$ be the initial distance, then it holds that $\forall t\geq 0$,
\begin{equation}\la{op0i9}
\mathrm{dist}\,(y_1,y_2)(t)\geq d_0/2,
\end{equation}
provided that the velocity field $u$ satisfies
\begin{equation}\label{C46}
\int_0^\infty\|u(\cdot,t)\|_{L^\infty(\Omega)}dt\leq d_0/4.
\end{equation}
\end{lemma}
\begin{proof}
We observe that $y_1$ and $y_2$ in Lagrangian coordinates correspond to $X_t(y_1)$ and $X_t(y_2)$ in Eulerian coordinates respectively. Thus in view of Remark \ref{R25}, we argue that
$$\mathrm{dist}\,(y_1,y_2)(t)=\mathrm{dist}\big(X_t(y_1),X_t(y_2)\big).$$
where the right hand side is the Euclidean distance in Eulerian coordinates. Moreover, we have
\begin{equation*}
\frac{d}{dt}\big(X_t(y_1)-X_t(y_2)\big)=u(X_t(y_1),t)-u(X_t(y_2),t).
\end{equation*}
Multiplying the above equation by $(X_t(y_1)-X_t(y_2))$ yields that
\begin{equation}\label{57}
\frac{d}{dt}\,\mathrm{dist}\big(X_t(y_1),X_t(y_2)\big)=\big(u(X_t(y_1),t)-u(X_t(y_2),t)\big)\cdot\frac{(X_t(y_1)-X_t(y_2))}{\mathrm{dist} \big(X_t(y_1),X_t(y_2)\big)},
\end{equation}
since $\mathrm{dist}\big(X_t(y_1),X_t(y_2)\big)^2= \sum_i\big(X_t(y_1)_i-X_t(y_2)_i\big)^2 .$

Integrating \eqref{57} with respect to $t$, we utilize \eqref{C46} to declare that
\begin{equation*}
\begin{split}
\mathrm{dist}\big(X_t(y_1),X_t(y_2)\big)&\geq\mathrm{dist}\big(X_0(y_1),X_0(y_2)\big)-2\int_0^t \|u\|_{L^\infty}ds\geq d_0/2,
\end{split}
\end{equation*} 
which gives \eqref{op0i9} and finishes the proof.
\end{proof}

With the help of Lemma \ref{lm43}, we obtain the lower bound of the distance between free boundaries $F_\alpha$ with $\alpha\in[0,10)$.  

\begin{lemma}\label{512}
Suppose that the velocity field $u$ satisfies 
\begin{equation}\label{C314}
\int_0^\infty\|u(\cdot,t)\|_{L^\infty(\Omega)}dt<1/100.
\end{equation}
Then if $\alpha,\beta\in[0,10]$ and $|\alpha-\beta|\geq 1/2$, it holds that $\forall t\geq 0$,
\begin{equation*}
\mathrm{dist}\big(F_\alpha,F_\beta\big)(t)>1/4.
\end{equation*}
\end{lemma}
\begin{proof}
According to \eqref{CC213}$_1$, the initial distance between $F_\alpha$ and $F_\beta$ is $|\alpha-\beta|\geq 1$, thus by virtue of Lemma \ref{lm43} and \eqref{C314}, the distance between $F_\alpha$ and $F_\beta$ will be strictly large than $1/4$, and the proof is therefore completed.
\end{proof}

\section{A priori estimates I: Basic energy estimates}\label{sec5}
\quad Let $(\rho, u)$ be the strong solution to \eqref{11}--\eqref{CC102} on $\Omega\times[0, T]$ given by Lemma \ref{u21}. We establish the following lemma which extends the local solution $(\rho,u)$ to the global one. 
 
\begin{lemma}\label{T71}
Under the conditions of Theorem \ref{T1}, there is a constant $\varepsilon$ depending only on $\underline{\rho}$, $M_0$, and $\Omega$, such that if $\forall t\in[0,T]$ it holds that 
\begin{equation}\label{C501}
\begin{split}
\|\rho(\cdot,t)\|&_{L^\infty(\Omega)}\leq4M_0,~~\|\nabla \rho(\cdot,t)\|_{L^4(B)}\leq 4M_0,~~\|\rho^{-1}(\cdot,t)\|_{L^\infty(\Omega_{10})}\leq 4 \underline{\rho}^{-1},
\end{split}
\end{equation}
then we argue that, in fact $\forall t\in[0,T]$,
\begin{equation}\label{C502}
\begin{split}
\|\rho(\cdot,t)\|&_{L^\infty(\Omega)}\leq 2M_0,~~\|\nabla \rho(\cdot,t)\|_{L^4(B)}\leq 2M_0,~~\|\rho^{-1}(\cdot,t)\|_{L^\infty(\Omega_{10})}\leq 2 \underline{\rho}^{-1},
\end{split}
\end{equation}
provided that the initial energy satisfies $C_0\le\varepsilon$.
\end{lemma} 
We will carry out the proof after finishing a priori estimates in Section \ref{sec5} and Section \ref{S4}.
The constants $C$ below may vary from line to line, but they are all independent of $T$.

\begin{lemma}\label{41}
Suppose that \eqref{C501} is valid, then there are  constants $C$ and $\kappa_1$ depending only on  $M_0$ and $\Omega$, such that $\forall t\in[0,T]$,
\begin{equation}\label{C301}
\begin{split}
&\|\sqrt{\rho}u(\cdot,t)\|_{L^2(\Omega)}^2+\|(\rho-\bar{\rho})(\cdot,t)\|_{L^2(\Omega)}^2\leq CC_0e^{-\kappa_1t},\\
&\qquad\quad\int_0^te^{\kappa_1\tau/2}\,\|\nabla u(\cdot,\tau)\|^2_{L^2(\Omega)}\,d\tau\leq CC_0.
\end{split}
\end{equation}
\end{lemma}

\begin{proof}
Multiplying \eqref{11}$_2$ by $u$ and integrating  over $\Omega$ yield that
\begin{equation}\label{lj04}
\begin{split}
&\frac{d}{dt}\int_\Omega\left(\frac{1}{2}\rho|u|^{2}+G(\rho)\right)dx+\int_\Omega\big(|\nabla u|^{2}+(\mathrm{div}u)^2\big)\,dx=0,
\end{split}
\end{equation}
where $G(\rho)$ is given by \eqref{C108}.
Integrating \eqref{lj04} with respect to $t$ gives
\begin{equation}\label{lj07}
\sup_{0\leq\tau\leq t}\int_\Omega\left(\frac{1}{2}\rho|u|^{2}+(\rho-\bar{\rho})^2\right)dx+\int_0^t\int_\Omega|\na u|^{2}dxd\tau\leq CC_0,
\end{equation}
since there is a constant $\mathcal{C}_1$ depending only on $M_0$, such that $\forall\rho\in[0,M_0],$
\begin{equation}\label{CC301}
\mathcal{C}_1^{-1}(\rho-\bar{\rho})^2\leq G(\rho)\leq \mathcal{C}_1(\rho-\bar{\rho})^2.
\end{equation}

Next, multiplying \eqref{11}$_2$ by $-\mathcal{B}(\rho-\bar{\rho})$ with $\mathcal{B}$ given by Lemma \ref{lemmaz01} and utilizing Lemma \ref{lemmaz01} together with \eqref{C501} and \eqref{lj07} lead to
\begin{equation}\label{402}
\begin{split}
& \int_\Omega(\rho-\bar{\rho})^2\,dx-\frac{d}{dt}\int_\Omega\rho u\cdot \mathcal{B}(\rho-\bar{\rho})\,dx\\
&=\int_\Omega\rho u\cdot \mathcal{B}(\mathrm{div}(\rho u))\,dx -\int_\Omega\rho u_iu_j\cdot\partial_{i}\mathcal{B}_j(\rho-\bar{\rho})\,dx
+\int_\Omega\nabla u\cdot\nabla\mathcal{B}(\rho-\bar{\rho})\,dx\\
&\quad+\int_\Omega \mathrm{div}u\cdot(\rho-\bar{\rho})\,dx\\
&\leq C\bigg(\|\rho u\|_{L^2(\Omega)}^2+\|u\|_{L^4(\Omega)}^2\|\rho-\bar{\rho}\|_{L^2(\Omega)}
+\|\nabla u\|_{L^2(\Omega)}^2\bigg)+\frac{1}{2}\|\rho-\bar{\rho}\|_{L^2(\Omega)}^2\\
&\leq C\|\nabla u\|_{L^2(\Omega)}^2+\frac12\|\rho-\bar{\rho}\|_{L^2(\Omega)}^2.
\end{split}
\end{equation}
In view of Lemma \ref{lemmaz01}, \eqref{CC301} and Sobolev's inequality, there is a constant $\mathcal{C}_2$ depending on $M_0$ and $\Omega$, such that
\begin{equation}\label{C304}
\mathcal{C}_2^{-1} \int_\Omega|\rho u\cdot \mathcal{B}(\rho-\bar{\rho})|\,dx\leq \int_\Omega\big(\frac{1}{2}\rho |u|^2+G(\rho)\big)\,dx
\leq \mathcal{C}_2\left(\|\na u\|_{L^2(\Omega)}^2+\|\n-\bar\n\|_{L^2(\Omega)}^2\right).
\end{equation}

Thus adding \eqref{402} mutilpied by $(2\mathcal{C}_2)^{-1}$ to \eqref{lj04} and utilizing \eqref{C304} lead to
\begin{equation}\label{C308}
\begin{split}
\mathbf{E}'(t)+\kappa_1\mathbf{E}(t)\leq 0,~~\mathbf{E}(t)\triangleq
\int_\Omega\left(\frac{1}{2}\rho|u|^{2}+G(\rho)-
\frac{\rho u\cdot\mathcal{B }(\rho-\bar{\rho})}{2\mathcal{C}_2}\right)\,dx,
\end{split}
\end{equation}
for some $\kappa_1$ depending only on $\Omega$ and $M_0$.
By virtue of \eqref{CC301} and \eqref{C304}, we apply Gr\"{o}nwall's inequality in \eqref{C308} to infer that $\forall t\in[0,T]$,
\begin{equation*}
\mathbf{E}(t)\leq CC_0e^{-\kappa_1 t},
\end{equation*}
which along with \eqref{CC301} and \eqref{C304} gives \eqref{C301}$_1$
\begin{equation*} 
\|\sqrt{\rho}u(\cdot,t)\|_{L^2(\Omega)}^2+\|(\rho-\bar{\rho})(\cdot,t)\|_{L^2(\Omega)}^2\leq CC_0e^{-\kappa_1 t}.
\end{equation*}

Finally, by multiplying \eqref{lj04} by $e^{\kappa_1\tau/2}$ and integrating with respect to $\tau$, we take advantage of \eqref{C501} and \eqref{lj07} to deduce that $\forall t\in[0,T]$
\begin{equation*}
\begin{split}
\int_0^te^{\kappa_1\tau/2}\,\|\nabla u(\cdot,\tau)\|^2_{L^2(\Omega)}\,d\tau
&\leq CC_0+\int_0^t\int_\Omega e^{\kappa_1\tau/2}\,\big(\rho|u|^2+(\rho-\bar{\rho})^2\big)\,dx d\tau\leq CC_0,
\end{split}
\end{equation*} 
which gives \eqref{C301}$_2$ and finishes the proof.
\end{proof}

Next, by modifying some ideas due to Hoff   \cite{H3, Hof2, Ho3, HS1}, we have the following estimates on the material derivatives $\dot u$.
\begin{lemma}\label{21}
Suppose that \eqref{C501} is valid, then
there are constants $\varepsilon_0$, $\kappa_2$, and $C$ depending on $M_0$ and $\Omega$, such that for any $t\in[0,T]$ and $\sigma(t)\triangleq\min\{t,1\},$
\begin{equation}\label{C305}
\begin{split}
&\sup_{0\leq\tau\leq t}e^{\kappa
_2\tau}\|\nabla u(\cdot,\tau)\|_{L^2(\Omega)}^2
+\int_0^te^{\kappa_2\tau}\|\sqrt{\rho}\dot{u}(\cdot,\tau)\|_{L^2(\Omega)}^2\,d\tau\leq C,\\
&~~\sup_{0\leq\tau\leq t}\sigma\,\|\sqrt{\rho}\dot{u}(\cdot,\tau)\|_{L^2(\Omega)}^2+\int_0^t\sigma\|\nabla\dot{u}(\cdot,\tau)\|_{L^2(\Omega)}^2\,d\tau\leq C,\\
\end{split}
\end{equation}
provided that $C_0\leq\varepsilon_0$.
\end{lemma}
\begin{proof}
Let us rewrite \eqref{11}$_2$ as
\begin{equation}\label{302}
\rho\dot{u}-\Delta u-\nabla\mathrm{div} u+\nabla\rho=0,
\end{equation}
then multiplying \eqref{302} by $\dot{u}$ and integrating over $\Omega$ give that
\begin{equation}\label{CC302}
\begin{split}
\int_\Omega\rho|\dot{u}|^2\,dx
&=\int_\Omega\Delta u\cdot\dot{u}\,dx+\int_\Omega\nabla\div u\cdot\dot{u}\,dx-\int_\Omega\dot{u}\cdot\nabla(\rho-\bar{\rho})\,dx.
\end{split}
\end{equation}

The first term on the right side of \eqref{CC302} is treated by
\begin{equation*} 
\begin{split}
&\int_\Omega\Delta u\cdot\dot{u}\,dx\\
&=-\int_\Omega\nabla u\cdot\nabla\dot{u}\,dx\\
&=-\frac{d}{dt}\left(\frac{1}{2}\int_\Omega|\nabla u|^{2}\,dx\right)-\int_\Omega\nabla u\cdot\nabla(u\cdot\nabla u)\,dx \\
&=-\frac{d}{dt}\left(\frac{1}{2}\int_\Omega|\nabla u|^{2}\,dx\right)+\frac{1}{2}\int_\Omega\mathrm{div}u\cdot|\nabla u|^2\,dx-\int_\Omega\partial_ju^k\cdot\partial_ju^i
\cdot\partial_iu^k\,dx\\
&\leq-\frac{d}{dt}\left(\frac{1}{2}\int_\Omega|\nabla u|^{2}\,dx\right)+C\int_\Omega|\nabla u|^3\,dx.
\end{split}
\end{equation*}
Similar process also yields that 
\begin{equation*} 
\begin{split}
\int_\Omega\nabla\div u\cdot\dot{u}\,dx
\leq-\frac{d}{dt}\left(\frac{1}{2}\int_\Omega(\mathrm{div}u)^{2}\,dx\right)+C\int_\Omega|\nabla u|^3\,dx.
\end{split}
\end{equation*}
While in view of \eqref{11}$_1$, \eqref{C501}, and \eqref{CC102}, we argue that the third term is handled via
\begin{equation*} 
\begin{split}
&-\int_\Omega\dot{u}\cdot\nabla(\rho-\bar{\rho})\,dx\\
&=\int_\Omega(\rho-\bar{\rho})\cdot\div u_{t}\,dx+ \int_\Omega(\rho-\bar{\rho})\,\mathrm{div}(u\cdot\nabla u)\,dx\\
&=\frac{d}{dt}\left(\int_\Omega(\rho-\bar{\rho})\cdot\mathrm{div}u\,dx\right)
+\int_\Omega\rho\,\partial_iu^j\cdot\partial_ju^i\,dx\\
&\leq\frac{d}{dt}\left(\int_\Omega(\rho-\bar{\rho})\cdot\mathrm{div}u\,dx\right)+C\|\nabla u\|_{L^2(\Omega)}^2.
\end{split}
\end{equation*}
By substituting these results into \eqref{CC302}, we add \eqref{C308} multiplied by $10\mathcal{C}_1$ to \eqref{CC302} and arrive at
\begin{equation}\label{3006}
\begin{split}
&\mathbf{F}'(t)+\int_\Omega\rho|\dot{u}|^2\,dx\leq C\int_\Omega(|\nabla u|^2+|\nabla u|^3)\,dx.\\
\end{split}
\end{equation}
The leading term $F(t)$ is defined by
\begin{equation*}
\begin{split}
\mathbf{F}(t)&\triangleq\int_\Omega\bigg(|\nabla u|^2+(\mathrm{div}u)^2-(\rho-\bar{\rho})\cdot\mathrm{div}u\bigg)\,dx+10\,\mathcal{C}_1\mathbf{E}(t)\\
&\geq\int_\Omega\big(|\nabla u|^2+(\rho-\bar{\rho})^2\big)\,dx,\\
\end{split}
\end{equation*}
where $\mathcal{C}_1$ and $\mathbf{E}(t)$ are given by \eqref{CC301} and \eqref{C308} respectively.

Next, operating $\dot u^j\big(\frac{\partial}{\partial t}+\div(u\cdot)\big)$ to $\eqref{302}$, summing with respect to $j$ and integrating over $\Omega$ leads to
\begin{align}\label{CC303}
\frac{d}{dt}\bigg(\frac{1}{2}\int_\Omega\rho|\dot u|^2\,dx\bigg) 
+\int_\Omega\bigg(|\nabla\dot u|^2+(\mathrm{div}\dot u)^2\bigg)\,dx=J_1+J_2+J_3.
\end{align}
The right side of \eqref{CC303} is given by  
\begin{equation*} 
\begin{split}
J_1&=-\int_\Omega\bigg(\nabla(u_i\cdot\partial_i u)\cdot\nabla\dot u-\Delta u\cdot(u_i\cdot\partial_i\dot u)\bigg)dx,\\
J_2&=-\int_\Omega\bigg(\mathrm{div}(u_i\cdot\partial_i u)\cdot\mathrm{div}\dot u-\nabla\mathrm{div}u\cdot(u_i\cdot\partial_i\dot u)\bigg)dx,\\
J_3&=\int_\Omega\bigg(\rho_t\cdot\mathrm{div}\dot u+\nabla \rho\cdot(u_i\cdot\partial_i\dot u)\bigg)dx.
\end{split}
\end{equation*}
Direct calculations yield that
\begin{equation}\label{CC304}
\begin{split}
J_1&=-\int_\Omega\bigg(\nabla(u_i\cdot\partial_i u)\cdot\nabla\dot u-\Delta u\cdot(u_i\cdot\partial_i\dot u)\bigg)dx\\
&\leq-\int_\Omega\bigg((u_i\cdot\partial_i)\nabla u\cdot\nabla\dot u+\nabla u\cdot (u_i\cdot\partial_i)\nabla\dot u\bigg)\,dx+C\int_\Omega|\nabla u|^2\cdot|\nabla\dot u|\,dx\\
&\leq C\int_\Omega|\nabla u|^2\cdot|\nabla\dot u|\,dx\leq C_\Omega\int|\nabla u|^4\,dx+\frac{1}{10}\int_\Omega|\nabla\dot u|^2\,dx.
\end{split}
\end{equation}
Similarly, we also infer that
\begin{equation}\label{CC305}
\begin{split}
J_2\leq C\int_\Omega|\nabla u|^4\,dx+\frac{1}{10}\int_\Omega|\nabla\dot u|^2\,dx.
\end{split}
\end{equation}
With the help of \eqref{11}$_1$ and \eqref{C501}, we deduce that
\begin{equation}\label{CC306}
\begin{split}
J_3&=\int_\Omega\bigg(\rho_t\cdot\mathrm{div}\dot u+\nabla\rho\cdot(u_i\cdot\partial_i\dot u)\bigg)dx\\
&=\int_\Omega\bigg(\rho_t\cdot\mathrm{div}\dot u-\big(\rho-\bar{\rho})\big)\big(u_i\cdot\partial_i\mathrm{div}\dot u+ \partial_ju_i\cdot\partial_i\dot u_j\big)\bigg)\,dx\\
&\leq\int_\Omega\bigg(-\bar{\rho}\,\mathrm{div}u\cdot\mathrm{div}\dot u+|\rho-\bar{\rho}|\cdot|\nabla u|\cdot|\nabla\dot u|\bigg)\,dx\\
&\leq C\int_\Omega\bigg(|\nabla u|^4+|\nabla u|^2+ |\rho-\bar{\rho}|^4\bigg)\, dx+\frac{1}{10}\int|\nabla\dot u|^2\,dx.
\end{split}
\end{equation}
Combining \eqref{C501}, \eqref{CC304}, \eqref{CC305} and \eqref{CC306} shows that
\begin{equation*} 
\begin{split}
\frac{d}{dt}\int_\Omega\rho|\dot{u}|^2dx+\int_\Omega|\nabla\dot{u}|^2\,dx
\leq C\int\left(|\nabla u|^4+|\nabla u|^2+|\rho-\bar{\rho}|^2\right)\,dx.
\end{split}
\end{equation*}
which multiplied by $\sigma(t)$ gives that
\begin{equation}\label{C103}
\begin{split}
&\frac{d}{dt}\left(\sigma\int_\Omega\rho|\dot{u}|^2\,dx\right)+\sigma \int_\Omega|\nabla\dot{u}|^2\,dx\\
&\leq C \int_\Omega\sigma\left(|\nabla u|^4+|\nabla u|^2+(\rho-\bar{\rho})^2\right)dx+\sigma'\int_\Omega\rho|\dot{u}|^2\,dx.
\end{split}
\end{equation}

According to the decomposition \eqref{CC214}, elliptic estimate \eqref{355} and \eqref{GN}, we argue that
\begin{equation}\label{3008}
\begin{split}
&\int_\Omega|\nabla u|^3\,dx\\
&\leq C\int_\Omega|\nabla w|^3\,dx+C\int_\Omega|\nabla v|^3\,dx\\
&\leq C\bigg(\|\nabla w\|_{L^2}^{3/2}\|\nabla^2 w\|^{3/2}_{L^2(\Omega)}+\|\nabla w\|_{L^2(\Omega)}^3+\|\rho-\bar{\rho}\|^3_{L^3(\Omega)}\bigg)\\
&\leq C\bigg(\|\nabla u\|_{L^2(\Omega)}+\|\rho-\bar{\rho}\|_{L^2(\Omega)}\bigg)^{3/2}\|\,\rho\dot{u}\|^{3/2}_{L^2(\Omega)}
+C\bigg(\|\nabla u\|_{L^2(\Omega)}^3+\|\rho-\bar{\rho}\|^3_{L^3(\Omega)}\bigg),\\
&\int_\Omega|\nabla u|^4\,dx\\
&\leq C\int_\Omega|\nabla w|^4\,dx+C\int_\Omega|\nabla v|^4\,dx\\
&\leq C\bigg(\|\nabla w\|_{L^2(\Omega)}\|\nabla^2 w\|^3_{L^2(\Omega)}+\|\nabla w\|_{L^2(\Omega)}^4+\|\rho-\bar{\rho}\|^4_{L^4(\Omega)}\bigg)\\
&\leq C\bigg(\|\nabla u\|_{L^2(\Omega)}+\|\rho-\bar{\rho}\|_{L^2(\Omega)}\bigg)\,\|\rho\dot{u}\|^3_{L^2(\Omega)}+C\,\bigg(\|\nabla u\|_{L^2(\Omega)}^4
+\|\rho-\bar{\rho}\|^4_{L^4(\Omega)}\bigg).\\
\end{split}
\end{equation}
 
Substituting \eqref{3008}$_1$ into \eqref{3006} and applying \eqref{C501} yield 
\begin{equation}\label{C313}
\begin{split}
&\mathbf{F}'(t)+\int_\Omega\rho|\dot{u}|^{2}dx \\
&\leq C\|\nabla u\|_{L^2(\Omega)}^2
+C\left(\big(\|\nabla u\|_{L^2(\Omega)}+\|\rho-\bar{\rho}\|_{L^2(\Omega)}\big)^{3/2}\,\|\rho\dot{u}\|^{3/2}_{L^2(\Omega)}
+\|\rho-\bar{\rho}\|^3_{L^3(\Omega)}\right)\\
&\leq \frac{1}{2}\int_\Omega\rho|\dot{u}|^2dx+C\bigg(\|\nabla u\|_{L^2(\Omega)}^6
+\|\nabla u\|_{L^2(\Omega)}^2+\|\rho-\bar{\rho}\|_{L^2(\Omega)}^2\bigg).
\end{split}
\end{equation}
Let $\kappa_2=\kappa_1/10$ with $\kappa_1$ given in Lemma \ref{41}, then we multiply \eqref{C313} by $e^{\kappa_2t}$ and make use of \eqref{3006} to derive that
\begin{equation}\label{C34}
\begin{split}
& \frac{d}{dt}\big(e^{\kappa_2 t}\,\mathbf{F}(t)\big)+e^{\kappa_2 t}/2\int_\Omega\rho|\dot{u}|^{2}\,dx\\
&\leq C\|\nabla u\|_{L^2(\Omega)}^2\big(e^{\kappa_2 t}\,\mathbf{F}(t)\big)^2+Ce^{\kappa_2t}\bigg(\|\nabla u\|_{L^2(\Omega)}^2+\|\rho-\bar{\rho}\|_{L^2(\Omega)}^2\bigg).
\end{split}
\end{equation}
Consequently, by choosing $C_0\leq\varepsilon_0$ for some $\varepsilon_0$ determined by $M_0$ and $\Omega$, we apply Lemma \ref{41} and Gr\"{o}nwall's inequality in \eqref{C34} to declare  
\begin{equation*}
\sup_{0\leq\tau\leq t}e^{\kappa_2\tau}\|\nabla u(\cdot,\tau)\|_{L^2(\Omega)}^2
+\int_0^te^{\kappa_2\tau}\|\sqrt{\rho}\dot{u}(\cdot,\tau)\|_{L^2(\Omega)}^2\,d\tau\leq C, 
\end{equation*}
which provides \eqref{C305}$_1$ and finishes the first part.

Next, substituting \eqref{3008}$_2$ into \eqref{C103} shows that
\begin{equation}\label{CC318}
\begin{split}&
\frac{d}{dt}\bigg(\sigma \int_\Omega\rho|\dot{u}|^2\,dx\bigg)+\sigma\int_\Omega|\nabla\dot{u}|^{2}\,dx\\
&\leq C\big(\|\nabla u\|_{L^2(\Omega)}+\|\rho-\bar{\rho}\|_{L^2(\Omega)}\big)\,\|\rho\dot{u}\|_{L^2(\Omega)}\left(\sigma \int_\Omega\rho|\dot{u}|^2\,dx\right)\\
&\quad+C\bigg(\|\nabla u\|_{L^2(\Omega)}^4+\|\nabla u\|^2_{L^2(\Omega)}+\|\rho-\bar{\rho}\|^4_{L^4(\Omega)}+\sigma'\,\|\sqrt{\rho}\dot{u}\|_{L^2(\Omega)}^2\bigg).
\end{split}
\end{equation}
Thus we utilize Lemma \ref{41}, \eqref{C305}$_1$, and Gr\"{o}nwall's inequality to derive that 
\begin{equation*}
\sup_{0\leq\tau\leq t}\sigma\|\sqrt{\rho}\dot{u}(\cdot,\tau)\|_{L^2(\Omega)}^2+\int_0^t\sigma\|\nabla\dot{u}(\cdot,\tau)\|_{L^2(\Omega)}^2\,d\tau\leq C,
\end{equation*}
which gives \eqref{C305}$_2$ and finishes the proof.
\end{proof}
Recalling the decomposition $u=w+v$ given by \eqref{CC214}, we utilize Lemmas \ref{41} and \ref{21} to derive following estimates on the $L^1\big(0,T;L^\infty(\Omega)\big)$-norm of both $u$ and $\nabla w$.
\begin{lemma}\label{31}
Suppose that \eqref{C501} and \eqref{C305} are valid, then there is a constant $C$ depending only on $M_0$ and $\Omega$, such that  
\begin{equation}\label{C527}
\int_{0}^T\left(\|u(\cdot,t)\|_{L^\infty(\Omega)}+\|\nabla w(\cdot,t)\|_{L^{\infty}(\Omega)}\right) dt\leq CC_0^{1/14}.
\end{equation}
\end{lemma}
\begin{proof}
According to \eqref{GN}, \eqref{355}, Lemmas \ref{41} and \ref{21}, we check directly that
\begin{equation}\label{3113}
\begin{split}&
\int_{0}^{T}\|\nabla w(\cdot,t)\|_{L^{\infty}(\Omega)}\,dt\\
&\leq C\int_0^T\left(\|\nabla w\|^{1/7}_{L^2(\Omega)}\|\nabla^2 w\|^{6/7}_{L^4(\Omega)}+\|\nabla w\|_{L^2(\Omega)}\right)\,dt\\
&\leq C\int_0^T\left(\big(\|\nabla u\|_{L^2(\Omega)}+\|\rho-\bar{\rho}\|_{L^2(\Omega)}\big)^{1/7}
\|\nabla\dot{u}\|^{6/7}_{L^2(\Omega)}
+\|\nabla u\|_{L^2(\Omega)}+\|\rho-\bar{\rho}\|_{L^2(\Omega)}\right)\,dt\\
&\leq C\int_0^T\sigma^{-3/7}\big(\|\nabla u\|_{L^2(\Omega)}+\|\rho-\bar{\rho}\|_{L^2(\Omega)}\big)^{1/7}
\big(\sigma^{1/2}\|\nabla\dot{u}\|_{L^2(\Omega)}\big)^{6/7}\,dt\\
& \quad+C\int_0^T\big(\|\nabla u\|_{L^2(\Omega)}+\|\rho-\bar{\rho}\|_{L^2(\Omega)}\big)\,dt\leq CC_0^{1/14},
\end{split}
\end{equation}
where we have taken advantage of the following facts
\begin{equation*}
\begin{split}
&\int_0^1\sigma^{-3/7}\big(\|\nabla u\|_{L^2(\Omega)}+\|\rho-\bar{\rho}\|_{L^2(\Omega)}\big)^{1/7}
\big(\sigma^{1/2}\|\nabla\dot{u}\|_{L^2(\Omega)}\big)^{6/7}\,dt\\
&\leq\big(\int_0^1t^{-6/7}\,dt\big)^{1/2}\left(\int_0^1 2\big(\|\nabla u\|_{L^2(\Omega)}^2+\|\rho-\bar{\rho}\|_{L^2(\Omega)}^2\big)\,dt\right)^{1/14}\left(\int_0^1\sigma\|\nabla\dot{u}\|_{L^2(\Omega)}^2\,dt\right)^{3/7}\\
&\leq CC_0^{1/14},\\
&\int_{1}^{T}\big(\|\nabla u\|_{L^2(\Omega)}+\|\rho-\bar{\rho}\|_{L^2(\Omega)}\big)^{1/7}
\big(\sigma^{1/2}\|\nabla\dot{u}\|_{L^2(\Omega)}\big)^{6/7}\,dt\\
&\leq\big(\int_1^Te^{-\kappa_1t/14}\,dt\big)^{1/2}\left(\int_1^T2e^{\kappa_1t/2}\big(\|\nabla u\|_{L^2(\Omega)}^2+\|\rho-\bar{\rho}\|_{L^2(\Omega)}^2\big)\,dt\right)^{1/14}\left(\int_1^T\|\nabla\dot{u}\|_{L^2(\Omega)}^2\,dt\right)^{3/7}\\
&\leq CC_0^{1/14}.
\end{split}
\end{equation*}

Now the estimate upon $\|u\|_{L^{\infty}}$ is direct. Observing that $w=0$ and $v=0$ on $\partial\Omega$, thus we apply Sobolev's inequality, \eqref{355}, \eqref{C501}, \eqref{3113}, and Lemma \ref{41} to declare that
\begin{equation*}
\begin{split}
&\int_0^T\|u(\cdot,t)\|_{L^{\infty}(\Omega)}\,dt\\
&\leq C\int_0^T\big(\|\nabla w(\cdot,t)\|_{L^{4}(\Omega)}+\|\nabla v(\cdot,t)\|_{L^4(\Omega)}\big)\,dt\\
&\leq C\int_0^T\big(\|\nabla w(\cdot,t)\|_{L^{\infty}(\Omega)}+\|(\rho-\bar{\rho})(\cdot,t)\|_{L^4(\Omega)}\big)\,dt\\
&\leq C\int_0^T\|\nabla w(\cdot,t)\|_{L^{\infty}(\Omega)}\,dt+CC_0^{1/4}\int_0^Te^{-\kappa_1t}\,dt
\leq CC_0^{1/14},
\end{split}
\end{equation*}
which together with \eqref{3113} gives \eqref{C527} and finishes the proof.
\end{proof}

Based on Lemma \ref{31}, we apply Lemma \ref{512} to provide lower bounds of the distance between free boundaries $F_r$.
\begin{lemma}\label{L35}
Suppose that \eqref{C527} is valid, we can find a constant $\varepsilon_1$ depending only on $M_0$ and $\Omega$, such that if the initial energy satisfies $C_0<\varepsilon_1$, then $\forall\alpha,\beta\in[0,10]$ with $|\alpha-\beta|\geq 1/2$,  
\begin{equation}\label{CC307}
\mathrm{dist}\big(F_\alpha,F_\beta\big)(t)>1/4,~~\forall t\geq 0.
\end{equation}
\end{lemma}
\begin{proof}
In view of \eqref{C527}, by choosing $\varepsilon_1\triangleq\min\{\varepsilon_0,(100C)^{-7/50}\}$ with $\varepsilon_0$ given by Lemma \ref{21}, we observe that the condition \eqref{C314} holds provided $C_0\leq\varepsilon_1$. Therefore Lemma \ref{512} implies that \eqref{CC307} is valid and finishes the proof.
\end{proof}

\section{A priori estimates II: Near boundary estimates}\label{S4}
\quad Supposing that \eqref{C501} and arguments in Lemmas \ref{41}--\ref{L35} are valid, we establish the local energy estimates near $\partial\Omega$.  Subsections \ref{SS41} and \ref{SS42} deal with the intermediate layer $\Gamma=\Omega^1_9$ and the boundary domain $B=\Omega_5$, then we prove Lemma \ref{T71} in Subsection \ref{SS43}.
  
\subsection{The estimates in the intermediate layer $\Gamma=\Omega^1_9$}\label{SS41}
\quad
We first consider the following $C^\alpha$ estimates of $\rho$ in $\Omega_9$ which is the crucial step.
\begin{lemma}\label{lm51}
Suppose that \eqref{C501} and Lemmas \ref{41}--\ref{L35} are valid, then there are constants $\varepsilon_2$ and $C$ depending on $\underline{\rho}$, $M_0$, and $\Omega$, such that if the initial energy satisfies $C_0\leq\varepsilon_2$, then $\forall t\in[0,T]$,
\begin{equation}\label{C12241}
\sup_{0\leq\tau\leq t}\|\rho(\cdot,\tau)\|_{C^\alpha(\Omega_9)}\leq C,
\end{equation}
which also implies that
\begin{equation}\label{618}
\int_0^T\|\nabla u(\cdot,t)\|_{L^\infty(\Omega_8)}\,dt\leq C.
\end{equation}
\end{lemma}
\begin{proof}
Recalling the modified norm $C^\alpha(\Omega_{10};\Omega^1_9)$ given by \eqref{D51},  
$$\|\rho\|_{C^\alpha(\Omega_{10}^{1/2};\Omega^1_9)}=\sup_{x\in\Omega_{10}^{1/2},y\in\Omega^1_9}\frac{|\rho(x)-\rho(y)|}{|x-y|^\alpha}+\sup_{x\in\Omega_{10}^{1/2}}|\rho(x)|,$$
we mention that $\|\rho\|_{C^\alpha(\Omega^1_9)}  $ is controlled by the above norm.

\textit{Step 1. $C^\alpha$ estimates along the flow.}

Let $X_t(x)$ and $X_t(y)$ be two flow lines satisfying $X_t(x)\subset \Omega_{10}$ and $X_t(y)\subset\Omega^1_{9}$ for all $t>0$. By virtue of \eqref{57}, we declare that
\begin{equation}\label{CC311}
\begin{split}
\frac{d}{dt}\log\big(|X_t(x)-X_t(y)|\big)
&=a(t)\cdot\frac{X_t(x)-X_t(y)}{|X_t(x)-X_t(y)|},~~a(t)\triangleq \frac{u(X(t),t)-u(X_t(y),t)}{|X_t(x)-X_t(y)|}.\\
\end{split}
\end{equation}
According to Lemma \ref{le52}, \eqref{C301}, and Lemmas \ref{41}--\ref{L35}, we argue that $\forall\varepsilon\in(0,1/1000)$,
\begin{equation*}
\begin{split}
&\int_{s}^{t}|a(\tau)|\,d\tau\\
&\leq\int_{s}^{t}\|\nabla w(\cdot,\tau)\|_{L^\infty(\Omega)}\,d\tau+\int_{s}^{t}
\frac{|v(X_\tau(x),\tau)-v(X_\tau(y),\tau)|}{|X_\tau(x)-X_\tau(y)|}\,d\tau\\
&\leq \int_{s}^{t}\|\nabla w(\cdot,\tau)\|_{L^\infty(\Omega)}\,d\tau+C\int_{s}^{t}\bigg(\varepsilon^{\alpha/2}\|\rho\|_{C^\alpha(\Omega_{10}^{1/2};\Omega^1_9)}
+\varepsilon^{-2}\|\rho-\bar{\rho}\|_{L^4(\Omega)}\bigg)\,d\tau\\
&\leq (t-s)\,\varepsilon^{\alpha/2}\cdot\left(C \sup_{0\leq\tau\leq t}\|\rho(\cdot,\tau)\|_{C^\alpha(\Omega_{10}^{1/2};\Omega^1_9)}\right)+\varepsilon^{-2}C\, C_0^{1/4}.
\end{split}
\end{equation*}
Thus by integrating \eqref{CC311} with respect to $t$, we deduce that for $0\leq s\leq t\leq T$,
\begin{equation}\label{C649}
\begin{split}
&\frac{|X_{s}(x)-X_{s}(y)|}{|X_{t}(x)-X_{t}(y)|}\\
&=\exp\left(-\int_{s}^{t}a(\tau)\cdot\frac{X_\tau(x)-X_\tau(y)}{|X_\tau(x)-X_\tau(y)|}\,d\tau\right)\\
&\leq\exp\left(\int_{s}^{t}|a(\tau)|\,d\tau\right)\\
&\leq  \exp\left((t-s)\,\varepsilon^{\alpha/2}\cdot\left(C\sup_{0\leq\tau\leq t}\|\rho(\cdot,\tau)\|_{C^\alpha(\Omega_{10};\Omega^1_9)}\right)+\varepsilon^{-2} C\,C_0^{1/4}\right)\\
&\leq e^{(t-s)/2}\cdot\exp\left(C C_0^{1/4}\sup_{0\leq\tau\leq t}\|\rho(\cdot,\tau)\|^{4/\alpha}_{C^\alpha(\Omega_{10}^{1/2};\Omega^1_9)}\right),
\end{split}
\end{equation}
where we may take $\varepsilon^{\alpha/2}= \big(2 C\sup_{0\leq\tau\leq t}\|\rho(\cdot,\tau)\|_{C^\alpha(\Omega_{10};\Omega^1_9)}\big)^{-1}$ in the last but two line.

Setting $F\triangleq 2\,\mathrm{div}u-(\rho-\bar{\rho})$, we rewrite \eqref{11}$_1$ as
\begin{equation*}
\frac{d}{dt}\log\rho\big(X_t(x)\big)+\frac{1}{2}\big(\rho\big(X_t(x)\big)-\bar{\rho}\big)=-\frac{1}{2}F\big(X_t(x)\big).
\end{equation*}
Since $\rho\geq\underline{\rho}/4$ in $\Omega_{10}$, it holds that for $L=\underline{\rho}/8$,
\begin{equation*}
\begin{split}
&\frac{d}{dt}\big|\log\rho\big(X_t(x)\big)-\log\rho\big(X_t(y)\big)\big|+L\,\big|\log\rho(X_t(x))-\log\rho(X_t(y))\big|\\
&\leq \big|F\big(X_t(x)\big)-F\big(X_t(y)\big)\big|.
\end{split}
\end{equation*}
Supposing $L>1$, we integrate the above inequality on $[0,t]$ to infer that 
\begin{equation*} 
\begin{split}
&\big|\log\rho\big(X_t(x)\big)-\log\rho\big(X_t(y)\big)\big|\\
&\leq e^{-t}\big|\log\rho_0(x)-\log\rho_0(y)\big|+\int_0^te^{-(t-s)}\big|F\big(X_s(x)\big)-F\big(X_s(y)\big)\big|\,ds.
\end{split}
\end{equation*}
Multiplying the above equation by $1/|X_t(x)-X_t(y)|^\alpha$ gives
\begin{equation}\label{619}
\begin{split}
\frac{\big|\log\rho\big(X_t(x)\big)-\log\rho\big(X_t(y)\big)\big|}{|X_t(x)-X_t(y)|^\alpha}\leq \mathbf{A}(0)\cdot \|\log\rho_0\|_{C^\alpha(\Omega_{10})}+\int_0^t\mathbf{A}(s)\cdot\|F(\cdot,s)\|_{C^\alpha(\Omega_{10}^{1/2})} ds.\\ 
\end{split}
\end{equation}
In view of \eqref{C649}, $\mathbf{A}(s)$ is given by
\begin{equation}\label{CA401}
\begin{split}
\mathbf{A}(s)&\triangleq e^{-(t-s)}\cdot\left(\frac{|X_{s}(x)-X_{s}(y)|}{|X_{t}(x)-X_{t}(y)|}\right)^\alpha\\
&\leq e^{-(t-s)/2}\cdot\exp\left(CC_0^{1/4}\sup_{0\leq\tau\leq t}\|\rho(\cdot,\tau)\|^{4/\alpha}_{C^\alpha(\Omega_{10}^{1/2};\Omega^1_9)}\right).
\end{split}
\end{equation}
Writing $F=2\,\mathrm{div}w+F_v$ with $F_v\triangleq 2\,\mathrm{div}v-(\rho-\bar{\rho})$, we invoke Sobolev's inequality, \eqref{355}, \eqref{C501}, and Lemmas \ref{41}--\ref{L35} to check that, for $\theta=(9+6\alpha)/11$ with $\alpha<1/3$, 
\begin{equation}\label{C650}
\begin{split}
&\int_0^t\|\mathrm{div} w(\cdot,s)\|_{C^\alpha(\Omega^{1/2}_{10})}e^{-(t-s)/2}ds\\
&\leq C\int_0^t\|\nabla w\|^{1-\theta}_{L^2(\Omega)}\cdot\|\nabla w\|^{\theta}_{C^{1/3}(\Omega)}e^{-(t-s)/2}\,ds\\
&\leq C\int_0^t\|\nabla w\|^{1-\theta}_{L^2(\Omega)}\cdot\|\nabla w\|^{\theta}_{W^{1,4}(\Omega)}e^{-(t-s)/2}\,ds\\
&\leq C\int_0^t\big(\|\nabla u\|_{L^2(\Omega)}+\|\rho-\bar{\rho}\|_{L^2(\Omega)}\big)^{1-\theta}\cdot\|\sqrt{\rho}\dot{u}\|^{\theta/4}_{L^2(\Omega)}\cdot\|\nabla\dot{u}\|^{3\theta/4}_{L^2(\Omega)}\,e^{-(t-s)/2}ds\\
&\leq C C_0^{(1-\theta)/2}.
\end{split}
\end{equation}
The interpolation in the second line is left to the next step.
Moreover, since $F_v$ is harmonic in $\Omega$ (see \eqref{CC223}), we apply \eqref{CC215} together with \eqref{C301} and \eqref{CC307} to declare that
\begin{equation}\label{6618}
\begin{split}
\int_0^t\|F_v\|_{C^\alpha(\Omega^{1/2}_{10})}\,e^{-(t-s)/2}\,ds\leq
C\int_0^t\|\rho-\bar{\rho}\|_{L^2(\Omega)}\,e^{-(t-s)/2}\,ds\leq CC_0^{1/2}.
\end{split}
\end{equation}
Collecting \eqref{CA401}--\eqref{6618} leads to
\begin{equation*}
\begin{split}
\int_0^t\mathbf{A}(s)\cdot\|F(\cdot,s)\|_{C^\alpha(\Omega_{10}^{1/2})} ds 
&\leq C\int_0^t\mathbf{A}(s)\cdot\left(\|\mathrm{div} w(\cdot,s)\|_{C^\alpha(\Omega_{10}^{1/2})}+\|F_v(\cdot,s)\|_{C^\alpha(\Omega_{10}^{1/2})}\right)ds\\
&\leq CC_0^{\frac{1-\theta}{2}}\cdot\exp\left(CC_0^{1/4}\sup_{0\leq\tau\leq t}\|\rho(\cdot,\tau)\|^{4/\alpha}_{C^\alpha(\Omega_{10}^{1/2};\Omega^1_9)}\right).
\end{split}
\end{equation*}
Similarly, we also deduce that
\begin{equation*} 
\begin{split} 
\mathbf{A}(0)\cdot \|\log\rho_0\|_{C^\alpha(\Omega_{10})}
&\leq C\|\rho_0\|_{C^\alpha(\Omega_{10})}\cdot\exp\left(CC_0^{1/4}\sup_{0\leq\tau\leq t}\|\rho(\cdot,\tau)\|^{4/\alpha}_{C^\alpha(\Omega_{10}^{1/2};\Omega^1_9)}\right).\\
\end{split}
\end{equation*}
 
Combining these estimates and the fact $\rho\geq\underline{\rho}/4$, we take the supremum over all $X_t(x)$ and $X_t(y)$ in \eqref{619} to infer that $\forall t\in[0,T]$,
\begin{equation*}
\|\rho(\cdot,t)\|_{C^\alpha(\Omega_{10}^{1/2};\,\Omega^1_9)}\leq C\left(C_0^{\frac{1-\theta}{2}}+\|\rho_0\|_{C^\alpha(\Omega_{10})}\right)
\exp\left(CC_0^{1/4}\sup_{0\leq\tau\leq t}\|\rho(\cdot,\tau)\|^{4/\alpha}_{C^\alpha(\Omega_{10}^{1/2};\,\Omega^1_9)}\right).
\end{equation*}
By choosing $C_0\leq\varepsilon_2$ for some $\varepsilon_2$ determined by $\underline{\rho}$, $M_0$, and $\Omega$, we arrive at 
$$\sup_{0\leq\tau\leq t}\|\rho(\cdot,\tau)\|_{C^\alpha(\Omega^{1}_{9})}\leq 2M_0,$$
which together with the fact $\|\rho\|_{C^\alpha(\Omega_3)}\leq C\,\|\rho\|_{W^{1,4}(B)}$  due to Sobolev's inequality and \eqref{CC307} yields \eqref{C12241} and finishes the first step. 

\textit{Step 2. The interpolation assertions.}

Supposing $\Gamma'\subset\Gamma$ with $\mathrm{dist}\big(\Gamma',\,\partial\Gamma\big)>d$, we declare that there is a constant $C$ depending on $d$, such that for $f\in C^\alpha(\Gamma)$ and $\theta=(3+2\beta)/(3+2\alpha)$ with $0\leq\beta<\alpha$, 
\begin{equation}\label{620}
\begin{split}
\|f\|_{C^\beta(\Gamma')}\leq C\,\|f\|_{L^2(\Gamma)}^{1-\theta}\cdot\|f\|_{C^\alpha(\Gamma)}^{\theta}+C\,\|f\|_{L^2(\Gamma)}.
\end{split}
\end{equation}
In fact, we may assume $\|f\|_{C^\alpha(\Gamma)}> 10\,d^{-3/2\theta}\|f\|_{L^2(\Gamma)}$, otherwise \eqref{620} is valid directly. 
By taking $\delta=(\|f\|_{L^2(\Gamma)}\|f\|_{C^{\alpha}(\Gamma)}^{-1})^{2\theta/(3+2\beta)}<d$, it holds that
\begin{equation}\label{CC314}
\begin{split}
|f(x)|
&\leq C\,\delta^{-3}\big(\int_{B_\delta(x)}\big|f(x)-f(y)\big|\,dy
+\int_{B_\delta(x)}|f(y)|\,dy\big)\\
&\leq C\,\big(\delta^{-3}\|f\|_{C^{\alpha}(\Gamma)}\cdot\int_{B_\delta}|y|^{\alpha}\,dy
+\delta^{-3/2}\cdot\|f\|_{L^2(\Gamma)}\big)\\
&\leq C\big(\delta^\alpha\|f\|_{C^{\alpha}(\Gamma)}+\delta^{-3/2}\|f\|_{L^2(\Gamma)}\big)\leq C\,\|f\|_{L^2(\Gamma)}^{1-\theta}\cdot\|f\|_{C^\alpha(\Gamma)}^{\theta},
\end{split}
\end{equation}
where we have applied $\|f\|_{L^2(\Gamma)}\leq C\|f\|_{C^\alpha(\Gamma)}$ in the last inequality. Therefore when $|x-y|\geq\delta$, we utilize \eqref{CC314} to deduce that
\begin{equation*}
\begin{split}
\frac{|f(x)-f(y)|}{|x-y|^\beta}&\leq \delta^{-\beta}\bigg(|f(x)|+|f(y)|\bigg)\\
&\leq C\bigg(\delta^{\alpha-\beta}\|f\|_{C^{\alpha }(\Gamma)}+\delta^{-3/2-\beta}\|f\|_{L^2(\Gamma)}\bigg)\\
&\leq C\,\|f\|_{L^2(\Gamma)}^{1-\theta}\cdot\|f\|_{C^\alpha(\Gamma)}^{\theta}, 
\end{split}
\end{equation*}
when $|x-y|\leq\delta$, we directly check that
\begin{equation*}
\begin{split}
\frac{|f(x)-f(y)|}{|x-y|^\beta}&\leq \frac{|f(x)-f(y)|}{|x-y|^\alpha}\cdot\delta^{\alpha-\beta}\leq C\,\|f\|_{L^2(\Gamma)}^{1-\theta}\cdot\|f\|_{C^\alpha(\Gamma)}^{\theta}. 
\end{split}
\end{equation*}
Combining these two cases gives \eqref{620} and proves the second line of \eqref{C650}.

Then we apply \eqref{GN2}, \eqref{C406}, Lemmas \ref{41}--\ref{L35},  \eqref{C12241}, and \eqref{CC314} to declare \eqref{618} by
\begin{equation*} 
\begin{split}
&\int_0^T\|\nabla u(\cdot,t)\|_{L^\infty(\Omega_8)}dt\\
&\leq\int_0^T\|\nabla v(\cdot,t)\|_{L^\infty(\Omega_8)}\,dt+\int_0^T\|\nabla w(\cdot,t)\|_{L^\infty(\Omega_8)}\,dt\\
&\leq C\int_0^T\bigg(\|\rho-\bar{\rho}\|_{L^2(\Omega_9)}^{1-\theta'}\cdot\|\rho\|_{C^\alpha(\Omega_9)}^{\theta'}+\|\rho-\bar{\rho}\|_{L^2(B)}^{1/7}\cdot\|\nabla\rho\|_{L^4(B)}^{6/7}\bigg)\,dt\\
&\quad+C\int_0^T\bigg(\|\nabla w\|_{L^\infty(\Omega)}+\|\rho-\bar{\rho}\|_{L^2(\Omega_9)}\bigg)\,dt\leq C,
\end{split}
\end{equation*}
where in the third line, we have applied  \eqref{C406}$_3$ and \eqref{620} in $\Omega^2_8\subset\Omega_9$ to treat $\|\nabla v\|_{L^\infty(\Omega^2_8)}$ with $\theta'=3/(3+2\alpha)$, and applied \eqref{GN2} together with \eqref{C406}$_2$ in $\Omega_3\subset B$ to handle $\|\nabla v\|_{L^\infty(\Omega_3)}$. The proof of Lemma \ref{lm51} is therefore completed.
\end{proof}

With the help of Lemma \ref{lm51}, we can shrink the domain and establish $W^{1,4}$ estimates of $\rho$ in $\Omega_8$ and $\Omega^2_8$. 
\begin{lemma} 
Suppose that \eqref{C501} and Lemmas \ref{41}--\ref{L35} are valid, then there exists a constant $C$ depending only on $\underline{\rho}$, $M_0$, and $\Omega$, such that $\forall t\in[0,T]$,
\begin{equation}\label{614}
\sup_{0\leq\tau\leq t}\|\nabla\rho(\cdot,\tau)\|_{L^4(\Omega_8)}^4+\int_0^t\|\nabla\rho(\cdot,\tau)\|_{L^4(\Omega^2_8)}^4\,d\tau\leq C.
\end{equation}
\end{lemma}
\begin{proof}
Recalling that \eqref{C675} gives
\begin{equation*} 
D_t\rho_{,k}=-\rho_{,k}\,u^s_{,s}-\rho\,u^s_{,sk}.
\end{equation*}
For $\Omega_8=B\cup\Omega^2_8$, we invoke the above equation to  calculate that
\begin{equation}\label{C672}
\begin{split}
&D_t\int_{\Omega^2_8}\big(g^{ij}\,\rho_{,i}\rho_{,j}\big)^2\,d\nu\\
&=\int_{\Omega^2_8}2\,g^{ij}\rho_{,i}\rho_{,j}\,D_t\big(g^{kl}\rho_{,k}\rho_{,l}\big)\,d\nu
+\int_{\Omega^2_8}\big(g^{ij}\,\rho_{,i}\rho_{,j}\big)^2\,D_td\nu\\
&=\int_{\Omega^2_8}g^{ij}\,\rho_{,i}\rho_{,j}\,\big(2D_t(g^{kl})-3g^{kl}u^s_{,s}\big)\,\rho_{,k}\rho_{,l}\,d\nu-4\int_{\Omega^2_8}\rho\, g^{ij}g^{kl}\,u^s_{,sk}\,\rho_{,i}\rho_{,j}\rho_{,l}\,d\nu.\\
\end{split}
\end{equation}
By taking advantage of \eqref{CC208} and \eqref{CC211}, we argue that
\begin{equation}\label{C652}
\begin{split}
&\int_{\Omega^2_8}g^{ij}\,\rho_{,i}\rho_{,j}\big(2D_t(g^{kl})-3g^{kl}u^s_{,s}\big)\,\rho_{,k}\rho_{,l}\,d\nu\leq C\|\nabla u\|_{L^\infty(\Omega^2_8)}\int_{\Omega^2_8}\big(g^{ij}\,\rho_{,i}\rho_{,j}\big)^2\,d\nu.
  \end{split}
\end{equation}
While in view of \eqref{C501}, the principal term in \eqref{C672} is treated by
\begin{equation}\label{C653}
\begin{split}
&\int_{\Omega^2_8}\rho\, g^{ij}g^{kl}\,u^s_{,sk}\,\rho_{,i}\rho_{,j}\rho_{,l}\,d\nu\\
&=\int_{\Omega^2_8}\rho\, g^{ij}g^{kl}\,w^s_{,sk}\,\rho_{,i}\rho_{,j}\rho_{,l}\,d\nu+\int_{\Omega^2_8}\rho\,g^{ij}g^{kl}\,v^s_{,sk}\,\rho_{,i}\rho_{,j}\rho_{,l}\,d\nu\\
&=\int_{\Omega^2_8}\rho\, g^{ij}g^{kl}\,\big(w^s_{,s}+ F_v/2\big)_{,k}\,\rho_{,i}\rho_{,j}\rho_{,l}\,d\nu+\frac{1}{2}\int_{\Omega^2_8}\rho\,\big(g^{ij}\rho_{,i}\rho_{,j}\big)^2\,d\nu,\\
&\leq C\big(\|\nabla F_v\|_{L^4(\Omega^2_8)}+\|\nabla^2 w\|_{L^4(\Omega^2_8)}\big)\cdot\|\nabla\rho\|_{L^4(\Omega^2_8)}^3+\frac{1}{2}\int_{\Omega^2_8}\rho\,\big(g^{ij}\rho_{,i}\rho_{,j}\big)^2\,d\nu,
\end{split}
\end{equation}
where $F_v=2u^s_{,s}-(\rho-\bar{\rho})$.
Thus by substituting \eqref{C652}--\eqref{C653} into \eqref{C672} and recalling that $\rho\geq\underline{\rho}/4$ in $\Omega_{10}$, we 
apply \eqref{355} and \eqref{CC215} to declare that
\begin{equation}\label{C434}
\begin{split}
&\frac{d}{dt}\int_{\Omega^2_8}\big(g^{ij}\,\rho_{,i}\rho_{,j}\big)^2\,d\nu+\int_{\Omega^2_8}\big(g^{ij}\,\rho_{,i}\rho_{,j}\big)^2\,d\nu\\
&\leq C\,\|\nabla u\|_{L^\infty(\Omega^2_8)}\int_{\Omega^2_8}\big(g^{ij}\,\rho_{,i}\rho_{,j}\big)^2\,d\nu+C\bigg(\|\nabla F_v\|_{L^4(\Omega^2_8)}+\|\nabla^2 w\|_{L^4(\Omega^2_8)}\bigg)\,\|\nabla\rho\|_{L^4(\Omega^2_8)}^3\\ 
&\leq C\,\|\nabla u\|_{L^\infty(\Omega^2_8)}\int_{\Omega^2_8}\big(g^{ij}\,\rho_{,i}\rho_{,j}\big)^2\,d\nu+C\,\bigg(
\|\rho-\bar{\rho}\|_{L^2(\Omega)}+\|\rho\dot{u}\|_{L^4(\Omega)}\bigg)\,\|\nabla\rho\|_{L^4(\Omega^2_8)}^3.
\end{split}
\end{equation} 
According to \eqref{C501} and Lemma \ref{21}, it holds that
\begin{equation}\label{CC315} 
\begin{split}
&\int_0^1\|\rho\dot{u}(\cdot,t)\|_{L^4(\Omega)}\,dt
\leq C\int_0^1\sigma^{-3/8}\,\|\rho\dot{u}\|^{1/4}_{L^2(\Omega)}
\bigg(\sigma^{1/2}\|\nabla\dot{u}\|_{L^2(\Omega)}\bigg)^{3/4}dt\leq C,\\
&\int_1^T\|\rho\dot{u}(\cdot,t)\|_{L^4(\Omega)}\,dt
\leq C\int_1^Te^{-\kappa_2t/4}
\bigg(e^{\kappa_2t}\|\rho\dot{u}\|_{L^2(\Omega)}\bigg)^{1/4}\,\|\nabla\dot{u}\|_{L^2(\Omega)}^{3/4}\,dt\leq C.
\end{split}
\end{equation}
Consequently, we integrate \eqref{C434} with respect to $t$ and take advantage of \eqref{C301}, \eqref{618}, \eqref{CC315}, and Gr\"{o}nwall's inequality to deduce that
\begin{equation*}
\sup_{0\leq\tau\leq t}\|\nabla\rho(\cdot,\tau)\|_{L^4(\Omega^2_8)}^4+\int_0^t\|\nabla\rho(\cdot,\tau)\|_{L^4(\Omega^2_8)}^4\,d\tau\leq C,
\end{equation*}
which together with \eqref{C501} gives \eqref{614} and finishes the proof.
\end{proof}
 
Observing that $\Omega^4_6\subset\Omega^3_7 \subset\Omega^2_8,$ the next assertion improves Lemma \ref{21} near the boundary.
In fact, the weight $\sigma$ in \eqref{C305}$_2$ can be dropped, since no vacuum occurs in $\Omega_{10}$ and $\sqrt{\rho}\dot{u}$ admits a proper initial value.  
\begin{lemma}\label{T42}
Suppose that \eqref{C501} and Lemmas \ref{41}--\ref{L35} are valid, then there exists a constant $C$ depending only on $\underline{\rho}$, $M_0$, and $\Omega$, such that $\forall t\in[0,T]$,
\begin{equation}\label{C107}
\sup_{0\leq\tau\leq t}\|\sqrt{\rho}\dot{u}(\cdot,\tau)\|_{L^2(\Omega_7)}^2+\int_0^t\|\nabla\dot{u}(\cdot,\tau)\|_{L^2(\Omega_7)}^2\,d\tau\leq C.
\end{equation}
\end{lemma}
\begin{proof}
Let us first construct a truncation function $\phi(y)$ in Lagrangian coordinates satisfying
\begin{equation*}
\begin{cases}
\phi=1~~\forall y\in\Omega_{7.5},\\
\phi=0~~\forall y\in\Omega\setminus\Omega_8,\\
0\leq\phi\leq 1~~\forall y\in\Omega.
\end{cases}
\end{equation*}
If we go back to Eulerian coordinates,  
the above truncation function corresponds to $\phi\big(X_t^{-1}(x)\big)$, for convenience, we still denote it as $\phi$. Observe that $\partial_t\phi+u\cdot\nabla\phi=0$ in $\Omega$, moreover in view of Lemma \ref{lm29}, \eqref{CC212}, and \eqref{618}, we derive that
\begin{equation}\label{CC308}
\sup_{0\leq\tau\leq t}\|\nabla\phi\|_{L^\infty(\Omega)}\leq C,
\end{equation}
for some constant $C$ depending only on $\underline{\rho}$, $M_0$, and $\Omega$.

Then, by operating $\phi^2\,\dot u^j\left(\frac{\partial}{\partial t}+\div(u\cdot)\right)$ to \eqref{11}$_2$ and summing with respect to $j$, we utilize the fact $\partial_t\phi^2+u\cdot\nabla\phi^2=0$ to obtain that  
\begin{equation}\label{C116}
\begin{split}
&\frac{d}{dt}\int_\Omega\rho|\dot{u}|^2\cdot\phi^2\,dx\\
&=\int_\Omega\big(\Delta u_t+\partial_j(u_j\Delta u)\big)\,\dot{u}\cdot\phi^2\,dx
+\int_\Omega\big(\nabla\mathrm{div}u_t+\partial_j(u_j\nabla\mathrm{div}u)\big)\,\dot{u}\cdot\phi^2\,dx\\
&\quad-\int_\Omega\big(\nabla\rho_t+\partial_j(u_j\nabla\rho)\big)\,\dot{u}\cdot\phi^2\,dx\triangleq I_1+I_2+I_3.
\end{split}
\end{equation}
For $I_1$ and $I_2$, we argue that
\begin{equation*} 
\begin{split}
I_1&=\int_\Omega\big(\Delta u_t+\partial_j(u_j\Delta u)\big)\,\dot{u}\cdot\phi^2\,dx
=-\int_\Omega|\nabla\dot{u}|^2\phi^2\,dx+\mathcal{R}_1,\\
I_2&=\int_\Omega\big(\nabla\mathrm{div}u_t+\partial_j(u_j\nabla\mathrm{div}u)\big)\,\dot{u}\cdot\phi^2\,dx=-\int_\Omega(\mathrm{div}\dot{u})^2\phi^2\,dx+\mathcal{R}_1,\\
\end{split}
\end{equation*}
where the remaining term $\mathcal{R}_1$ is given by
\begin{equation*}
\begin{split}
\mathcal{R}_1=&\int_\Omega\big(\nabla u\cdot\nabla u+u\cdot\nabla^2u\big)\cdot\big(\nabla\dot{u}\cdot\phi^2+\dot{u}\cdot\nabla\phi\cdot\phi\big)\,dx 
+\int_\Omega\nabla\dot u\cdot\dot{u}\cdot\nabla\phi\cdot\phi\,dx,\\
\end{split}
\end{equation*}
where $\nabla^k u$ denote a term involving only the $k^{th}$ order derivatives of $u$, then let us carefully check $\mathcal{R}_1$.
We apply \eqref{C501} and \eqref{CC308} to infer that
\begin{equation}\label{C126}
\begin{split}
&\int_\Omega\nabla u\cdot\nabla u\cdot\big(\nabla\dot{u}\cdot\phi^2+\dot{u}\cdot\nabla\phi\cdot\phi\big)\,dx\\
&=\int_\Omega\nabla^2 u\cdot\nabla u\cdot\dot{u}\cdot\phi^2\,dx
+\int_\Omega\nabla u\cdot\nabla u\cdot\dot{u}\cdot\nabla\phi\cdot\phi\,dx\\
&\leq\|\nabla u\|_{L^\infty(\Omega_8)}\int_\Omega\big(|\nabla^2u|^2+\rho|\dot{u}|^2\big)\cdot\phi^2\,dx
+C\,\|\nabla u\|_{L^\infty(\Omega_8)}\|\nabla u\|_{L^2}^2,\\
\end{split}
\end{equation}
where we have applied the fact $\rho\geq \underline{\rho}/4$ in $\Omega_{10}$. By similar process, it also holds that
\begin{equation}\label{C114}
\begin{split}
&\int_\Omega u\cdot\nabla^2u \cdot\big(\nabla\dot{u}\cdot\phi^2+\dot{u}\cdot\nabla\phi\cdot\phi\big)\,dx+\int_\Omega\nabla\dot u\cdot\dot{u}\cdot\nabla\phi\cdot\phi\,dx\\
&\leq \|u\|_{L^\infty(\Omega)}^2\int_\Omega|\nabla^2u|^2\cdot\phi^2\,dx+\frac{1}{2}\int_\Omega |\nabla\dot{u}|^2\cdot\phi^2\,dx+C\int_\Omega\rho|\dot{u}|^2\,dx,\\
\end{split}
\end{equation}
where we have made use of \eqref{C501}, \eqref{CC308} and the fact $\rho\geq\underline{\rho}/4$ in $\Omega_{10}$.
Moreover, according to Lemma \ref{lm42}, \eqref{C406}$_2$, \eqref{C301}, \eqref{C305}, \eqref{CC307}, and \eqref{614}, we check that 
\begin{equation*}
\begin{split}
&\int_\Omega|\nabla^2u|^2\cdot\phi^2\,dx\\
&\leq \int_\Omega|\nabla^2w|^2\cdot\phi^2\,dx
+\int_\Omega|\nabla^2v|^2\cdot\phi^2\,dx\\
&\leq C\big(\int_\Omega\rho|\dot{u}|^2\cdot\phi^2\,dx
+\|\nabla w\|^2_{L^2(\Omega)}
+\|\rho-\bar{\rho}\|_{L^2(\Omega)}^2+\|\nabla\rho\|^2_{L^2(\Omega_8)}\big)\\
&\leq C\big(\int_\Omega\rho|\dot{u}|^2\cdot\phi^2\,dx+1\big),
\end{split}
\end{equation*}
which together with \eqref{C126} and \eqref{C114} gives that
\begin{equation}\label{CC309}
\begin{split}
&I_1+I_2+\frac{1}{2}\int_\Omega|\nabla\dot{u}|^2\cdot\phi\,dx\\
&\leq C\big(\|\nabla u\|_{L^\infty(\Omega_8)}+\|u\|_{L^\infty(\Omega)}^2\big)\big(\int_\Omega\rho|\dot{u}|^2\cdot\phi^2\,dx+1\big)
+C\,\|\sqrt{\rho}\dot{u}\|_{L^2(\Omega)}^2.\\
\end{split}
\end{equation}
While for $I_3$, we utilize \eqref{614} to derive that
\begin{equation}\label{CC310}
\begin{split}
I_3&=-\int_\Omega\big(\nabla\rho_t+\partial_j(u_j\,\nabla\rho)\big)\,\dot{u}\cdot\phi^2\,dx\\
&=\int_\Omega\big(u\cdot\nabla\rho+\rho\cdot\nabla u\big)\cdot\big(\nabla\dot{u}\cdot\phi^2+\dot{u}\cdot\nabla\phi\cdot\phi\big)\,dx\\
&\leq C\big(\|u\|_{L^4(\Omega)}^2\|\nabla\rho\|_{L^4(\Omega_8)}^2+\|\nabla u\|_{L^2(\Omega)}^2+\|\sqrt{\rho}\dot{u}\|_{L^2(\Omega)}^2\big)+\frac{1}{4}\int_\Omega |\nabla\dot{u}|^2\cdot\phi^2\,dx\\
&\leq C\big(\|\nabla u\|_{L^2(\Omega)}^2+\|\sqrt{\rho}\dot{u}\|_{L^2(\Omega)}^2\big)+\frac{1}{4}\int_\Omega |\nabla\dot{u}|^2\cdot\phi^2\,dx.
\end{split}
\end{equation}
Collecting \eqref{C116}, \eqref{CC309}, and \eqref{CC310} yields that
\begin{equation}\label{C117}
\begin{split}
&\frac{d}{dt}\int_\Omega\rho|\dot{u}|^2\cdot\phi^2\,dx+\int_\Omega|\nabla\dot{u}|^2\cdot\phi^2\,dx\\
&\leq C(\|\nabla u\|_{L^\infty(\Omega_8)}+\|u\|_{L^\infty(\Omega)}^2\big)\big(\int_\Omega\rho|\dot{u}|^2\cdot\phi^2\,dx+1\big)+C\big(\|\nabla u\|_{L^2(\Omega)}^2+\|\sqrt{\rho}\dot{u}\|_{L^2(\Omega)}^2\big).
\end{split}
\end{equation}
Moreover, we apply Sobolev's inequality, \eqref{355}, \eqref{C501}, and Lemmas \ref{41}--\ref{21} to obtain
\begin{equation}\label{C115}
\begin{split}
&\int_0^T\|u(\cdot,t)\|_{L^\infty}^2\,
dt\\
&\leq C\int_0^T\big(\|w(\cdot,t)\|_{H^{2,2}}^2+\|v(\cdot,t)\|_{W^{1,4}}^2\big)\,dt\\
&\leq C\int_0^T\big(\|\sqrt{\rho}\dot{u}(\cdot,t)\|_{L^2}^2+\|(\rho-\bar{\rho})(\cdot,t)\|_{L^4}^2\big)dt\leq C.
\end{split}\end{equation}

By virtue of Lemmas \ref{41}--\ref{21}, \eqref{618}, and \eqref{C115},  we utilize Gr\"{o}nwall's inequality in \eqref{C117} to deduce that $\forall t\in[0,T]$
\begin{equation*}
\sup_{0\leq\tau\leq t}\int_\Omega\rho|\dot{u}|^2\cdot\phi^2\,dx+\int_0^t\int_\Omega|\nabla\dot{u}|^2\cdot\phi^2\,dxd\tau
\leq C,
\end{equation*}
where the initial value of $\sqrt{\rho}\dot{u}\cdot\phi$ can be taken as
$$\sqrt{\rho}\dot{u}\cdot\phi\bigg|_{t=0}
=\frac{\phi(x,0)}{\sqrt{\rho_0}}
\big(\Delta u_0+\nabla\mathrm{div}u_0-\nabla\rho_0\big),$$
due to the fact $\rho_0\geq\underline{\rho}$ in $\Omega_{10}$. The proof of Lemma \ref{T42} is therefore completed.
\end{proof}
 
Now, we can complete the local $H^2$ estimates of $\rho$ in $\Omega^4_6$. It gives proper controls  over the domain linking the inner and boundary parts of $\Omega$.
\begin{lemma}\label{lm53}
Suppose that \eqref{C501} and Lemmas \ref{41}--\ref{L35} are valid, then there exists a constant $C$ depending only on $\underline{\rho}$, $M_0$, and $\Omega$, such that $\forall t\in[0,T]$, 
\begin{equation}\label{C12243}
\sup_{0\leq\tau\leq t}\|\nabla^2\rho(\cdot,\tau)\|_{L^2(\Omega^4_6)}^2+\int_0^t\|\nabla^2\rho(\cdot,\tau)\|_{L^2(\Omega^4_6)}^2\,d\tau\leq C.
\end{equation}
\end{lemma}
\begin{proof} 
By taking advantage of Lemma \ref{411}, \eqref{C665}, and \eqref{CC211}, we deduce that 
\begin{equation}\label{C673}
\begin{split}
&\frac{d}{dt}\bigg(\frac{1}{2}\int_{\Omega^4_6}g^{ij}g^{kl}\,\rho_{,ik}\rho_{,jl}\,d\nu\bigg)\\
&=\int_{\Omega^4_6}D_t (g^{ij})g^{kl}\rho_{,ik}\rho_{,jl}\,d\nu +\int_{\Omega^4_6}g^{ij}g^{kl}D_t (\rho_{,ik})\rho_{,jl}\,d\nu
+\frac{1}{2}\int_{\Omega^4_6}g^{ij}g^{kl}\rho_{,ik}\rho_{,jl}\cdot D_td\nu\\
&=\int_{\Omega^4_6}D_t (g^{ij})g^{kl}\,\rho_{,ik}\rho_{,jl}\,d\nu-\int_{\Omega^4_6}g^{ij}g^{kl}\big((\rho_{,s}u^s_{,ki}+2\rho_{,i} u^s_{,sk})
+\rho_{,ik}u^s_{,s}/2+\rho u^s_{,sik}\big)\rho_{,jl}\,d\nu\\
&\triangleq J_1+J_2+J_3,\\
\end{split}
\end{equation}
where the remaining terms $J_i$ are given by
\begin{equation*}
\begin{split}
&J_1\triangleq\int_{\Omega^4_6}g^{kl}\left(D_t (g^{ij})
-g^{ij}u^s_{,s}/2\right)\rho_{,ik}\rho_{,jl}\,d\nu,\\
&J_2\triangleq-\int_{\Omega^4_6}g^{ij}g^{kl}\,\big(u^s_{,ki}\rho_{,s}+2\rho_{,i} u^s_{,sk}\big)\rho_{,jl}\,d\nu,\\
&J_3\triangleq-\int_{\Omega^4_6}\rho\,g^{ij}g^{kl}\, u^s_{,sik}\rho_{,jl}\,d\nu.\\
\end{split}
\end{equation*}
According to Lemma \ref{411}, $J_1$ and $J_2$ are bounded by
\begin{equation}\label{C657}
\begin{split}
|J_1|&\leq C\,\|\nabla u\|_{L^\infty(\Omega^4_6)}\int_{\Omega^4_6}|\nabla^2\rho|^2\,d\nu,\\
|J_2|&\leq C\,\|\nabla^2u\|_{L^4(\Omega^4_6)}\,\|\nabla\rho\|_{L^4(\Omega^4_6)}\,\|\nabla^2\rho\|_{L^2(\Omega^4_6)}.\\
\end{split}
\end{equation}
Then we turn to the principal term $J_3$ in \eqref{C673}.
\begin{equation}\label{C658}
\begin{split}
J_3&=-\int_{\Omega^4_6}\rho\,g^{ij}g^{kl}\,u^s_{,sik}\rho_{,jl}\,d\nu\\
&=-\int_{\Omega^4_6}\rho\,g^{ij}g^{kl}\big(w^s_{,s}+F_v\big)_{,ik}\rho_{,jl}\,d\nu-\frac{1}{2}
\int_{\Omega^4_6}\rho\,g^{ij}g^{kl}\,\rho_{,ik}\rho_{,jl}\,d\nu\\
&\leq C\bigg(\|\nabla^3 w\|_{L^2(\Omega^4_6)}+\|F_v\|_{L^2(\Omega^4_6)}\bigg)\|\nabla^2\rho\|_{L^2(\Omega^4_6)}-\frac{1}{2}
\int_{\Omega^4_6}\rho\,g^{ij}g^{kl}\,\rho_{,ik}\rho_{,jl}\,d\nu.
\end{split}
\end{equation}
Thus, in view of $\rho\geq\underline{\rho}/4$ in $\Omega_{10}$, we substitute \eqref{C657} and \eqref{C658} into \eqref{C673} and arrive at
\begin{equation}\label{C661}
\begin{split}
&\frac{d}{dt}\int_{\Omega^4_6}|\nabla^2\rho|^2\,d\nu+\int_{\Omega^4_6}|\nabla^2\rho|^2\,d\nu \leq C\bigg(\|\nabla u\|_{L^{\infty}(\Omega^4_6)}\cdot\int_{\Omega^4_6}|\nabla^2\rho|^2\,d\nu+\mathcal{R}_2\bigg),
\end{split}
\end{equation}
where the remaining term $\mathcal{R}_2$ is given by
$$\mathcal{R}_2(t)\triangleq \|\nabla^2u\|_{L^4(\Omega^4_6)}^2\,\|\nabla\rho\|_{L^4(\Omega^4_6)}^2+\|\nabla^3 w\|_{L^2(\Omega^4_6)}^2
+\|F_v\|_{L^2(\Omega^4_6)}^2.$$

Let us utilize \eqref{GN2} together with \eqref{C406} in $\Omega^4_6\subset\Omega^{7/2}_{13/2}\subset\Omega^3_7$ to deduce that
\begin{equation}\label{CC317}
\begin{split} 
&\|\nabla^2u\|_{L^4(\Omega^4_6)}\\
&\leq C\bigg(\|\nabla^2w\|_{L^4(\Omega^4_6)}
+\|\nabla^2v\|_{L^4(\Omega^4_6)}\bigg)\\
&\leq C\bigg(\|\dot{u}\|_{H^1(\Omega^3_7)}+\|\rho-\bar{\rho}\|_{W^{1,4}(\Omega^3_7)}
+\|\nabla u\|_{L^2(\Omega)}
\bigg),\\
&\|\nabla^3w\|_{L^2(\Omega^4_6)}\\
&\leq C\bigg(\|\rho\,\nabla\dot{u}\|_{L^2(\Omega^{7/2}_{13/2})}+\|\nabla\rho\cdot\dot{u}\|_{L^2(\Omega^{7/2}_{13/2})}+\|\nabla w\|_{L^2(\Omega)}\bigg)\\
&\leq C\bigg(\|\dot{u}\|_{H^1(\Omega^3_7)}\,\|\rho-\bar{\rho}\|_{W^{1,4}(\Omega^3_7)}+
\|\nabla u\|_{L^2(\Omega)}+\|\rho-\bar{\rho}\|_{L^2(\Omega)}\bigg),\\
\end{split}
\end{equation}
which along with \eqref{CC215}, \eqref{614}, \eqref{C107}, and Lemmas \ref{41}--\ref{L35} yields that 
\begin{equation}\label{CC316}
\int_0^T\mathcal{R}_2(t)\,dt\leq C.
\end{equation}
Therefore,
by collecting \eqref{618}, \eqref{C661}, and \eqref{CC316}, we apply Gr\"{o}nwall's inequality to derive that
\begin{equation*}
\sup_{0\leq t\leq T}\int_{\Omega^4_6}|\nabla^2\rho|^2\,d\nu+\int_0^T\int_{\Omega^4_6}|\nabla^2\rho|^2\,d\nu dt\leq C,
\end{equation*}
which is \eqref{C12243} and finishes the proof of Lemma \ref{lm53}.
\end{proof}

\subsection{The estimates in the boundary domain $B$}\label{SS42}
\quad In this section, we will establish the local $H^2$ energy estimates in the boundary domain $B$, and the main result is given by the next lemma.
 
\begin{lemma}\label{lm54}
Suppose that \eqref{C501}, Lemmas \ref{41}--\ref{L35}, and Lemma \ref{lm53} are valid, then there is a constant $C$ depending only on $\underline{\rho}$, $M_0$, and $\Omega$, such that $\forall t\in[0,T]$, 
\begin{equation}\label{687}
\sup_{0\leq\tau\leq t}\|\nabla^2\rho(\cdot,\tau)\|_{L^2(\Omega_6)}\leq C. 
\end{equation}
\end{lemma}

According to the decomposition \eqref{C43}, we argue that
\begin{equation}\label{CC435}
\begin{split}
|\nabla^2\rho|^2 
&=g^{ij}g^{kl}\,\rho_{,ik}\rho_{,jl}\\
&=(\gamma^{ij}+N^iN^j)(\gamma^{kl}+N^kN^l)\,\rho_{,ik}\rho_{,jl}\\
&=\gamma^{ij}\gamma^{kl}\,\rho_{,ik}\rho_{,jl}+2\gamma^{ij}N^kN^l\,\rho_{,ik}\rho_{,jl}+(N^iN^k\,\rho_{,ik})^2.
\end{split}
\end{equation}
The above equation consists of the tangential, mixed, and normal derivatives. We will
obtain estimates on each terms in following lemmas.

To begin with, let us carry out some calculations  common for each lemma. We introduce a truncation function $\psi(y)\in C^\infty(\Omega)$ satisfying
\begin{equation*}
\begin{cases}
\psi=1~~\forall y\in\Omega_{4.5},\\
\psi=0~~\forall y\in\Omega\setminus\Omega_5,\\
0\leq\psi\leq 1,~\sum_i\big(\partial\psi/\partial y^i\big)^2\leq 2~~\forall y\in\Omega.
\end{cases}
\end{equation*}
According to Lemma \ref{lm29}, \eqref{CC212}, and \eqref{618}, we have $D_t\psi=0,\,|\nabla\psi|\leq C$ for some $C$ depending only on $\underline{\rho}$, $M_0$, and $\Omega$.  Suppose that $\mathbf{\Phi}$
and $\mathbf{\Psi}$ are two symmetric tensors of type $(2,0)$, then in view of \eqref{C665} and \eqref{C673}, it holds that
\begin{equation}\label{CC401}
\begin{split}
D_t\bigg(\frac{1}{2}\int_\Omega\mathbf{\Phi}^{ij}\mathbf{\Psi}^{kl}\,\rho_{,ik}\rho_{,jl}\,\psi^2\,d\nu\bigg)
&=K_1+K_2+K_3,\\
\end{split}
\end{equation}
where  $K_i$ are given by
\begin{equation}\label{CC402}
\begin{split}
&K_1\triangleq\frac{1}{2}\int_\Omega \big(D_t(\mathbf{\Phi}^{ij} \mathbf{\Psi}^{kl})-\mathbf{\Phi}^{ij} \mathbf{\Psi}^{kl} u^s_{,s}/2\big)\,\rho_{,ik}\rho_{,jl}\,\psi^2\,d\nu,\\
&K_2\triangleq-\int_{\Omega}\mathbf{\Phi}^{ij}\mathbf{\Psi}^{kl}\,\big(u^s_{,ki}\rho_{,s}+2\rho_{,i} u^s_{,sk}\big)\rho_{,jl}\,\psi^2\,d\nu,\\
&K_3\triangleq-\int_{\Omega}\rho\,\mathbf{\Phi}^{ij}\mathbf{\Psi}^{kl}\, u^s_{,sik}\rho_{,jl}\,\psi^2\,d\nu.\\
\end{split}
\end{equation}
If $\mathbf{\Phi}$ and $\mathbf{\Psi}$ are taken as $\gamma$ or $N\otimes N$, in view of Lemma \ref{411} and \eqref{CC211}, we argue that
\begin{equation}\label{CC403}
\begin{split}
|K_1|&\leq C\,\|\nabla u\|_{L^\infty(B)}\int_{\Omega}|\nabla^2\rho|^2\,\psi^2\,d\nu,\\
|K_2|&\leq C\,\|\nabla^2u\|_{L^4(B)}\,\|\nabla\rho\|_{L^4(B)}\,\bigg(\int_\Omega|\nabla^2\rho|^2\,\psi^2\,d\nu\bigg)^{1/2},\\
\end{split}
\end{equation}
while $K_3$ is treated by
\begin{equation}\label{CC404}
\begin{split}
K_3&=-\int_{\Omega}\rho\,\mathbf{\Phi}^{ij}\mathbf{\Psi}^{kl}\, u^s_{,sik}\rho_{,jl}\,\psi^2\,d\nu,\\
&=-\int_{\Omega}\rho\,\mathbf{\Phi}^{ij}\mathbf{\Psi}^{kl}\, w^s_{,sik}\rho_{,jl}\,\psi^2\,d\nu-\int_{\Omega}\rho\,\mathbf{\Phi}^{ij}\mathbf{\Psi}^{kl}\, v^s_{,sik}\rho_{,jl}\,\psi^2\,d\nu\\
&\leq C\,\|\nabla^3w\|_{L^2(B)}\bigg(\int_\Omega|\nabla^2\rho|^2\,\psi^2\,d\nu\bigg)^{1/2}-\int_{\Omega}\rho\,\mathbf{\Phi}^{ij}\mathbf{\Psi}^{kl}\, v^s_{,sik}\rho_{,jl}\,\psi^2\,d\nu.
\end{split}
\end{equation}
Substituting \eqref{CC402}--\eqref{CC404} into \eqref{CC401}, we arrive at
\begin{equation}\label{CC405}
\begin{split}
&D_t\bigg(\frac{1}{2}\int_\Omega\mathbf{\Phi}^{ij}\mathbf{\Psi}^{kl}\,\rho_{,ik}\rho_{,jl}\,\psi^2\,d\nu\bigg)\\
&\leq-\int_{\Omega}\rho\,\mathbf{\Phi}^{ij}\mathbf{\Psi}^{kl}\, v^s_{,sik}\rho_{,jl}\,\psi^2\,d\nu+\|\nabla u\|_{L^\infty(B)}\int_{\Omega}|\nabla^2\rho|^2\,\psi^2\,d\nu+\mathcal{R}_3.
\end{split}
\end{equation}
The first part in the above line is the \textbf{principal term}, and the remainder $\mathcal{R}_3$ is given by
\begin{equation*} 
\begin{split}
\mathcal{R}_3(t)&=\big(\|\nabla^2u\|_{L^4(B)}\,\|\nabla\rho\|_{L^4(B)}+\|\nabla^3w\|_{L^2(B)}\big)\,\bigg(\int_\Omega|\nabla^2\rho|^2\,\psi^2\,d\nu\bigg)^{1/2}.\\
\end{split}
\end{equation*}

In addition, let us provide some local estimates which will be repeatedly used. By virtue of \eqref{GN2}, \eqref{355}, \eqref{C406},  \eqref{C314}, \eqref{C301}, \eqref{CC307}, and  \eqref{C12243}, we infer that
\begin{equation}\label{CC408}
\begin{split}
&\|\rho-\bar{\rho}\|_{H^2(B)}+\|v\|_{H^3(B)}\\
&\leq C\,\bigg(\|\rho-\bar{\rho}\|_{H^2(\Omega_6)}+\|\rho-\bar{\rho}\|_{L^2(\Omega_6)}\bigg)\\
&\leq C\,\bigg(\int_\Omega|\nabla^2\rho|^2\,\psi^2\,d\nu+1\bigg)^{1/2},\\
&\|\rho-\bar{\rho}\|_{H^1(B)}^2+\|v\|_{H^2(B)}^2\\
&\leq C\,\bigg(\|\rho-\bar{\rho}\|_{L^2(\Omega_{13/2})}\|\nabla^2\rho\|_{L^2(\Omega_{13/2})}+\|\rho-\bar{\rho}\|_{L^2(\Omega_{13/2})}^2\bigg)\\
&\leq C\,\|\rho-\bar{\rho}\|_{L^2(\Omega)}\,  \left(\int_\Omega|\nabla^2\rho|^2\,\psi^2\,d\nu+1\right)^{1/2},
\end{split}
\end{equation}
As a direct corollary of \eqref{CC408}, we make use of \eqref{GN2}, \eqref{C501}, and \eqref{CC307} to argue that
\begin{equation*}
\begin{split}
\|\rho-\bar{\rho}\|_{W^{1,4}(\Omega_{11/2})}&\leq C\,\bigg(\|\rho-\bar{\rho}\|_{L^2(\Omega_6)}^{1/8}\|\nabla^2\rho\|_{L^2(\Omega_6)}^{7/8}+\|\rho-\bar{\rho}\|_{L^4(\Omega_6)}\bigg)\\
&\leq C\,\|\rho-\bar{\rho}\|_{L^2(\Omega)}^{1/8}\bigg(\int_\Omega|\nabla^2\rho|^2\,\psi^2\,d\nu+1\bigg)^{7/16},
\end{split}
\end{equation*}
which along with \eqref{CC317} provides proper estimates on the remaining term $\mathcal{R}_3$,
\begin{equation}\label{CC406}
\begin{split}
\mathcal{R}_3(t) 
&\leq C\,\bigg(\|\dot{u}\|_{H^1(\Omega_{11/2})}+\|\rho-\bar{\rho}\|_{W^{1,4}(\Omega_{11/2})}+\|\nabla u\|_{L^2(\Omega)}\bigg)\,\bigg(\int_\Omega|\nabla^2\rho|^2\,\psi^2\,d\nu\bigg)^{1/2}\\
&\leq C\,\bigg(\|\nabla\dot{u}\|_{L^2(\Omega_{7})}+\|\sqrt{\rho}\dot{u}\|_{L^2(\Omega)}+\|\nabla u\|_{L^2(\Omega)}\bigg)\,\bigg(\int_\Omega|\nabla^2\rho|^2\,\psi^2\,d\nu\bigg)^{1/2}\\
&\quad+C\,\|\rho-\bar{\rho}\|_{L^2(\Omega)}^{1/8}\bigg(\int_\Omega|\nabla^2\rho|^2\,\psi^2\,d\nu+1\bigg)^{15/16}.
\end{split}
\end{equation}

After these preparations, let us turn to the estimates on the tangential directions.
\begin{lemma}
Suppose that \eqref{C501}, Lemmas \ref{41}--\ref{L35}, and Lemma \ref{lm53} are valid, then there is a constant $C$ depending only on $\underline{\rho}$, $M_0$, and $\Omega$, such that $\forall t\in[0,T]$, 
\begin{equation}\label{CC409}
\begin{split}
&\frac{d}{dt}\int\gamma^{ij}\gamma^{kl}\rho_{,ik}\rho_{,jl}\psi^2d\nu
+\int_\Omega\gamma^{ij}\gamma^{kl}\,\rho_{,ik}\rho_{,jl}\,\psi^2\,d\nu\\
&\leq C\,\big(\|\rho-\bar{\rho}\|_{L^\infty(B)}
^{1/2}+\|\nabla u\|_{L^\infty(B)}\big)\bigg(\int|\nabla^2\rho|^2\psi^2d\nu+1\bigg)
+\mathcal{R}_3.
\end{split}
\end{equation}
\end{lemma}
\begin{proof}
Taking $\mathbf{\Phi}=\mathbf{\Psi}=\gamma$ in the \textbf{principal term} of \eqref{CC405}, we check that
\begin{equation}\label{C611}
\begin{split}&
\int_\Omega\rho\,\gamma^{ij}\gamma^{kl}\,v^s_{,sik}\rho_{,jl}\,\psi^2\,d\nu\\
&=\int_\Omega(\rho-\bar{\rho})\,\gamma^{ij}\gamma^{kl}\, v^s_{,sik}\rho_{,jl}\,\psi^2\,d\nu+\bar{\rho}\int_\Omega\gamma^{ij}\gamma^{kl}\, v^s_{,sik}\rho_{,jl}\,\psi^2\,d\nu\\
&=L_1+\bar{\rho}\,L_2,
\end{split}
\end{equation}
where in view of \eqref{CC408}, $L_1$ is directly bounded by
\begin{equation}\label{C0613}
\begin{split}
|L_1|
&\leq C\,\|\rho-\bar{\rho}\|_{L^\infty(B)}\bigg(\|v\|_{H^3(B)}+\|\rho-\bar{\rho}\|_{H^2(B)}\bigg)\\
&\leq C\,\|\rho-\bar{\rho}\|_{L^\infty(B)}\,\bigg(\int_\Omega|\nabla^2\rho|^2\,\psi^2\,d\nu+1\bigg).
\end{split}
\end{equation}
Inspired by \cite{MN2}, the estimates on $L_2$ is divided into several steps.

\textit{Step 1. There is a constant $C$ depending only on $\underline{\rho}$, $M_0$, and $\Omega$, such that} 
\begin{equation}\label{CC410}
\begin{split}
-L_2\leq-\int_\Omega \gamma^{ij}\gamma^{kl}g^{mn}v_{s,mik}v^s_{,njl}\,\psi^2\,d\nu+C\,\|\rho-\bar{\rho}\|_{L^2(\Omega)}^{1/2}\,\left(\int_\Omega|\nabla^2\rho|^2\,\psi^2\,d\nu+1\right)^{3/4}.
\end{split}
\end{equation}
Let us take the second order covariant derivatives on \eqref{62} and apply \eqref{C401} to argue that
\begin{equation*}
g^{mn}\,v_{s,mnjl}+v^m_{,msjl}=\rho_{,sjl} .\\
\end{equation*}
Multiplying the above equation by $\gamma^{ij}\gamma^{kl}\,v^s_{,ik}$ leads to
\begin{equation}\label{CC433}
\gamma^{ij}\gamma^{kl}\,\big(v^s_{,ik}\,\rho_{,sjl}-(g^{mn}\,v_{s,mnjl}\, v^s_{,ik} +v^m_{,msjl}\,v^s_{,ik})\big)=0,
\end{equation}
thus we can rewrite $L_2$ as 
\begin{equation}\label{CC432}
\begin{split}
L_2&=\int_\Omega\gamma^{ij}\gamma^{kl}\,\big(\big(v^s_{,sik}\,\rho_{,jl}+v^s_{,ik}\,\rho_{,sjl}\big)-(g^{mn}v_{s,mnjl}\,v^s_{,ik}
+v^m_{,msjl}\,v^s_{,ik})\big)\,\psi^2\,d\nu.\\
\end{split}
\end{equation}

For the first term in \eqref{CC432}, we calculate that
\begin{equation}\label{C615}
\begin{split}
&\int_\Omega\gamma^{ij}\gamma^{kl}\,\big(v^s_{,sik}\,\rho_{,jl}+v^s_{,ik}\,\rho_{,sjl}\big)\,\psi^2\,d\nu\\
&=\int_\Omega(\gamma^{ij}\gamma^{kl}v^s_{,ik}\,\rho_{,jl}\psi^2)_{,s}d\nu
-2\int_\Omega\big(\gamma^{ij}_{,s}\gamma^{kl}\,v^s_{,ik}\,\rho_{,jl}\,\psi^2+\gamma^{ij}\gamma^{kl}\,v^s_{,ik}\,\rho_{,jl}\,\psi\,\psi_{,s}\big)\,d\nu.
\end{split}
\end{equation}
According to Lemma \ref{411} and \eqref{CC408}, we argue that the second term of \eqref{C615} is bounded by
\begin{equation}\label{C616}
\begin{split}
&\left|\int_\Omega\big(\gamma^{ij}_{,s}\gamma^{kl}\,v^s_{,ik}\,\rho_{,jl}\,\psi^2+\gamma^{ij}\gamma^{kl}\,v^s_{,ik}\,\rho_{,jl}\,\psi\,\psi_{,s}\big)\,d\nu\right|\\
&\leq C\|v\|_{H^2(B)}\|\rho-\bar{\rho}\|_{H^2(B)}\\
&\leq C\,\|\rho-\bar{\rho}\|_{L^2(\Omega)}^{1/2}\, \left(\int_\Omega|\nabla^2\rho|^2\,\psi^2\,d\nu+1\right)^{3/4}.\\
\end{split}
\end{equation}
Moreover, we apply \eqref{CC216} in $\Omega$ to handle the first term of \eqref{C615} by
\begin{equation}\label{C96}
\begin{split}
&\int_\Omega(\gamma^{ij}\gamma^{kl}\,v^s_{,ik}\rho_{,jl}\,\psi^2)_{,s}d\nu\\
&=\int_{\partial\Omega} N_s\big(\gamma^{ij}\gamma^{kl}\,v^s_{,ik}\big)\,\rho_{,jl}\,d\nu
=-\int_{\partial\Omega}N_s\big(\gamma^{ij}\gamma^{kl}_{,i}\,v^s_{,k}\big)\,\rho_{,jl}\,d\nu,
\end{split}
\end{equation} 
where in the last equality we have taken advantage of \eqref{CC218} to deduce that
\begin{equation}\label{CC413}
\gamma^{ij}\gamma^{kl}v^s_{,ik}=\gamma^{ij}\big(\gamma^{kl}v^s_{,k}\big)_{,i}-\gamma^{ij}\gamma^{kl}_{,i}\,v^s_{,k}=-\gamma^{ij}\gamma^{kl}_{,i}\,v^s_{,k},
\end{equation}
which is due to the fact $v=0$ on $\partial\Omega$.
In addition, we further change \eqref{C96} into
\begin{equation*} 
\begin{split}
&\int_{\partial\Omega}N_s\big(\gamma^{ij}\gamma^{kl}_{,i}\,v^s_{,k}\big)\,\rho_{,jl}\,d\nu\\
&=\int_{\partial\Omega}\gamma^{ij}\left(N_s\gamma^{kl}_{,i}v^s_{,k}\,\rho_{,l}\,\psi^2\right)_{,j}\,d\nu-\int_{\partial\Omega}\gamma^{ij}\left(N_{s}\gamma^{kl}_{,i}v^s_{,k}\right)_{,j}\rho_{,l}d\nu. 
\end{split}
\end{equation*}
Now let us utilize \eqref{C409}, \eqref{CC217}, \eqref{C401}, and \eqref{CC408} to infer that
\begin{equation}\label{CC407}
\begin{split}
&\bigg|\int_{\partial\Omega}\gamma^{ij}\left(N_s\gamma^{kl}_{,i}v^s_{,k}\,\rho_{,l}\,\psi^2\right)_{,j}\,d\nu-\int_{\partial\Omega}\gamma^{ij}\left(N_{s}\gamma^{kl}_{,i}v^s_{,k}\right)_{,j}\rho_{,l}d\nu\bigg|\\
&\leq C\,\bigg(\|\rho-\bar{\rho}\|_{H^2(B)}\|v\|_{H^2(B)}+\|\rho-\bar{\rho}\|_{H^1(B)}\|v\|_{H^3(B)}\bigg)\\
&\leq C\,\|\rho-\bar{\rho}\|_{L^2(\Omega)}^{1/2}\,  \left(\int_\Omega|\nabla^2\rho|^2\,\psi^2\,d\nu+1\right)^{3/4}.
\end{split}
\end{equation}
Combining \eqref{C615}, \eqref{C616}, and \eqref{CC407} yields that
\begin{equation}\label{CC416}
\begin{split}
\int_\Omega\gamma^{ij}\gamma^{kl}\,\bigg(v^s_{,sik}\,\rho_{,jl}+v^s_{,ik}\,\rho_{,sjl}\bigg)\,\psi^2\,d\nu&\leq C\,\|\rho-\bar{\rho}\|_{L^2(\Omega)}^{1/2}\,\left(\int_\Omega|\nabla^2\rho|^2\,\psi^2\,d\nu+1\right)^{3/4}.
\end{split}
\end{equation}
 
Next we turn to the second term in \eqref{CC432} which is given by
$$\int_\Omega\gamma^{ij}\gamma^{kl}\,\big(g^{mn}v_{s,mnjl}\,v^s_{,ik}
+v^m_{,msjl}\,v^s_{,ik}\big)\,\psi^2\,d\nu.$$ 
Compared with \eqref{C615}, there is an additional good term involved. Let us check that
\begin{equation}\label{CC411}
\begin{split}
&\int_\Omega g^{mn}\gamma^{ij}\gamma^{kl}\,v_{s,mnjl}\,v^s_{,ik}\,\psi^2\,d\nu+\int_\Omega\gamma^{ij}\gamma^{kl}g^{mn}\,v_{s,mjl}\,v^s_{,nik}\,\psi^2\,d\nu\\
&=\int_\Omega g^{mn}\big(\gamma^{ij}\gamma^{kl}\,v_{s,mjl}\,v^s_{,ik}\,\psi^2\big)_{,n}\,d\nu-\int_\Omega g^{mn}\big(\gamma^{ij}\gamma^{kl}\,\psi^2\big)_{,n}v_{s,mjl}\,v^s_{,ik}\,d\nu.
\end{split}
\end{equation}
We apply \eqref{C401} and \eqref{CC408} to bound the second term of \eqref{CC411} by
\begin{equation}\label{CC412}
\begin{split}
&\left|\int_\Omega g^{mn}\big(\gamma^{ij}\gamma^{kl}\,\psi^2\big)_{,n}v_{s,mjl}\,v^s_{,ik}\,d\nu\right|\\
&\leq C\|v\|_{H^3(B)}\|v\|_{H^2(B)}\\
&\leq C\,\|\rho-\bar{\rho}\|_{L^2(\Omega)}^{1/2}\, \left(\int_\Omega|\nabla^2\rho|^2\,\psi^2\,d\nu+1\right)^{3/4}.
\end{split}
\end{equation}
While in view of \eqref{CC216}, the first term of \eqref{CC411} satisfies that
\begin{equation}\label{CC414}
\begin{split}
&\int_\Omega g^{mn}\big(\gamma^{ij}\gamma^{kl}\,v_{s,mjl}\,v^s_{,ik}\,\psi^2\big)_{,n}\,d\nu\\
&=\int_{\partial\Omega}N_n\,g^{mn}\big(\gamma^{ij}\gamma^{kl}\,v^s_{,ik}\big)\,v_{s,mjl}\,d\nu=-\int_{\partial\Omega}N_n\,g^{mn}\big(\gamma^{ij}\gamma^{kl}_{,i}\,v^s_{,k}\big)\,v_{s,mjl}\,d\nu,
\end{split}
\end{equation}
where we have applied \eqref{CC413} in the last line. In addition, we rewrite \eqref{CC414} as
\begin{equation*}
\begin{split}
&\int_{\partial\Omega}N_n\,g^{mn}\big(\gamma^{ij}\gamma^{kl}_{,i}\,v^s_{,k}\big)\,v_{s,mjl}\,d\nu
\\
&=\int_{\partial\Omega}\gamma^{ij}\left(N_n\,g^{mn}\gamma^{kl}_{,i}\,v^s_{,k}\,v_{s,ml}\,\psi^2\right)_{,j}\,d\nu-\int_{\partial\Omega}\gamma^{ij}\left(N_n\,g^{mn}\gamma^{kl}_{,i}\,v^s_{,k}\right)_{,j}\,v_{s,ml}\,d\nu.
\end{split}
\end{equation*}
Therefore, let us utilize  \eqref{C409}, \eqref{CC217}, \eqref{C401}, and \eqref{CC408} to infer that
\begin{equation}\label{CC415}
\begin{split}
&\bigg|\int_{\partial\Omega}\gamma^{ij}\left(N_n\,g^{mn}\gamma^{kl}_{,i}\,v^s_{,k}\,v_{s,ml}\right)_{,j}\,d\nu-\int_{\partial\Omega}\gamma^{ij}\left(N_n\,g^{mn}\gamma^{kl}_{,i}\,v^s_{,k}\right)_{,j}\,v_{s,ml}\,d\nu\bigg|\\
&\leq C\,\|v\|_{H^3(B)}\|v\|_{H^2(B)}\\
&\leq C\,\|\rho-\bar{\rho}\|_{L^2(\Omega)}^{1/2}\,  \left(\int_\Omega|\nabla^2\rho|^2\,\psi^2\,d\nu+1\right)^{3/4}.
\end{split}
\end{equation}
Combining \eqref{CC412} and \eqref{CC415}, we deduce that
\begin{equation}\label{CC417}
\begin{split}
&\int_\Omega g^{mn}\gamma^{ij}\gamma^{kl}\,v_{s,mnjl}\,v^s_{,ik}\,\psi^2\,d\nu+\int_\Omega\gamma^{ij}\gamma^{kl}g^{mn}\,v_{s,mjl}\,v^s_{,nik}\,\psi^2\,d\nu\\
&\leq C\,\|\rho-\bar{\rho}\|_{L^2(\Omega)}^{1/2}\,  \left(\int_\Omega|\nabla^2\rho|^2\,\psi^2\,d\nu+1\right)^{3/4}.
\end{split}
\end{equation}
The same process of \eqref{CC411}--\eqref{CC417} also leads to
\begin{equation}\label{CC434}
\begin{split}
&\int_\Omega\gamma^{ij}\gamma^{kl}\,v^m_{,msjl}v^s_{,ik}\,\psi^2\,d\nu+\int_\Omega\gamma^{ij}\gamma^{kl}\,v^m_{,mjl}v^s_{,sik}\,\psi^2\,d\nu\\
&\leq C\,\|\rho-\bar{\rho}\|_{L^2(\Omega)}^{1/2}\,  \left(\int_\Omega|\nabla^2\rho|^2\,\psi^2\,d\nu+1\right)^{3/4},
\end{split}
\end{equation}
which together with \eqref{CC416} and \eqref{CC417} gives \eqref{CC410} and finishes the first step.

\textit{Step 2. There is a constant $C$ depending  on $\underline{\rho}$, $M_0$, and $\Omega$, such that}
\begin{equation}\label{CC418}
\begin{split}
&\int_\Omega\gamma^{ij}\gamma^{kl}\,\rho_{,ik}\rho_{,jl}\,\psi^2\,d\nu\\
&\leq C\int_\Omega \gamma^{ij}\gamma^{kl}g^{mn}v_{s,mik}v^s_{,njl}\,\psi^2\,d\nu+C\,\|\rho-\bar{\rho}\|_{L^2(\Omega)}^{1/2}\,  \left(\int_\Omega|\nabla^2\rho|^2\,\psi^2\,d\nu+1\right)^{3/4}.
\end{split}
\end{equation}
Recalling that
 $\gamma=X_t^*\widetilde{\gamma},$ we apply Eulerian coordinates and deduce that \eqref{CC418} is equivalent to
\begin{equation}\label{CC419}
\begin{split}
&\int_\Omega\widetilde{\gamma}^{ij}\widetilde{\gamma}^{kl}\,\partial_{ik}\rho\,\partial_{jl}\rho \,\psi^2\,dx\\
&\leq C\int_\Omega\widetilde{\gamma}^{ij}\widetilde{\gamma}^{kl}\,\partial_{ik}\nabla v\cdot \partial_{jl}\nabla v\,\psi^2\,dx+C\,\|\rho-\bar{\rho}\|_{L^2(\Omega)}^{1/2}\,  \left(\int_\Omega|\nabla^2\rho|^2\,\psi^2\,dx+1\right)^{3/4},
\end{split}
\end{equation}
where $\partial_i$ and $\nabla$ are the usual space derivatives, and $\psi$ denote the truncation function $\psi\big(X_t^{-1}(x)\big)$. According to Lemma \ref{lm29}, \eqref{CC212}, and \eqref{618}, we have $\partial_t\psi+u\cdot\nabla\psi=0,\,|\nabla\psi|\leq C$ for some constant $C$ depending only on $\underline{\rho}$, $M_0$, and $\Omega$.

Note that the support of $\psi$ is strictly containing in $\Omega_7$, while in view of \eqref{57}, \eqref{C527}, and \eqref{CC307}, it also holds that $\mathrm{dist}(F_7,\,\partial\Omega)\leq 8$. Thus according to \eqref{CC221}, we admit $\sum_j(\tilde{N}^j)^2=1$ within the support of $\psi$, which implies that
\begin{equation}\label{CC421}
\begin{split}
\sum_{j}\widetilde{\gamma}^{ij}\widetilde{\gamma}^{jk}\,\psi&=\sum_j\big(\delta^{ij}-\tilde{N}^i\tilde{N}^j\big)\big(\delta^{jk}-\tilde{N}^j\tilde{N}^k\big)\,\psi\\
&=\big(\delta^{ik}-2\tilde{N}^i\tilde{N}^k+\sum_{j}(\tilde{N}^j)^2\tilde{N}^i\tilde{N}^k\big)\psi\\
&=\big(\delta^{ik}-\tilde{N}^i\tilde{N}^k\big)\psi=\widetilde{\gamma}^{ik}\,\psi.
\end{split}
\end{equation}
Then we introduce the notation
$$\mathcal{A}_{jl}\triangleq \mathcal{B}\big(\widetilde{\gamma}^{js}\widetilde{\gamma}^{lr}\partial_{sr}\rho\,\psi-\int_\Omega\widetilde{\gamma}^{js}\widetilde{\gamma}^{lr}\partial_{sr}\rho\,\psi\,dx\big),$$
where $\mathcal{B}(\cdot)$ is provided by Lemma \ref{lemmaz01}. According to \eqref{CC220}, we deduce that
\begin{equation}\label{CC422}
\begin{split}
\mathrm{div}&(\mathcal{A}_{jl})=\widetilde{\gamma}^{js}\widetilde{\gamma}^{lr}\partial_{sr}\rho\,\psi-\int_\Omega\widetilde{\gamma}^{js}\widetilde{\gamma}^{lr}\partial_{sr}\rho\,\psi\,dx,\\
&\|\mathcal{A}_{jl}\|_{H^1(\Omega)}\leq \tilde{C}\int_\Omega\widetilde{\gamma}^{js}\widetilde{\gamma}^{lr}\,\partial_{jl}\rho\,\partial_{sr}\rho\,\psi^2\,dx,
\end{split}
\end{equation}
where in the second line, we have applied \eqref{CC421} to deduce that
\begin{equation*}
\begin{split}
\sum_{j,l}(\widetilde{\gamma}^{js}\widetilde{\gamma}^{lr}\partial_{sr}\rho)^2\,\psi^2=\sum_{j,l}\big(\widetilde{\gamma}^{js}\widetilde{\gamma}^{jm}\big)\,\big(\widetilde{\gamma}^{lr}\widetilde{\gamma}^{ln}\big)\,\partial_{sr}\rho\,\partial_{mn}\rho\,\psi^2=\widetilde{\gamma}^{js}\widetilde{\gamma}^{lr}\,\partial_{jl}\rho\,\partial_{sr}\rho\,\psi^2.
\end{split}
\end{equation*}
In addition, we utilize \eqref{CC222} and \eqref{C401} to check that
\begin{equation}\label{CC423}
\begin{split}
\int_\Omega\widetilde{\gamma}^{js}\widetilde{\gamma}^{lr}\partial_{sr}\rho\,\psi\,dx
&=\int_\Omega\partial_s\big(\widetilde{\gamma}^{js}\widetilde{\gamma}^{lr}\partial_{r}\rho\,\psi)\,dx
-\int_\Omega\partial_s\big(\widetilde{\gamma}^{js}\widetilde{\gamma}^{lr}\,\psi)\,\partial_{r}\rho\,dx\\
&=\int_{\partial\Omega}\big(\tilde{N}_s\cdot\widetilde{\gamma}^{js}\big)\,\widetilde{\gamma}^{lr}\partial_{r}\rho\,dS-\int_\Omega\partial_s\big(\widetilde{\gamma}^{js}\widetilde{\gamma}^{lr}\,\psi)\,\partial_{r}\rho\,dx\\
&=-\int_\Omega\partial_s\big(\widetilde{\gamma}^{js}\widetilde{\gamma}^{lr}\,\psi)\,\partial_{r}\rho\,dx\leq\int_B|\nabla\rho|\,dx.
\end{split}
\end{equation}

Now, let us operate $\psi\cdot\widetilde{\gamma}^{ij}\widetilde{\gamma}^{kl}\,\partial_{ik}$ on \eqref{CC201}$_2$ and arrive at
\begin{equation}\label{CC420}
\widetilde{\gamma}^{ij}\widetilde{\gamma}^{kl}\Delta\partial_{ik}v\,\psi+\widetilde{\gamma}^{ij}\widetilde{\gamma}^{kl}\nabla\mathrm{div}(\partial_{ik}v)\,\psi=\widetilde{\gamma}^{ij}\widetilde{\gamma}^{kl}\,\nabla\partial_{ik}\rho\,\psi.
\end{equation}
Multiplying \eqref{CC420} by $\mathcal{A}_{jl}$ and integrating over $\Omega$, then the right hand side of \eqref{CC420} is given by  
\begin{equation}\label{CC424}
\begin{split}
&\int_\Omega\widetilde{\gamma}^{ij}\widetilde{\gamma}^{kl}\,\nabla\partial_{ik}\rho\cdot\mathcal{A}_{jl}\,\psi\,dx\\
&=-\int_\Omega\widetilde{\gamma}^{ij}\widetilde{\gamma}^{kl}\,\partial_{ik}\rho\cdot\mathrm{div}(\mathcal{A}_{jl})\,\psi\,dx 
-\int_\Omega\nabla\big(\widetilde{\gamma}\cdot\widetilde{\gamma}\cdot\psi\big)\cdot\nabla^2\rho\cdot\mathcal{A}\,dx\\
&=-\int_\Omega\widetilde{\gamma}^{ij}\widetilde{\gamma}^{kl}\,\partial_{ik}\rho\cdot\mathrm{div}(\mathcal{A}_{jl})\,\psi\,dx
+\int_\Omega\nabla\rho\cdot\nabla\left(\mathcal{A}\cdot\nabla\big(\widetilde{\gamma}\cdot\widetilde{\gamma}\cdot\psi\big)\right)\,dx\\
&\leq-\int_\Omega\widetilde{\gamma}^{ij}\widetilde{\gamma}^{kl}\,\partial_{ik}\rho\cdot\mathrm{div}(\mathcal{A}_{jl})\,\psi\,dx+C\int_{B}|\nabla\rho|\cdot\big(|\mathcal{A}|+|\nabla\mathcal{A}|\big)\,dx.
\end{split}
\end{equation}
Moreover, \eqref{CC422} and \eqref{CC423} ensure  that
\begin{equation}\label{CC425}
\begin{split}
&-\int_\Omega\widetilde{\gamma}^{ij}\widetilde{\gamma}^{kl}\,\partial_{ik}\rho\cdot\mathrm{div}(\mathcal{A}_{jl})\,\psi\,dx\\
&=-\int_\Omega\widetilde{\gamma}^{ij}\widetilde{\gamma}^{kl}\,\partial_{ik}\rho\,\partial_{jl}\rho \,\psi^2\,dx+\bigg(\int_\Omega\widetilde{\gamma}^{js}\widetilde{\gamma}^{lr}\partial_{sr}\rho\,\psi\,dx\bigg)^2\\
&\leq-\int_\Omega\widetilde{\gamma}^{ij}\widetilde{\gamma}^{kl}\,\partial_{ik}\rho\,\partial_{jl}\rho \,\psi^2\,dx+\bigg(\int_B|\nabla\rho|\,dx\bigg)^2.
\end{split}
\end{equation}
Therefore, by combining \eqref{CC424} and \eqref{CC425}, we take advantage of \eqref{CC408} and \eqref{CC422} to derive that
\begin{equation}\label{CC426}
\begin{split}
&\int_\Omega\widetilde{\gamma}^{ij}\widetilde{\gamma}^{kl}\,\nabla\partial_{ik}\rho\cdot\mathcal{A}_{jl}\,\psi\,dx\\
&\leq-\int_\Omega\widetilde{\gamma}^{ij}\widetilde{\gamma}^{kl}\,\partial_{ik}\rho\,\partial_{jl}\rho \,\psi^2\,dx+C\,\|\rho-\bar{\rho}\|_{L^2(\Omega)}^{1/2}\,  \left(\int_\Omega|\nabla^2\rho|^2\,\psi^2\,d\nu+1\right)^{3/4}.
\end{split}
\end{equation}

The left hand side of \eqref{CC420} follows in the similar process. We check that
\begin{equation*}
\begin{split}
&\int_\Omega\widetilde{\gamma}^{ij}\widetilde{\gamma}^{kl}\Delta\partial_{ik}v\cdot\mathcal{A}_{jl}\,\psi\,dx\\
&=-\int_\Omega\widetilde{\gamma}^{ij}\widetilde{\gamma}^{kl}\,\nabla\partial_{ik}v\cdot\nabla\mathcal{A}_{jl}\,\psi\,dx-\int_\Omega\nabla\big(\widetilde{\gamma}\cdot\widetilde{\gamma}\cdot\psi\big)\cdot\nabla^3v\cdot\mathcal{A}\,dx\\
&=-\int_\Omega\widetilde{\gamma}^{ij}\widetilde{\gamma}^{kl}\,\nabla\partial_{ik}v\cdot\nabla\mathcal{A}_{jl}\,\psi\,dx
+\int_\Omega\nabla^2v\cdot\nabla\left(\mathcal{A}\cdot\nabla\big(\widetilde{\gamma}\cdot\widetilde{\gamma}\cdot\psi\big)\right)\,dx\\
&\geq-\epsilon^{-1}\int_\Omega\widetilde{\gamma}^{ij}\widetilde{\gamma}^{kl}\,\nabla\partial_{ik}v\cdot\nabla\partial_{jl}v\,\psi^2\,dx-\epsilon\|\mathcal{A}\|_{H^1(\Omega)}^2-C\int_{B}|\nabla^2 v|\cdot\big(|\mathcal{A}|+|\nabla\mathcal{A}|\big)\,dx,
\end{split}
\end{equation*}
Then, for $\tilde{C}$ in \eqref{CC422}, we take $\epsilon=(4\tilde{C})^{-1}$ and apply \eqref{CC408} together with \eqref{CC422} to deduce that
\begin{equation}\label{CC427}
\begin{split}
&\int_\Omega\widetilde{\gamma}^{ij}\widetilde{\gamma}^{kl}\Delta\partial_{ik}v\cdot\mathcal{A}_{jl}\,\psi\,dx\\
&\geq-4\int_\Omega\widetilde{\gamma}^{ij}\widetilde{\gamma}^{kl}\,\nabla\partial_{ik}v\cdot\nabla\partial_{jl}v\,\psi^2\,dx-\frac{1}{4}\int_\Omega\widetilde{\gamma}^{ij}\widetilde{\gamma}^{kl}\,\partial_{ik}\rho\,\partial_{jl}\rho \,\psi^2\,dx\\
&\quad-C\,\|\rho-\bar{\rho}\|_{L^2(\Omega)}^{1/2}\,  \left(\int_\Omega|\nabla^2\rho|^2\,\psi^2\,d\nu+1\right)^{3/4}.
\end{split}
\end{equation}
By the same method, we also have
\begin{equation*} 
\begin{split}
&\int_\Omega\widetilde{\gamma}^{ij}\widetilde{\gamma}^{kl}\nabla\mathrm{div}(\partial_{ik}v)\cdot\mathcal{A}_{jl}\,\psi\,dx\\
&\geq-4\int_\Omega\widetilde{\gamma}^{ij}\widetilde{\gamma}^{kl}\,\partial_{ik}\mathrm{div}v\cdot\partial_{jl}\mathrm{div}v\,\psi^2\,dx-\frac{1}{4}\int_\Omega\widetilde{\gamma}^{ij}\widetilde{\gamma}^{kl}\,\partial_{ik}\rho\,\partial_{jl}\rho \,\psi^2\,dx\\
&\quad-C\,\|\rho-\bar{\rho}\|_{L^2(\Omega)}^{1/2}\,  \left(\int_\Omega|\nabla^2\rho|^2\,\psi^2\,d\nu+1\right)^{3/4}.
\end{split}
\end{equation*}
which along with \eqref{CC426} and \eqref{CC427} gives
\begin{equation*} 
\begin{split}
&\int_\Omega\widetilde{\gamma}^{ij}\widetilde{\gamma}^{kl}\,\partial_{ik}\rho\,\partial_{jl}\rho \,\psi^2\,dx\\
&\leq C\int_\Omega\widetilde{\gamma}^{ij}\widetilde{\gamma}^{kl}\,\partial_{ik}\nabla v\cdot \partial_{jl}\nabla v\,\psi^2\,dx+C\,\|\rho-\bar{\rho}\|_{L^2(\Omega)}^{1/2}\,  \left(\int_\Omega|\nabla^2\rho|^2\,\psi^2\,dx+1\right)^{3/4},
\end{split}
\end{equation*}
Thus \eqref{CC418} and \eqref{CC419} are valid. The second step is therefore finished.

Finally, by collecting \eqref{CC410} and \eqref{CC418}, we arrive at
\begin{equation}\label{CC428}
\begin{split}
-L_2\leq-C\int_\Omega\widetilde{\gamma}^{ij}\widetilde{\gamma}^{kl}\,\partial_{ik}\rho\,\partial_{jl}\rho \,\psi^2\,dx+C\,\|\rho-\bar{\rho}\|_{L^2(\Omega)}^{1/2}\,  \left(\int_\Omega|\nabla^2\rho|^2\,\psi^2\,dx+1\right)^{3/4}.
\end{split}
\end{equation}
Let us set $\mathbf{\Phi}=\mathbf{\Psi}=\gamma$ in \eqref{CC405} and utilize \eqref{C611}, \eqref{C0613}, and \eqref{CC428} to infer that
\begin{equation*} 
\begin{split}
&\frac{d}{dt}\int\gamma^{ij}\gamma^{kl}\rho_{,ik}\rho_{,jl}\psi^2d\nu
+\int_\Omega\gamma^{ij}\gamma^{kl}\,\rho_{,ik}\rho_{,jl}\,\psi^2\,d\nu\\
&\leq C\,\big(\|\rho-\bar{\rho}\|_{L^\infty(B)}
^{1/2}+\|\nabla u\|_{L^\infty(B)}\big)\bigg(\int|\nabla^2\rho|^2\psi^2d\nu+1\bigg)
+\mathcal{R}_3,
\end{split}
\end{equation*}
which gives \eqref{CC409} and the proof is therefore completed.
\end{proof}
 
Next we turn to the estimates on the mixed directions.
\begin{lemma}
Suppose that \eqref{C501}, Lemmas \ref{41}--\ref{L35}, and Lemma \ref{lm53} are valid, then there are constant $C$ and $C_1$ depending only on $\underline{\rho}$, $M_0$, and $\Omega$, such that $\forall t\in[0,T]$, 
\begin{equation}\label{CC429}
\begin{split}
&\frac{d}{dt}\int_\Omega\gamma^{ij}N^kN^l\,\rho_{,ik}\rho_{,jl}\,\psi^2\,d\nu
+\int_\Omega\gamma^{ij}N^kN^l\,\rho_{,ik}\rho_{,jl}\,\psi^2\,d\nu\\
&\leq C\|\nabla u\|_{L^\infty(B)}\int_\Omega|\nabla^2\rho|^2\,\psi^2\,d\nu+\mathcal{R}_3
+C_1\int_\Omega \gamma^{ij}\gamma^{kl}
g^{mn}\,v^s_{,ikm}\,v_{s,jln}\,\psi^2\,d\nu,
\end{split}
\end{equation}
\end{lemma}
\begin{proof}
Taking $\mathbf{\Phi}=\gamma$ and $\mathbf{\Psi}=N\otimes N$ in \eqref{CC405}, then we check that the \textbf{principal term} satisfies
\begin{equation}\label{C634}
\begin{split}
&\int_\Omega\rho\,\gamma^{ij}N^kN^l\, v^s_{,sik}\rho_{,jl}\,\psi^2\,d\nu\\
&=\frac{1}{2}\int_\Omega\rho\,\gamma^{ij}N^kN^l\,F _{v,ik}\,\rho_{,jl}\,\psi^2\,d\nu 
+\frac{1}{2}\int_\Omega\rho\,\gamma^{ij}N^kN^l\, \rho_{,ik}\rho_{,jl}\,\psi^2\,d\nu.\\
\end{split}
\end{equation}
Note that \eqref{63} ensures that
\begin{equation*}
\begin{split}
g^{mn}(v_{k,m}-v_{m,k})_{,ni}&=-F_{v,ik},
\end{split}
\end{equation*}
therefore the first term in \eqref{C634} is handled by
\begin{equation*}
\begin{split}
&\int_\Omega\rho\,\gamma^{ij}N^kN^l\,F_{v,ik}\rho_{,jl}\,\psi^2\,d\nu\\
&=-\int_\Omega\rho\,\gamma^{ij}N^kN^lg^{mn}\,(v_{k,m}-v_{m,k})_{,ni}\rho_{,jl}\,\psi^2\,d\nu\\
&=-\int_\Omega\rho\,\gamma^{ij}N^kN^l(\gamma^{mn}+N^mN^n)(v_{k,m}-v_{m,k})_{,ni}\rho_{,jl}\,\psi^2\,d\nu\\
&=-\int_\Omega\rho\,\gamma^{ij}N^kN^l\gamma^{mn}\,(v_{k,m}-v_{m,k})_{,ni}\rho_{,jl}\,\psi^2\,d\nu,\\
\end{split}
\end{equation*} 
where in the last line, we have used the fact
\begin{equation*}
N^kN^lN^mN^n(v_{k,m}-v_{m,k})_{,ni}=0.
\end{equation*}
Consequently, in view of \eqref{C501}, we deduce that
\begin{equation}\label{C635}
\begin{split}
&\left|\int_\Omega\rho\,\gamma^{ij}N^kN^l\gamma^{mn}\,(v_{k,m}-v_{m,k})_{,ni}\rho_{,jl}\,\psi^2\,d\nu\right|\\
&\leq C\int_\Omega \gamma^{ij}\gamma^{kl}
g^{mn}\,v^s_{,ikm}\,v_{s,jln}\,\psi^2\,d\nu+\frac{1}{4}\int\rho\gamma^{ij}N^kN^l\rho_{,ik}\rho_{,jl}\psi^2d\nu.
\end{split}
\end{equation}

In conclusion, by taking $\mathbf{\Phi}=\gamma$ and $\mathbf{\Psi}=N\otimes N$ in \eqref{CC405} and recalling the fact $\rho\geq\underline{\rho}/4$ in $\Omega_{10}$, we apply \eqref{C634} and \eqref{C635} to deduce that   
\begin{equation*} 
\begin{split}
&\frac{d}{dt}\int_\Omega\gamma^{ij}N^kN^l\,\rho_{,ik}\rho_{,jl}\,\psi^2\,d\nu
+\int_\Omega\gamma^{ij}N^kN^l\,\rho_{,ik}\rho_{,jl}\,\psi^2\,d\nu\\
&\leq C\|\nabla u\|_{L^\infty(B)}\int_\Omega|\nabla^2\rho|^2\,\psi^2\,d\nu+\mathcal{R}_3
+C_1\int_\Omega \gamma^{ij}\gamma^{kl}
g^{mn}\,v^s_{,ikm}\,v_{s,jln}\,\psi^2\,d\nu,
\end{split}
\end{equation*}
which yields \eqref{CC429} and finishes the proof.
\end{proof}

The estimates on normal directions are similar.  
\begin{lemma}
Suppose that \eqref{C501}, Lemmas \ref{41}--\ref{L35}, and Lemma \ref{lm53} are valid, then there are constants $C$ and $C_2$ depending only on $\underline{\rho}$, $M_0$, and $\Omega$, such that $\forall t\in[0,T]$, 
\begin{equation}\label{CC430}
\begin{split}
&\frac{d}{dt}\int_\Omega\big(N^iN^k\rho_{,ik}\big)^2\,\psi^2\,d\nu
+\int_\Omega\big(N^iN^k\rho_{,ik}\big)^2\,\psi^2\,d\nu\\
&\leq C\|\nabla u\|_{L^\infty(B)}\int_\Omega|\nabla^2\rho|^2\,\psi^2\,d\nu+\mathcal{R}_3
+C_2\int_\Omega\gamma^{ij}\gamma^{kl}
(v^s_{,sik}v^m_{,mjl}+\rho_{,ik}\rho_{,jl})\,\psi^2\,d\nu.
\end{split}
\end{equation}
\end{lemma}
\begin{proof}
Taking $\mathbf{\Phi}=\mathbf{\Psi}=N\otimes N$ in \eqref{CC405}, then we check that the \textbf{principal term} satisfies
\begin{equation}\label{C639}
\begin{split}
&\int_\Omega\rho\,N^iN^jN^kN^l\,v^s_{,sik}\,\rho_{,jl}\,\psi^2\,d\nu\\
&=\frac{1}{2}\int_\Omega\rho\,N^iN^jN^kN^l\,F_{v,ik}\,\rho_{,jl}\,\psi^2\,d\nu+\frac{1}{2}\int_\Omega\rho\, (N^iN^k\rho_{,ik})^2\,\psi^2\,d\nu.\\
\end{split}
\end{equation}
Note that by virtue of \eqref{CC223}, $F_v$ is harmonic and satisfies
\begin{equation*}
g^{ik}F_{v,ki}=(\gamma^{ik}+N^iN^k)F_{v,ki}=0~\Rightarrow~N^iN^k\,F_{v,ki}=-\gamma^{ik}\,F_{v,ki},
\end{equation*}
which implies that the first term in \eqref{C639} is treated by
\begin{equation}\label{C640}
\begin{split}
&\int_\Omega\rho\, N^jN^l\big(N^iN^k\,F_{v,ik}\big)\,\rho_{,jl}\,\psi^2\,d\nu\\
&=-\int_\Omega\rho\,N^jN^l\big(\gamma^{ik}\,F_{v,ik}\big)\,\rho_{,jl}\,\psi^2\,d\nu\\
&=-2\int_\Omega\rho\,\gamma^{ik}N^jN^l\,v^s_{,sik}\,\rho_{,jl}\,\psi^2\,d\nu
+\int_\Omega\rho\,\gamma^{ik}N^jN^l\,\rho_{,ik}\,\rho_{,jl}\,\psi^2\,d\nu\\
&\geq-\frac{1}{4}\int_\Omega(N^iN^k\rho_{,ik})^2\,\psi^2\,d\nu-C\int_\Omega\gamma^{ij}\gamma^{kl}\,(v^s_{,sik}v^m_{,mjl}+\rho_{,ik}\rho_{,jl})\,\psi^2\,d\nu.
\end{split}
\end{equation}

In conclusion, by taking $\mathbf{\Phi}=\mathbf{\Psi}=N\otimes N$ in \eqref{CC405} and recalling the fact $\rho\geq\underline{\rho}/4$ in $\Omega_{10}$, we apply \eqref{C639} and \eqref{C640} to deduce that   
\begin{equation*} 
\begin{split}
&\frac{d}{dt}\int_\Omega\big(N^iN^k\rho_{,ik}\big)^2\,\psi^2\,d\nu
+\int_\Omega\big(N^iN^k\rho_{,ik}\big)^2\,\psi^2\,d\nu\\
&\leq C\|\nabla u\|_{L^\infty(B)}\int_\Omega|\nabla^2\rho|^2\,\psi^2\,d\nu+\mathcal{R}_3
+C_2\int_\Omega\gamma^{ij}\gamma^{kl}
(v^s_{,sik}v^m_{,mjl}+\rho_{,ik}\rho_{,jl})\,\psi^2\,d\nu,
\end{split}
\end{equation*}
which gives \eqref{CC430} and finishes the proof.
\end{proof}
 
Finally, we can close the estimates of $\nabla^2\rho$ in the boundary domain $B$.\\\\
\textit{\underline{Proof of Lemma \ref{lm54}.}}
Let us first illustrate that there is a generic constant $C_3$ depending only on $\underline{\rho}$, $M_0$, and $\Omega$, such that
\begin{equation}\label{CC431}
\begin{split}
&\int_\Omega \gamma^{ij}\gamma^{kl}g^{mn}\,v_{s,mik}\,v^s_{,njl}\,\psi^2\,d\nu\\
&\leq C_3\int_\Omega\gamma^{ij}\gamma^{kl}\,\rho_{,ik}\rho_{,jl}\,\psi^2\,d\nu+C\,\|\rho-\bar{\rho}\|_{L^2(\Omega)}^{1/2}\,  \left(\int_\Omega|\nabla^2\rho|^2\,\psi^2\,d\nu+1\right)^{3/4}.
\end{split}
\end{equation}
In fact, according to \eqref{CC433}, \eqref{CC416}, \eqref{CC417}, and \eqref{CC434}, we argue that
\begin{equation*} 
\begin{split}
&\int_\Omega \gamma^{ij}\gamma^{kl}\big(g^{mn}v_{s,mik}\,v^s_{,njl}+v^s_{,sik}\,v^m_{,mjl}\big)\,\psi^2\,d\nu\\
&\leq \int_\Omega\gamma^{ij}\gamma^{kl}\,v^s_{,sik}\,\rho_{,jl}\,\psi^2\,d\nu+C\,\|\rho-\bar{\rho}\|_{L^2(\Omega)}^{1/2}\,  \left(\int_\Omega|\nabla^2\rho|^2\,\psi^2\,d\nu+1\right)^{3/4},
\end{split}
\end{equation*}
which together with H\"{o}lder's inequality gives \eqref{CC431}. 

Now, let us add \eqref{CC409}, \eqref{CC429} multiplied by $(4C_1C_3)^{-1}$, and \eqref{CC430} multiplied by $(4C_2C_3)^{-1}$ together, and utilize \eqref{CC431} to deduce that
\begin{equation*}
\begin{split}
&\Theta'(t)+\Theta(t)\leq C\,\big(\|\rho-\bar{\rho}\|_{L^\infty(B)}
^{1/2}+\|\nabla u\|_{L^\infty(B)}\big)\bigg(\int|\nabla^2\rho|^2\psi^2d\nu+1\bigg)
+\mathcal{R}_3,\\
\end{split}
\end{equation*}
where the leading term $\Theta(t)$ is given by
$$\Theta(t)\triangleq\int_\Omega\left(\gamma^{ij}\gamma^{kl}\rho_{,ik}\rho_{,jl}
+\frac{1}{4C_1C_3}\gamma^{ij}N^kN^l\rho_{,ik}\rho_{,jl}
+\frac{1}{4C_2C_3}\big(N^iN^k\rho_{,ik}\big)^2\right)\psi^2\,d\nu.$$
Note that in view of \eqref{CC435}, we can find a generic constant $C$, such that
\begin{equation}\label{CC436}
C^{-1}\int_\Omega|\nabla^2\rho(\cdot,t)|^2\,\psi^2\,d\nu\leq\Theta(t)\leq C\int_\Omega|\nabla^2\rho(\cdot,t)|^2\,\psi^2\,d\nu,
\end{equation}
which combined with \eqref{CC406} also implies that 
\begin{equation}\label{C671}
\begin{split}
&\Theta(t)'+\Theta(t)\leq C\,\mathcal{R}_4(t)\cdot\big(\Theta(t)+1\big),
\end{split}
\end{equation}
where $\mathcal{R}_4(t)$ is given by
$$\mathcal{R}_4(t)\triangleq\|\rho-\bar{\rho}\|_{L^\infty(B)}
^{1/8}+\|\nabla u\|_{L^\infty(B)}+\|\nabla\dot{u}\|_{L^2(\Omega_{7})}^2+\|\sqrt{\rho}\dot{u}\|_{L^2(\Omega)}^2+\|\nabla u\|_{L^2(\Omega)}^2.$$
In view of \eqref{GN2}, \eqref{C301}, \eqref{CC307}, and \eqref{614}, we deduce that
\begin{equation*}
\begin{split}
\int_0^T\|\rho-\bar{\rho}\|_{L^\infty(B)}
^{1/8}\,dt&\leq C\int_0^T\big(\|\rho-\bar{\rho}\|_{L^2(\Omega_6)}^{1/56}\|\nabla\rho\|_{L^4(\Omega_6)}^{55/56}+\|\rho-\bar{\rho}\|_{L^2(\Omega_6)}^{1/8}\big)\,dt\\
&\leq C\int_0^T\big(\|\rho-\bar{\rho}\|_{L^2(\Omega_6)}^{1/56}+\|\rho-\bar{\rho}\|_{L^2(\Omega_6)}^{1/8}\big)\,dt\leq CC_0^{1/112},
\end{split}
\end{equation*}
which combined with \eqref{C305}, \eqref{618}, and 
\eqref{C107} gives
$$\int_0^T\mathcal{R}_4(t)\,dt\leq C.$$
Consequently, we utilize Gr\"{o}nwall's inequality in \eqref{C671} to deduce that $\forall t\in[0,T]$
\begin{equation*}
\sup_{0\leq\tau\leq t}\Theta(\tau)+\int_0^t\Theta(\tau)\,d\tau\leq C,
\end{equation*}
which along with \eqref{C12243} and \eqref{CC436} yields \eqref{687} and finishes the proof of Lemma \ref{lm54}.
\thatsall

\subsection{The proof of the extension assertion---Lemma \ref{T71}}\label{SS43}
\quad$~$ With all preparations done, the proof of Lemma \ref{T71} is rather direct.\\\\
\textit{\underline{Proof of Lemma \ref{T71}.}}
Suppose that \eqref{C501} is valid, then Lemmas \ref{41}--\ref{L35} and Lemmas \ref{lm51}--\ref{lm54} are valid, provided that the initial energy satisfies $C_0\leq\min\{\varepsilon_1,\varepsilon_2\}$ with $\varepsilon_1$ and $\varepsilon_2$ given by Lemma \ref{21} and \ref{lm51} respectively. Then we apply \eqref{GN2} and \eqref{CC307} to deduce that
\begin{equation*}
\begin{split}
&\|\rho-\bar{\rho}\|_{L^\infty(B)}\leq C\bigg(\|\rho-\bar{\rho}\|_{L^2(\Omega_6)}^{1/4}\|\nabla^2\rho\|_{L^2(\Omega_6)}^{3/4}+\|\rho-\bar{\rho}\|_{L^2(\Omega_6)}\bigg)\leq CC_0^{1/8},\\
&\|\nabla\rho\|_{L^4(B)}\leq C\bigg(\|\rho-\bar{\rho}\|_{L^2(\Omega_6)}^{1/8}\|\nabla^2\rho\|_{L^2(\Omega_6)}^{7/8}+\|\rho-\bar{\rho}\|_{L^2(\Omega_6)}\bigg)\leq CC_0^{1/16}.\\
\end{split}
\end{equation*}
Consequently, if we select $\varepsilon_3=\min\{\big(\bar{\rho}/4C\big)^8,\,(M_0/C)^{16}\}$ and let $C_0\leq\varepsilon_3$, then in view of \eqref{CC103}, it holds that   $\forall t\in[0,T]$,
\begin{equation}\label{CC437}
\|\rho(\cdot,t)\|_{L^\infty(B)}\leq 2M_0,~~\|\nabla\rho(\cdot,t)\|_{L^4(B)}\leq 2M_0,~~\|\rho^{-1}(\cdot,t)\|_{L^\infty(B)}\leq 2\underline{\rho}^{-1}. 
\end{equation}

We still need to investigate the interior domain $\Omega\setminus B$ and deduce that $\forall t\in[0,T]$,
\begin{equation}\label{C70}
\|\rho(\cdot,t)\|_{L^\infty(\Omega\setminus B)}\leq 2M_0,~~\|\rho^{-1}(\cdot,t)\|_{L^\infty(\Omega_{10}\setminus B)}\leq 2\underline{\rho}^{-1}.
\end{equation}
To achieve the goal, let us rewrite \eqref{11}$_1$ as
\begin{equation}\label{C689}
D_t\rho+ \frac{1}{2}\rho\,(\rho-\bar{\rho})=-\rho F,~~F=2\mathrm{div}u-(\rho-\bar{\rho})=2\mathrm{div}w+F_v.
\end{equation}
Moreover, let us set $C_0^{1/14}\leq\varepsilon_4^{1/14}\triangleq\min\{1/(100C) ,\,\underline{\rho}/(100CM_0)\}$, then in view of \eqref{CC215}, \eqref{C301}, and \eqref{C527}, it holds that
\begin{equation}\label{627}
\begin{split}
&\int_0^T\|F\|_{L^\infty(\Omega\setminus\Omega_1)}\,dt\\
&\leq\int_0^T\bigg(\|\nabla w\|_{L^\infty(\Omega)}+\|\rho-\bar{\rho}\|_{L^2(\Omega)}\bigg)\,dt \leq CC_0^{1/14}\leq \min\{1/100,\,\underline{\rho}/(100M_0)\}.
\end{split}
\end{equation}
Now, we prove \eqref{C70} by contradiction. Suppose that the first part of \eqref{C70} fails, then there must be a flow line $X_t(x_1)$, such that for some $t_1\in[0,T]$,
\begin{equation}\label{C80}
\rho\big(X_{t_1}(x_1)\big)>2M_0~~\mathrm{and}~~X_t(x_1)\subset\Omega\setminus B~~\forall t\in[0,T].
\end{equation}
According to \eqref{CC103}, it is valid that
\begin{equation*}
\rho\big(X_0(x_1)\big)=\rho_0(x_1)\leq M_0<3M_0/2, 
\end{equation*}
Therefore, by the continuity of $\rho$ and the flow line, there must be some $t_0\in(0,t_1)$ such that
\begin{equation}\label{C802}
\rho\big(X_{t_0}(x_1)\big)=3M_0/2~~\mathrm{and}~~\rho\big(X_t(x_0)\big)\geq 3M_0/2~~ \forall t\in[t_0,t_1].
\end{equation}
Observe that $\bar\rho\leq 3M_0/2$ due to \eqref{CC103}, therefore let us integrate \eqref{C689} on $[t_0,t_1]$ and apply \eqref{C501} together with \eqref{627} to deduce that, 
\begin{equation*} 
\begin{split}
&\rho\big(X_{t_1}(x_1)\big)-\rho\big(X_{t_0}(x_1)\big)\\
&\leq -\int_{t_0}^{t_1}\rho F(X_t(x_1),t)\,dt\leq 4M_0\int_{t_0}^{t_1}\|F\|_{L^\infty(\Omega\setminus\Omega_1)}\,dt
\leq M_0/4,
\end{split}
\end{equation*}
which contradicts the fact
\begin{equation*}
\rho\big(X_{t_1}(x_1)\big)-\rho\big(X_{t_0}(x_1)\big)\geq 2M_0-3M_0/2=M_0/2,
\end{equation*}
due to \eqref{C80} and \eqref{C802}. Therefore we must have
$\|\rho(\cdot,t)\|_{L^\infty(\Omega\setminus B)}\leq 2M_0$ for all $t\in[0,T]$, which provides the first part of \eqref{C70}.

The proof of the second part of \eqref{C70} is quite similar. Suppose it is not the case, then there must be a flow line $X_t(x_0')$, such that for some $t_1'\in[0,T]$
\begin{equation}\label{C82}
\rho\big(X_{t_1'}(x_0')\big)<\underline{\rho}/2~~\mathrm{and}~~X_t(x_0')\subset\Omega_{10}\setminus B~~\forall t\in[0,T].
\end{equation}
According to \eqref{CC103}, it is true that 
\begin{equation*}
\rho\big(X_0(x_0')\big)=\rho_0(x_0')\geq\underline{\rho}>2\underline{\rho}/3.
\end{equation*}
Hence, by the continuity of $\rho$ and the flow line, there must be some $t'_0\in(0,t_1')$ such that
\begin{equation}\label{C85}
\rho\big(X_{t'_0}(x_0'))= 2\underline{\rho}/3~~\mathrm{and}~~\rho(X_{t}(x_0'))\leq 2\underline{\rho}/3~~\forall t\in [t'_0,t_1'].
\end{equation}
Observe that $\bar\rho\geq 2\underline{\rho}/3$ due to \eqref{CC103}, therefore let us integrate \eqref{C689} on $[t_0',t_1']$ and apply \eqref{C501} together with \eqref{627} to deduce that 
\begin{equation*} 
\begin{split}
&\rho\big(X_{t'_1}(x_0')\big)-\rho\big(X_{t'_0}(x_0')\big)\\
&\geq 
-\int_{t'_0}^{t_1'}\rho F\big(X_t(x_0'),t\big)\,dt\geq -4M_0\int_{t'_0}^{t_1'}\|F\|_{L^\infty(\Omega_{10}\setminus B)}\,dt\geq-\underline{\rho}/10,
\end{split}
\end{equation*}
which contradicts the fact
\begin{equation*}
\rho\big(X_{t'_1}(x_0')\big)-\rho\big(X_{t'_0}(x_0')\big)\leq\underline{\rho}/2-2\underline{\rho}/3=-\underline{\rho}/6,
\end{equation*}
due to \eqref{C82} and \eqref{C85}. Consequently we argue that $\|\rho^{-1}(\cdot,t)\|_{L^\infty(\Omega_{10}\setminus B)}\leq 2\underline{\rho}^{-1}$ for all $t\in[0,T]$
which finishes the proof of \eqref{C70}.
 
In conclusion, suppose that \eqref{C501} holds and the initial energy satisfies $C_0\leq\varepsilon$ with $\varepsilon\triangleq\min\{\varepsilon_1,\varepsilon_2,\varepsilon_3,\varepsilon_4\}$, then \eqref{CC437} together with \eqref{C70} implies  \eqref{C502}. Therefore the proof of Lemma \ref{T71} is finally completed. \thatsall
 
\section{A priori estimates III: Higher order estimates}\label{S5}
\quad Once we obtain the uniform upper bound of $\rho$, the rest of paper follows in the standard way as \cite{hlx, caili01}, for the sake of completeness, we sketch the proof here. The constants $C$ below vary from line to line, but they depend only on the initial value and $T$. Moreover all assertions are established in the Eulerian coordinates, and we denote $\partial_i=\partial/\partial x^i$.
 
\begin{lemma}\label{T81}
Under the conditions of Theorem \ref{T1} and \eqref{C501}, there is a constant $C$ depending on $(\rho_0,u_0)$, $\underline{\rho}$, $M_0$, $\Omega$, and $T$, such that
\begin{equation}\label{CC504}
\sup_{0 \leq t \leq T}\big(\|\rho(\cdot,t)\|_{W^{1,6}(\Omega)}+\|u(\cdot,t)\|_{H^2(\Omega)}\big)
+\int_{0}^{T}\|\nabla\dot{u}(\cdot,t)\|_{L^2(\Omega)}^2\,dt\leq C.
\end{equation}
\end{lemma}
\begin{proof}
Recalling that \eqref{CC318} yields 
\begin{equation*} 
\begin{split}&
\frac{d}{dt}\bigg(\int_\Omega\rho|\dot{u}|^2\,dx\bigg)+\int_\Omega|\nabla\dot{u}|^{2}\,dx\\
&\leq C\big(\|\nabla u\|_{L^2(\Omega)}+\|\rho-\bar{\rho}\|_{L^2(\Omega)}\big)\,\|\rho\dot{u}\|_{L^2(\Omega)}\left(\int_\Omega\rho|\dot{u}|^2\,dx\right)\\
&\quad+C\bigg(\|\nabla u\|_{L^2(\Omega)}^4+\|\nabla u\|^2_{L^2(\Omega)}+\|\rho-\bar{\rho}\|^4_{L^4(\Omega)}\bigg),
\end{split}
\end{equation*}
which along with the compatible condition \eqref{1116} and \eqref{C301} gives that
\begin{equation}\label{CC503}
\sup_{0\leq t\leq T}\|\sqrt{\rho}\dot{u}(\cdot,t)\|_{L^2(\Omega)}^2+\int_0^T\|\nabla\dot{u}(\cdot,t)\|_{L^2(\Omega)}^2\,dt\leq C.
\end{equation}

Now, taking $\partial_i$ on \eqref{11}$_1$ yields that
\begin{equation}\label{CC502}
\partial_t(\partial_i\rho)+\mathrm{div}(u\,\partial_i\rho)=-\mathrm{div}(\rho\,\partial_iu).
\end{equation}
We multiplying the above equation by $|\nabla\rho|^4\,\partial_i\rho$ and integrate over $\Omega$ to declare that
\begin{equation}\label{C71}
\begin{aligned}
\frac{d}{dt}\|\nabla\rho\|_{L^6(\Omega)}& \leq C\|\nabla u\|_{L^{\infty}(\Omega)}\,\|\nabla \rho\|_{L^{6}(\Omega)}+C\big(\|\nabla^2w\|_{L^{6}(\Omega)}+\|\nabla^2v\|_{L^{6}(\Omega)}\big)\\
&\leq C\big(1+\|\nabla u\|_{L^{\infty}(\Omega)}\big)\|\nabla \rho\|_{L^{6}(\Omega)}+C\|\rho \dot{u}\|_{L^{6}(\Omega)}\\
& \leq C\big(1+\|\nabla u\|_{L^{\infty}(\Omega)}\big)\|\nabla \rho\|_{L^{6}(\Omega)}+C\|\nabla \dot{u}\|_{L^{2}(\Omega)},
\end{aligned}
\end{equation}
where in the  second inequality we have used \eqref{355}. Recalling that \eqref{CC215} ensures that $F_v$ and $\mathrm{rot}v$ are harmonic in $\Omega$, which together with \eqref{C406} gives that
\begin{equation*}
\begin{aligned}
\|\mathrm{div} u\|_{L^{\infty}(\Omega)}
&\leq C\big(\|\nabla w\|_{L^{\infty}(\Omega)}+\|F_v\|_{L^{\infty}(\Omega\setminus\Omega_1)}+\|\nabla v\|_{L^{\infty}(\Omega_1)}+\|\rho-\bar{\rho}\|_{L^{\infty}(\Omega)}\big)\\
&\leq C\big(\|\nabla w\|_{L^{\infty}(\Omega)}+\|\rho\|_{C^\alpha(B)}+\|\rho-\bar{\rho}\|_{L^{\infty}(\Omega)}\big)\\
&\leq C\big(1+\|\nabla w\|_{L^{\infty}(\Omega)}\big),
\end{aligned}
\end{equation*}
where the last line is due to \eqref{C12241}, while we have applied \eqref{C406} and \eqref{CC307} in the domain $\Omega_1\subset B$ to obtain the third line. The same conclusion holds if we replace $\mathrm{div}u$ by $\mathrm{rot}u$, thus we utilize \eqref{24}$_1$ to deduce that
\begin{equation}\label{C77}
\begin{aligned}
\|\nabla u\|_{L^{\infty}(\Omega)}&\leq C\big(\|\mathrm{div} u\|_{L^{\infty}}
+\|\mathrm{rot} u\|_{L^{\infty}}\big)\log\big(e+\|\nabla^{2} u\|_{L^{6}(\Omega)}\big)+C\|\nabla u\|_{L^{2}(\Omega)}+C \\
& \leq C\big(1+\|\nabla w\|_{L^{\infty}(\Omega)}\big)\log\big(e+\|\nabla \dot{u}\|_{L^{2}(\Omega)}+\|\nabla \rho\|_{L^{6}(\Omega)}\big) \\
& \leq C\big(1+\|\nabla w\|_{L^{\infty}(\Omega)}\big)\,\big(1+\|\nabla \dot{u}\|_{L^{2}(\Omega)}\big)\,\log \big(e+\|\nabla \rho\|_{L^{6}(\Omega)}\big),
\end{aligned}
\end{equation}
where in the second line, we have made use of the fact that for $q\in[2,6]$,
\begin{equation}\label{CC501}
\begin{split}
\|\nabla^{2} u\|_{L^{q}(\Omega)}
&\leq C\big(\|\nabla^{2} w\|_{L^{q}(\Omega)}+\|\nabla^{2} v\|_{L^{q}(\Omega)}\big)\\
&\leq C\big(\|\rho\dot{u}\|_{L^{q}(\Omega)}+\|\nabla\rho\|_{L^{q}(\Omega)}\big)\\
&\leq C\big(\|\nabla\dot{u}\|_{L^2(\Omega)}+\|\nabla\rho\|_{L^q(\Omega)}\big),
\end{split}
\end{equation}
which is due to \eqref{355}. Thus, 
in view of Lemmas \ref{21}--\ref{31}, \eqref{CC503}, and \eqref{C77}, we take advantage of Gr\"{o}nwall's inequality in \eqref{C71} to deduce that
\begin{equation*} 
\sup_{0\leq t\leq T}\|\nabla \rho(\cdot,t)\|_{L^{6}(\Omega)}\leq C.
\end{equation*}
which along with \eqref{CC503} and \eqref{CC501}  leads to
\begin{equation*}
\begin{split}
\sup_{0\leq t\leq T}&\|u(\cdot,t)\|_{H^2(\Omega)}\leq \sup_{0\leq t\leq T}C\,\big(\|\sqrt{\rho}\dot{u}(\cdot,t)\|_{L^2(\Omega)}+\|\nabla\rho(\cdot,t)\|_{L^2(\Omega)}\big)\leq C,\\
\end{split}
\end{equation*}
and finishes the proof of \eqref{CC504}.
\end{proof}

Based on Lemma \ref{T81}, we can establish some further regularity assertions.

\begin{lemma}\label{T83}
Under the conditions of Theorem \ref{T1} and \eqref{C501}, there is a constant $C$ depending on $(\rho_0,u_0)$, $\underline{\rho}$, $M_0$, $\Omega$, and $T$, such that for $\sigma(t)=\min\{1,t\}$,
\begin{equation}\label{CC505}
\begin{split}
\sup_{0\leq t \leq T}&\big(\|\rho_{t}(\cdot,t)\|_{H^{1}(\Omega)}+\|\rho(\cdot,t)\|_{H^2(\Omega)}+\sigma\|u_{t}(\cdot,t)\|_{H^{1}(\Omega)}^2\big)\leq C,\\
&\int_{0}^{T}\big(\|\rho_{t t}(\cdot,t)\|_{L^{2}(\Omega)}^{2}+\sigma\|\sqrt{\rho} u_{t t}(\cdot,t)\|_{L^{2}(\Omega)}^{2}\big)\,dt\leq C.
\end{split}
\end{equation}
\end{lemma}
\begin{proof}
Taking $\partial_{ij}$ on \eqref{11}$_1$ yields that
\begin{equation}\label{CC514}
\begin{split}
\partial_t(\partial_{ij}\rho)+\mathrm{div}(u\,\partial_{ij}\rho)=
-\mathrm{div}(\partial_ju\,\partial_i\rho)-\mathrm{div}(\partial_j\rho\,\partial_iu)-\mathrm{div}(\rho\,\partial_{ij}u).
\end{split}
\end{equation}
Multiplying the above equation by $\partial_{ij}\rho$ and utilizing \eqref{355} together with \eqref{CC504} gives that
\begin{equation}\label{C74}
\begin{aligned}
&\frac{d}{dt}\|\nabla^{2}\rho\|_{L^{2}(\Omega)}^{2}\leq
C\big(1+\|\nabla u\|_{L^{\infty}(\Omega)}\big)\|\nabla^{2}\rho\|_{L^{2}(\Omega)}^{2}+C\,\big(\|\nabla \dot{u}\|_{L^{2}(\Omega)}^{2}+1\big).
\end{aligned}
\end{equation}
Thus, in view of Lemmas \ref{21}--\ref{31}, \eqref{CC504}, and \eqref{C77}, we apply Gr\"{o}nwall's inequality in \eqref{C74} to declare that
\begin{equation}\label{CC506}
\sup_{0 \leq t \leq T}\|\rho(\cdot,t)\|_{H^{2}(\Omega)}\leq C,
\end{equation}
which along with \eqref{11}$_1$, \eqref{CC504}, and \eqref{CC502} ensures that
\begin{equation}\label{CC507}
\sup_{0 \leq t \leq T}\|\rho_{t}(\cdot,t)\|_{H^{1}(\Omega)}
\leq \sup_{0 \leq t \leq T}C\,\big(\|u(\cdot,t)\|_{H^{2}(\Omega)}\|\rho(\cdot,t)\|_{H^{2}(\Omega)}\big)\leq C.
\end{equation}
Then, let us take derivatives with respect to $t$ twice on \eqref{11}$_1$ and arrive at
\begin{equation*} 
\rho_{t t}=-\rho_{t}\,\mathrm{div} u-\rho\,\mathrm{div} u_{t}-u_t\cdot\nabla\rho-u\cdot\nabla\rho_t.
\end{equation*}
Thus, we make use of \eqref{CC504}, \eqref{CC506}, and \eqref{CC507} to argue that
\begin{equation}\label{CC510}
\begin{split}
\int_0^T\|\rho_{tt}\|_{L^2(\Omega)}^2\,dt&\leq\int_0^T\big(\|\rho_{t}\|_{H^{1}(\Omega)}\|u\|_{H^{2}(\Omega)}+\|\nabla u_{t}\|_{H^{1}(\Omega)}\|\rho\|_{H^2(\Omega)}\big)^2\,dt\\
&\leq\int_0^TC\big(1+\|\nabla\dot{u}\|_{L^2(\Omega)}^2\big)\,dt\leq C.
\end{split}
\end{equation}

By taking derivative with respect to $t$ on \eqref{11}$_2$ and multiplying by $u_{tt}$, we integrate the resulting equation over $\Omega$ and obtain that
\begin{equation}\label{CC509}
\begin{aligned}
\frac{d}{dt}\int_\Omega\big(|\nabla u_t|^2+(\mathrm{div}u_t)^2\big)\,dx+2 \int_\Omega \rho|u_{tt}|^{2}\,dx=R_1(t)'+R_2.
\end{aligned}
\end{equation}
where the remaining terms are given by
\begin{equation*}
\begin{split}
R_1(t)&=\int_\Omega\bigg(-\rho_{t}\left|u_{t}\right|^{2}-2\rho_{t}\,u\cdot\nabla u\cdot u_{t} 
+2\rho_{t}\,\mathrm{div}u_{t}\bigg)\,dx,\\
R_2(t)&=\int_\Omega\bigg(\rho_{tt}\left|u_{t}\right|^{2}+2  \big(\rho_{t}\,u\cdot\nabla u\big)_t\,u_{t}-2\rho (u \cdot\nabla u)_t\,u_{tt}-2\rho_{tt}\,\mathrm{div}u_{t}\bigg)\,dx.
\end{split}
\end{equation*}
The standard calculations in \cite[Lemma 3.9]{hlx1} together with Lemmas \ref{41}--\ref{21} and \eqref{CC505}$_1$ provide  
\begin{equation}\label{CC508} 
\begin{aligned}
R_1(t)&\leq\frac{1}{2}\int_\Omega\big(|\nabla u_t|^2+(\mathrm{div}u_t)^2\big)\,dx+C,\\
R_2(t)&\leq \|\sqrt{\rho}u_{tt}\|_{L^2(\Omega)}^{2}+C\bigg(\|\nabla u_{t}\|_{L^2(\Omega)}^{4}+
\|\rho_{tt}\|_{L^2(\Omega)}^{2}+1\bigg).
\end{aligned}
\end{equation}
Substituting \eqref{CC508} into \eqref{CC509} and multiplying $\sigma$ yield that
\begin{equation*}
\begin{split}
\frac{d}{dt}\big(\sigma\mathbf{G}(t)&\big)+2\sigma\int_\Omega \rho|u_{tt}|^{2}\,dx \leq C\bigg(\|\nabla u_{t}\|_{L^2(\Omega)}^{2}+\|\rho_{tt}\|_{L^2(\Omega)}^{2}+1\bigg)+\sigma\|\nabla u_t\|_{L^2(\Omega)}^4,\\
&\mathbf{G}(t)\triangleq  \int_\Omega\big(|\nabla u_t|^2+(\mathrm{div}u_t)^2\big)\,dx-R_1(t)+C\geq\frac{1}{2}\int_\Omega|\nabla u_t|^2\,dx.
\end{split}
\end{equation*}
Consequently, we utilize Gr\"{o}nwall's inequality, \eqref{CC504}, and \eqref{CC510} to argue that
\begin{equation*}
\sup_{0\leq t\leq T}\sigma\|u_{t}(\cdot,t)\|_{H^{1}(\Omega)}^{2}+\int_0^T \sigma\|\sqrt{\rho}u_{tt}(\cdot,t)\|_{L^{2}(\Omega)}^{2}\,dt\leq C,
\end{equation*}
which combined with \eqref{CC506}--\eqref{CC510} provides \eqref{CC505} and finishes the proof.
\end{proof}

In view of Lemma \ref{T83}, the next corollary is direct.

\begin{lemma} 
Under the conditions of Theorem \ref{T1} and \eqref{C501}, there is a constant $C$ depending on $(\rho_0,u_0)$, $\underline{\rho}$, $M_0$, $\Omega$, and $T$, such that for $q\in(3,6)$ and $r\in(1,7/6)$,  
\begin{equation}\label{CC511}
\begin{split}
&\sup_{t\in[0,T]}\big(\|\rho(\cdot,t)\|_{W^{2,q}(\Omega)}+\sigma\|u(\cdot,t)\|_{H^3(\Omega)}^2\big)\leq C,\\
\int_0^T\big(\|u(&\cdot,t)\|_{H^3(\Omega)}^2+\|
u(\cdot,t)\|^{r}_{W^{3,q}(\Omega)}+\sigma\| u_t(\cdot,t)\|_{H^2(\Omega)}^2\big)\,dt\leq C.\\
\end{split}
\end{equation}
\end{lemma}
\begin{proof}
First, we utilize \eqref{355} and \eqref{CC505} to declare that
\begin{equation}\label{CC513}
\begin{split}
&\|u\|_{W^{3,q}(\Omega)}\leq C\big(\|\rho \dot{u}\|_{W^{1,q}(\Omega)}+\|\rho\|_{W^{2,q}(\Omega)}\big).
\end{split}
\end{equation}
In particular, by taking $q=2$ in \eqref{CC513} and utilizing \eqref{CC504} together with \eqref{CC505}, we infer that 
\begin{equation}\label{CC512}
\sup_{0\leq t\leq T}\sigma\|u(\cdot,t)\|_{H^{3}(\Omega)}^2+\int_0^T\|u(\cdot,t)\|_{H^3(\Omega)}^2\,dt\leq C.
\end{equation}

Next, taking derivatives with respect to $t$ on \eqref{11}$_2$ gives the elliptic system
\begin{equation*} 
\begin{cases}
\Delta u_t+\nabla\div u_t=(\rho\dot{u})_t+\nabla\rho_t~~&\mathrm{in}\ \Omega,\\
u_t=0~~&\mathrm{on}\ \partial\Omega.
\end{cases}
\end{equation*}
Thus \eqref{355} together with \eqref{CC504} and \eqref{CC505} ensures that
\begin{equation}\label{CC517}
\begin{split}
&\|\nabla^2u_t\|_{L^2(\Omega)}\\
&\leq C\big(\|(\rho\dot{u})_t\|_{L^2(\Omega)}+\|\nabla\rho_t\|_{L^2(\Omega)}\big)\\
&\leq C \|\sqrt{\rho}u_{tt}\|_{L^2(\Omega)}+C\big(\|\rho_t\|_{H^1(\Omega)}+1\big)\,\big(\|u_t\|_{H^1(\Omega)}+1\big)\\
&\leq C\big(\|\sqrt{\rho}u_{tt}\|_{L^2(\Omega)}+\|\nabla u_t\|_{L^2(\Omega)}+1\big),
\end{split}
\end{equation}
which along with \eqref{GN}, \eqref{CC504}, and \eqref{CC505} also yields that for $q\in (3,6)$ and $\theta=(6-q)/2q$,
\begin{equation}\label{816}
\begin{split}
&\|\nabla(\rho\dot u)\|_{L^q(\Omega)}\\
&\leq C\big(\|\nabla u_t\|_{L^q(\Omega)}+\|\nabla(u\cdot\nabla u)\|_{L^q(\Omega)}\big)\\
&\leq C\big(\|\nabla u_t\|_{L^2(\Omega)}^{\theta}\|\nabla^2 u_t\|_{L^2(\Omega)}^{1-\theta}+\|\nabla u_t\|_{L^2(\Omega)}+\|u\|_{H^3(\Omega)}\big)\\
&\leq C\big(1+\sigma^{-\frac{1}{2}}+\sigma^{-\frac{1}{2}}(\sigma\|\sqrt{\rho}u_{tt}\|_{L^2(\Omega)}^2)^{(1-\theta)/2}+\|u\|_{H^3(\Omega)}\big).
\end{split}
\end{equation}
Therefore, collecting \eqref{CC505} and \eqref{CC513}--\eqref{816} gives that
\begin{equation}\label{809}
\int_0^T\big(\|
u(\cdot,t)\|^{r}_{W^{3,q}(\Omega)}+\sigma\| u_t(\cdot,t)\|_{H^2(\Omega)}^2\big)\,dt\leq C,
\end{equation}
 
At last, by multiplying \eqref{CC514} by $|\nabla^2\rho|^{q-2}\,\partial_{ij}\rho$ and applying \eqref{CC505}, we carry out similar calculations as \eqref{C74} and infer that
\begin{equation}\label{C73}
\begin{split}
\frac{d}{dt} \|\nabla^2 \rho\|_{L^q(\Omega)}\leq& C\big(\|\nabla u\|_{L^\infty(\Omega)}\|\nabla^2 \rho\|_{L^q(\Omega)}
+\|u\|_{W^{3,q}(\Omega)}\big),\\
\end{split}
\end{equation}
Thus, in view of Lemmas \ref{21}--\ref{31}, \eqref{CC504}, \eqref{C77}, and \eqref{809}, we apply Gr\"{o}nwall's inequality in \eqref{C73} to argue that
\begin{equation*}
\begin{split}
\sup_{0\leq t\leq T}\|\rho(\cdot,t)\|_{W^{2,q}(\Omega)}\leq C,
\end{split}
\end{equation*}
which combined with \eqref{CC512} and \eqref{809} gives \eqref{CC511} and finishes the proof.
\end{proof}

The last regularity assertion of this section guarantees that the solution obtained in Theorem \ref{T1} is a classical one.
\begin{lemma}\label{T85}
Under the conditions of Theorem \ref{T1} and \eqref{C501}, there is a constant $C$ depending on $(\rho_0,u_0)$, $\underline{\rho}$, $M_0$, $\Omega$, and $T$, such that for $q\in(3,6),$ 
\begin{equation}\label{CC515}
\sup_{0\leq t\leq T}\sigma\big(\|u(\cdot,t)\|_{W^{3,q}(\Omega)}+\|u_t(\cdot,t)\|_{H^2(\Omega)}\big)
+\int_{0}^T\sigma^2\|u_{tt}(\cdot,t)\|_{H^1(\Omega)}^2\,dt\leq C.
\end{equation}
\end{lemma}
\begin{proof}
Let us take derivatives with respect to $t$ twice on \eqref{11}$_2$, then we multiply the resulting equation by $2u_{tt}$ and integrate over $\Omega$ to deduce that 
\begin{equation}\label{813}
\begin{split}
&\frac{d}{dt}\int_\Omega\rho|u_{tt}|^2\,dx+\int_{\Omega}\big(|\nabla u_{tt}|^2+(\mathrm{div}u_{tt})^2\big)\,dx\\
&=-8\int_{\Omega}\rho u_{tt}\cdot(u\cdot\nabla u_{tt})\,dx-2\int_{\Omega}(\rho u)_t\cdot \big(\nabla (u_t\cdot u_{tt})+2\nabla
u_t\cdot u_{tt}\big)\,dx\\
&\quad -2\int_\Omega\big(\rho_{tt}u+2\rho_t u_t\big)\cdot\nabla u\cdot u_{tt}\,dx-2\int_\Omega\rho_{tt}\cdot\nabla u\cdot u_{tt}\,dx\\
&\quad+2\int_{\Omega}\rho_{tt}\,\mathrm{div}u_{tt}\,dx.
\end{split}
\end{equation}
The standard calculations in \cite[Lemma 3.11]{hlx1} together with \eqref{CC505} make sure that
\begin{equation*}
\mathrm{The\ right\ hand\ side\ of\ \eqref{813}}\leq\frac{1}{2}\,\|\nabla u_{tt}\|_{L^2(\Omega)}^2+C\big(\sigma^{-3/2}+\|\sqrt{\rho}u_{tt}\|^2_{L^2(\Omega)}+\|\rho_{tt}\|^2_{L^2(\Omega)}\big),
\end{equation*}
which together with \eqref{813} leads to
\begin{equation}\label{CC516}
\begin{split}
&\frac{d}{dt}\int_\Omega\rho|u_{tt}|^2\,dx+\frac{1}{2}\int_{\Omega}|\nabla u_{tt}|^2\,dx\leq C\big(\sigma^{-3/2}+\|\sqrt{\rho}u_{tt}\|^2_{L^2(\Omega)}+\|\rho_{tt}\|^2_{L^2(\Omega)}\big).
\end{split}
\end{equation}
Therefore, we apply \eqref{CC505}, \eqref{CC511}, and Gr\"{o}nwall's inequality in \eqref{CC516} to infer that 
\begin{equation}\label{814}
\begin{split}
\sup_{0\leq t\leq T}\sigma^2\|\sqrt{\rho}u_{tt}(\cdot,t)\|_{L^2(\Omega)}^2+\int_{0}^T\sigma^2\|\nabla u_{tt}(\cdot,t)\|_{L^2(\Omega)}^2\,dt\leq C.
\end{split}
\end{equation}
In view of \eqref{CC511}, \eqref{CC513}, \eqref{CC517}, and \eqref{814}, we also declare that
\begin{equation*} 
\begin{split}
\sup_{0\leq t\leq T}\sigma\big(\|u(\cdot,t)\|_{W^{3,q}(\Omega)}+\|u_t(\cdot,t)\|_{H^2(\Omega)}\big)\leq C,\\
\end{split}
\end{equation*}
which together with \eqref{814} yields \eqref{CC515} and finishes the proof.
\end{proof}

\section{Proofs of Theorems \ref{T1} and \ref{T2}}\label{S6}
\quad$~$ With the help of a priori estimates given by previous sections, we can prove Theorems \ref{T1} and \ref{T2} by the standard method.\\

{\textit{Proof of Theorem \ref{T1}.}} Suppose that $\varepsilon$ is given by Theorem \ref{T71} and the initial energy satisfies $C_0<\varepsilon$. Then according to Lemma \ref{u21}, there exists 
$T_*>0$, such that \eqref{11}--\eqref{CC102} admits a unique classical solution $(\rho,u)$ on $\Omega\times(0,T_*]$. Let us define
\begin{equation*}
\begin{split}
\mathbf{H}_1(t)&\triangleq\|\rho(\cdot,t)\|_{L^{\infty}(\Omega)},~~\mathbf{H}_2(t)\triangleq\|\nabla\rho(\cdot,t)\|_{L^{4}(B)},~~
\mathbf{H}_3(t)\triangleq\|\rho^{-1}(\cdot,t)\|_{L^\infty(\Omega_{10})}.
\end{split}
\end{equation*}
Under the conditions of Theorem \ref{T1}, we mention that \eqref{CC103} guarantees  
\begin{equation*}
\mathbf{H}_1(0)<2M_0,~~\mathbf{H}_2(0)<2M_0,~~ \mathbf{H}_3(0)<2\underline{\rho}^{-1}.
\end{equation*}
Thus it is well defined to introduce
\begin{equation*}
T^*\triangleq\sup\big\{t>0\,\big|\,\forall\tau\in[0,t],\, \mathbf{H}_1(\tau)\leq 4M_0,\,\mathbf{H}_2(\tau)\leq 4M_0,\,\mathbf{H}_3(\tau)\leq 4\underline{\rho}^{-1}\big\}.
\end{equation*}
Moreover, in view of Lemma \ref{T81}-\ref{T85}, for any $0<\tau<T\leq T^*$ with $T$ finite, it holds that
\begin{equation}\label{CC801}
\begin{cases}
\rho\in C\big([0,T];W^{2,q}(\Omega)\big),\\
u\in C\big([0,T];H^2(\Omega)\big)\cap L^\infty\big(\tau,T; W^{3,q}(\Omega)\big),\\
u_t\in L^{\infty}\big(\tau,T; H^2(\Omega)\big)\cap H^1\big(\tau,T; H^1(\Omega)\big),\\
\sqrt{\n}u_t\in L^\infty\big(0,\infty;L^2(\Omega)\big).
\end{cases}
\end{equation}
Consequently, by virtue of the standard embedding
\begin{equation*}
L^\infty\big(\tau ,T;H^1(\Omega)\big)\cap H^1\big(\tau ,T;H^{-1}(\Omega)\big)\hookrightarrow
C\big([\tau ,T];L^q(\Omega)\big)~~\forall q\in[2,6), 
\end{equation*}
we also deduce that
\begin{equation}\label{92}
\sqrt{\rho}u_t,\,\sqrt{\rho}\dot{u}\in C\big([\tau,T];L^2(\Omega)\big).
\end{equation}

Now, we claim that $T^*=\infty$. If it is not true and $T^*<\infty$, by virtue of Theorem \ref{T71}, we deduce that $\forall\tau\in[0,T^*]$,
\begin{equation}\label{CC802}
\mathbf{H}_1(\tau)\leq 2M_0,~~
\mathbf{H}_2(\tau)\leq 2M_0,~~ 
\mathbf{H}_3(\tau)\leq 2\underline{\rho}^{-1}.
\end{equation}
In addition, according to \eqref{CC801} and \eqref{92}, $\big(\rho(x,T^*),u(x,T^*)\big)$ satisfies the initial conditions \eqref{1114} and \eqref{1116} with $g(x)\triangleq\sqrt{\rho}\dot{u}(x, T^*)$. Thus, Lemma \ref{u21} together with \eqref{CC802} implies that there exists some $T^{**}>T^*$ such that $\forall \tau\in(0,T^{**}]$
\begin{equation*} 
\mathbf{H}_1(\tau)\leq 4M_0,~~
\mathbf{H}_2(\tau)\leq 4M_0,~~ 
\mathbf{H}_3(\tau)\leq 4\underline{\rho}^{-1},
\end{equation*} 
which contradicts the definition of $ T^*$ and implies that $T^*=\infty$. In particular, note that \eqref{CC801} is valid for $0<\tau<T<\infty$, which along with Lemma \ref{u21} also ensures that the classical solution obtained above is unique.

We next illustrate the large time behaviour \eqref{1151}. Observe that \eqref{GN} and \eqref{355} ensure that $\forall p\in[1,6)$ and $\theta=(6-p)/2p$,   
\begin{equation*} 
\begin{split}
\|u\|_{W^{1,p}(\Omega)} 
&\leq C\big(\|\nabla w\|_{L^2(\Omega)}^{\theta}\|\nabla^2w\|_{L^2(\Omega)}^{1-\theta}
+\|\nabla w\|_{L^2(\Omega)}+\|\rho-\bar{\rho}\|_{L^p(\Omega)}\big)\\
&\leq C\bigg(\big(\|\nabla u\|_{L^2(\Omega)}
+\|\rho-\bar{\rho}\|_{L^2(\Omega)}\big)^\theta\|\rho\dot{u}\|_{L^2(\Omega)}^{1-\theta}
+\|\nabla u\|_{L^2(\Omega)}+\|\rho-\bar{\rho}\|_{L^6(\Omega)}\bigg),
\end{split}
\end{equation*}
which along with Lemmas \ref{T71}--\ref{21} gives \eqref{1151} and finishes the proof of Theorem \ref{T1}.\thatsall
\\

\textit{Proof of Theorem \ref{T2}.} Let us apply the method in \cite{lx}. Suppose that $X_t(x)$ is the flow line originated from $x\in\Omega$, then we utilize \eqref{11}$_1$ and \eqref{1118} to deduce that $\forall t\geq 0$,
\begin{equation}\label{93}
\rho\big(X_t(x)\big)=\rho_0\big(X_0(x)\big)\exp\big(-\int_0^t\div u(X_{\tau}(x),\tau)\,d\tau\big).
\end{equation}
If there is some $\xi\in \Omega\setminus\Omega_{10}$ such that $\rho_0(\xi)=0$, in view of \eqref{93}, we have $\rho\big(X_t(\xi)\big)\equiv 0$, $\forall t\geq 0$.
Thus, by virtue of \eqref{GN}, we also declare that $\forall r\in(3,\infty)$ and $\theta'=(2r-6)/(5r-6)$,
\begin{equation*}
\bar\n\leq\|\rho-\bar{\rho}\|_{C(\overline{\Omega})}
\leq C\|\rho-\bar{\rho}\|_{L^2(\Omega)}^{\theta'}\|\na \rho\|_{L^r(\Omega)}^{1-\theta'},
\end{equation*}
which together with Lemma \ref{41} implies that
\begin{equation*}
\|\nabla\rho (\cdot,t)\|_{L^r(\Omega)}\geq \hat{\mathbf{C}} e^{\hat{\eta} t}.
\end{equation*}
The proof of Theorem \ref{T2} is therefore completed.\thatsall

\section*{Appendix}
\qquad Let us finish the proof of Lemma \ref{le52}, and the arguments below are based on the standard potential estimates, see also \cite[Appendix]{PP1}. We recall that $v$ is given by the Lam\'{e} system
\begin{equation}\label{CA601}
\begin{cases} 
\Delta v+\nabla\mathrm{div}v=\nabla(\rho-\bar{\rho})&  \mathrm{in}\ \Omega, \\
v=0 & \mathrm{on}\ \partial\Omega.
\end{cases}
\end{equation}
If we define $F_v\triangleq 2\,\mathrm{div}v-(\rho-\bar{\rho})$, then $v$ solves the Dirichlet problem of the Poisson equation,
\begin{equation*}
\begin{cases}
\Delta v=\nabla\big((\rho-\bar{\rho})/2-F_v\big)\quad &\mathrm{in}\ \Omega, \\
v=0\quad &\mathrm{on}\ \partial\Omega.
\end{cases}
\end{equation*}
Consequently, we apply the Green function $G(x,y)$  of $\Delta$ on $\Omega$ for the Dirichlet problem to derive that $\forall x\in\Omega$,
\begin{equation}\label{CC203}
v(x)=\int_\Omega\partial_{y}G(x,y)\cdot\big((\rho-\bar{\rho})/2-F_v\big)\,dy.
\end{equation}

Suppose that the domain $\Gamma$ satisfies  
$\partial\Omega\subset\Gamma\subset\Omega$ and $\mathrm{dist}\,(\partial\Gamma\setminus\partial\Omega,\,\partial\Omega)>d$, then we utilize the method in \cite{PP1} to obtain the standard $C^\alpha$ estimates of $F_v$.

\textbf{Lemma.}
\textit{For any $\alpha\in(0,1)$, there is a constant $C$ depending only on $\alpha$, $d$, and $\Omega$ such that}
\begin{equation}\label{C92}
\|F_v\|_{C^\alpha(\Omega)}\leq C\big(\|\rho\|_{C^\alpha(\Gamma)}+\|\rho-\bar{\rho}\|_{L^2(\Omega)}\big).
\end{equation}

\begin{proof}
We consider the Poisson equation
\begin{equation}\label{C601}
\begin{cases}
\Delta v_1=\nabla(\rho-\bar{\rho})/2~~\mathrm{in}\ \Omega,\\
v_1=0~~\mathrm{on}\ \partial\Omega.
\end{cases}
\end{equation}
Then direct calculations yield that
\begin{equation}\label{CC602}
\begin{split}
F_{v_1}(x)&\triangleq\mathrm{div}v_1-(\rho-\bar{\rho})/2\\
&=\frac{1}{2}\int_\Omega\partial_{x_i}\partial_{y_i}G^*(x,y)\cdot(\rho-\bar{\rho})\,dy,
\end{split}
\end{equation}
where $G^*(x,y)=\Gamma(x-y)-G(x,y)$ and $\Gamma(x)$ is the fundamental solution of $\Delta$ in $\mathbb{R}^3$.

Set $v_2\triangleq v-v_1$, then we check that $v_2$ solves the Lam\'{e} system
\begin{equation}\label{C604}
\begin{cases} 
\Delta v_2+\nabla\mathrm{div}v_2=-\nabla F_{v_1}&  \mathrm{in}\ \Omega, \\
v_2=0 & \mathrm{on}\ \partial\Omega,
\end{cases}
\end{equation}
and it holds that $F_v=F_{v_1}+\mathrm{div}v_2$.  

Let us first treat $F_{v_1}$. According to \cite{soln1,soln2}, for any $\alpha,\beta\in\mathbb{N}$,
there is a constant $C$ depending only on $\alpha$, $\beta$, and $\Omega$, such that
\begin{equation}\label{C603}
|\partial^{\alpha}_x\partial^{\beta}_y G(x,y)|\leq C|x-y|^{-1-\alpha-\beta}.
\end{equation}
We define the domain $D=\{x\in\Omega|\ \mathrm{dist}(x,\partial\Omega)<d/200\}$ and introduce a truncation function $0\le\eta\in C^\infty$ satisfying
\begin{equation}\label{C47}
\begin{cases}
\eta(x)=1 & \mathrm{ for }\ \mathrm{dist}(x,\partial\Omega)\leq d/100,\\
\eta(x)=0 & \mathrm{ for }\ \mathrm{dist}(x,\partial\Omega)\geq d/50.
\end{cases}
\end{equation}

If $x\in D$, in view of \eqref{CC602}, we can split $F_{v_1}$ into two parts 
\begin{equation*}
\begin{split}
F_{v_1}(x) 
&=\frac{1}{2}\int_\Omega\partial_{x_i}\partial_{y_i}G^{*}(x,y)\cdot\big(\eta\cdot(\rho -\bar{\rho})\big)dy\\
&\quad+\frac{1}{2}\int_\Omega\partial_{x_i}\partial_{y_i}G^{*}(x,y)\cdot\big((1-\eta)\cdot(\rho-\bar{\rho})\big)dy\\
&=\hat{F}_{v_1}(x)+R_{v_1}(x).
\end{split}
\end{equation*}

By the same process as \eqref{C601}--\eqref{CC602}, we deduce that $\hat{F}_{v_1}=\mathrm{div}\hat{v}_1-(\rho-\bar{\rho})/2$, where $\hat{v}_1$ solves
the Poisson equation:
\begin{equation*}
\begin{cases}
\Delta\hat{v}_1=\nabla\big(\eta\cdot(\rho-\bar{\rho})/2\big)~~\mathrm{in}\ \Omega,\\
\hat{v}_1=0~~\mathrm{on}\ \partial\Omega.
\end{cases}
\end{equation*}
According to the standard $C^\alpha$ estimates of the elliptic equation(\cite[Theorems 8.1-8.3]{adn},  \cite[Theorems 4.8-4.12]{GT}), we argue that there is a constant $C$ depending on $\Omega$ such that
\begin{equation*}
\|\nabla\hat{v}_1\|_{C^\alpha(\Omega)}\leq C\|\eta\cdot(\rho-\bar{\rho})\|_{C^\alpha(\Omega)}\leq C\|\rho-\bar{\rho}\|_{C^\alpha(\Gamma)},
\end{equation*}
which also leads to the fact
\begin{equation}\label{72}
\|\hat{F}_{v_1}\|_{C^\alpha(\Omega)}\leq C\big(\|\nabla\hat{v}_1\|_{C^\alpha(\Omega)}+\|\eta(\rho-\bar{\rho})\|_{C^\alpha(\Omega)}\big)
\leq C\|\rho-\bar{\rho}\|_{C^\alpha(\Gamma)}.
\end{equation}

In addition, observe that $x\in D$ and \eqref{C47} ensures that the integral
\begin{equation*} 
R_{v_1}(x)=\frac{1}{2}\int_\Omega\partial_{x_i}\partial_{y_i}G^{*}(x,y)\cdot\big((1-\eta)\cdot(\rho-\bar{\rho})\big)dy
\end{equation*}
is out of singularity. Therefore, for $k\in\mathbb{N}$, we apply \eqref{C603} to derive that there is a constant $C$ depending on $k$ and $\Omega$ such that
\begin{equation*}
\begin{split}
&\big|\nabla^k R_{v_1}(x)\big|\\
&=\frac{1}{2}
\left|\int_\Omega\nabla^{k+2}G^{*}(x,y)\cdot\big((1-\eta)\cdot(\rho-\bar{\rho})\big)dy\right|\\
&\leq C\int_{\{|x-y|>d/200\}}\frac{|\rho(y)-\bar{\rho}|}{|x-y|^{k+3}}dy\\
&\leq C\|\rho-\bar{\rho}\|_{L^2(\Omega)},
\end{split}
\end{equation*}
which combining with \eqref{72} gives
\begin{equation}\la{93q}
  \|F_{v_1}\|_{C^\alpha(D)}\leq C\left(\|\rho\|_{C^\alpha(\Gamma)}+\|\rho-\bar{\rho}\|_{L^2(\Omega)}\right).
\end{equation}

If $x\in\Omega\setminus D$, we take $\mathrm{div}(\cdot)$ on both sides of \eqref{C601} and find that
$$\Delta F_{v_1}=0~~\mathrm{in}\ \Omega.$$
Consequently, the interior estimates of harmonic functions (\cite[Chapter 2]{GT}) together with the $L^p$ estimates of the Poisson equation  (\cite[Chapter 8]{GT}) lead to
\begin{equation}\label{71}
\|F_{v_1}\|_{C^\alpha(\Omega\setminus D)}\leq C \|F_s\|_{L^2(\Omega)} 
\leq C \|\rho-\bar{\rho}\|_{L^2(\Omega)}.
\end{equation}
 
Moreover, we utilize the standard $C^\alpha$ estimates of the elliptic system (\cite[Theorems 8.1-8.3]{adn}, \cite[Theorems 4.8-4.12]{GT}) to \eqref{C604} and declare that
\begin{equation*}
\|\nabla v_2\|_{C^\alpha(\Omega)}\leq C\|F_{v_1}\|_{C^\alpha(\Omega)}
\leq C\big(\|\rho\|_{C^\alpha(\Gamma)}+\|\rho-\bar{\rho}\|_{L^2(\Omega)}\big).
\end{equation*}
which together with \eqref{93q} and \eqref{71} gives \eqref{C92} and finishes the proof.
\end{proof}

Now we finish the proof of Lemma \ref{le52}. Recalling that  $\partial\Omega\subset\Gamma'\subset\Gamma
\subset\Omega$, $\mathrm{dist}(\Gamma',\partial\Gamma\setminus\partial\Omega)>d$, and $\mathrm{dist}\,(\partial\Gamma'\setminus\partial\Omega,\,\partial\Omega)>d$, and the modified norm $\|\rho\|_{C^\alpha(\Gamma,\Gamma')}$ is given by 
\begin{equation*} 
\|\rho\|_{C^\alpha(\Gamma,\Gamma')}
\triangleq\sup_{x\in\Gamma,\ y\in\overline{\Gamma'}}
\frac{|\rho(x)-\rho(y)|}{|x-y|^\alpha}+\sup_{x\in \Gamma}|\rho(x)|.
\end{equation*}
Then we declare the following $C^\alpha$ estimates.
 
\textbf{Lemma \ref{le52}.} \textit{Assume that $v$ is given by \eqref{CA601}, then for  
$x\in \Gamma$, $y\in \Gamma'$, and $\varepsilon\in(0,1/1000)$,
there is a constant $C$ determined by $d$ and $\Omega$ such that}
\begin{equation*}
\frac{|v(x)-v(y)|}{|x-y|}\leq C\,\varepsilon^{\frac{\alpha}{2}}\|\rho\|_{C^\alpha(\Gamma,\Gamma')}+
C\,\varepsilon^{-2}\|\rho-\bar{\rho}\|_{L^4(\Omega)}.
\end{equation*}
\begin{proof}
According to \eqref{CC203}, we obtain that
\begin{equation*}
\begin{split}
v(x) 
&=\frac{1}{2}\int_\Omega\partial_yG(x,y)
  \cdot(\rho-\bar{\rho})\,dy
 -\int_\Omega\partial_yG(x,y)\cdot F_v\,dy 
\triangleq V_1(x)+V_2(x).
\end{split}
\end{equation*}

\textit{Step 1. Estimates of $V_1(x)$.} 

If $|x-y|\geq \varepsilon d$, we make use of \eqref{C603} to infer that
\begin{equation}\label{61888}
\begin{split}
\frac{|V_1(x)-V_1(y)|}{|x-y|} 
&=\frac{1}{2|x-y|}\left|\int_\Omega\big(\partial_zG(x,z)-\partial_zG(y,z)\big)
\cdot(\rho(z)-\bar{\rho})\,dz\right|\\
&\leq\frac{C}{|x-y|}\cdot \|\partial_z G\|_{L^{\frac{4}{3}}(\Omega)}\|\rho-\bar{\rho}\|_{L^4(\Omega)} \leq C\varepsilon^{-1}\|\rho-\bar{\rho}\|_{L^4(\Omega)}.
\end{split}
\end{equation}

If $|x-y|<\varepsilon d,$ we calculate that
\begin{equation*}
\begin{split}
\frac{|V_1(x)-V_1(y)|}{|x-y|}
&\leq\frac{C}{|x-y|}\left|\int_{B^c_{10\varepsilon d}(y)\cap\Omega}\big(\partial_zG(x,z)-\partial_zG(y,z)\big)
(\rho(z)-\bar{\rho})dz\right|\\
&\quad+\frac{C}{|x-y|}\left|\int_{B_{10\varepsilon d}(y)\cap\Omega}\big(\partial_zG(x,z)-\partial_zG(y,z)\big)
(\rho(z)-\bar{\rho})dz\right|\\
&\triangleq I_1+I_2.
\end{split}
\end{equation*}
By virtue of \eqref{C603}, $I_1$ is directly handled by:
\begin{equation}\label{C520}
\begin{split}
I_1&=\frac{C}{|x-y|}\left|\int_{B^c_{10\varepsilon d}(y)\cap\Omega}\big(\partial_zG(x,z)-\partial_zG(y,z)\big)\cdot
(\rho(z)-\bar{\rho})dz\right|\\
&\leq\frac{C}{|x-y|}\int_{\{|z|>\varepsilon d\}}\frac{|x-y|}{|z|^3}\cdot|\rho(z)-\bar{\rho}|dz\\
&\leq C\varepsilon^{-\frac{3}{2}}\|\rho-\bar{\rho}\|_{L^2(\Omega)}.
 \end{split}
\end{equation}
In addition, we further decompose $I_2$ into two parts.
\begin{equation*}
\begin{split}
I_2&=\frac{C}{|x-y|}\left|\int_{B_{10\varepsilon d}(y)\cap\Omega}\big(\partial_zG(x,z)-\partial_zG(y,z)\big)
(\rho(z)-\bar{\rho})\,dz\right|\\
&\leq C\,\frac{|\rho(y)-\bar{\rho}|}{|x-y|}\,
\left|\int_{B_{10\epsilon d}(y)\cap\Omega}\big(\partial_zG(x,z)-\partial_zG(y,z)\big)\,dz\right|\\
&\quad+\frac{C}{|x-y|}\left|\int_{B_{10\varepsilon d}(y)\cap\Omega}\big(\partial_zG(x,z)-\partial_zG(y,z)\big)\cdot\big
(\rho(z)-\rho(y)\big)\,dz\right|\\
&\triangleq J_1+J_2.
\end{split}
\end{equation*}

To handle $J_1$, we observe that the regularity assumption on $\Omega$ ensures that there is a constant $C_\Omega$ determined by the shape of $\Omega$, such that $\forall r\leq\mathrm{diam}(\Omega)$ and $\forall x\in\Omega$, 
$$C_\Omega^{-1}\,r^3\leq|B_r(x)\cap\Omega|\leq C_\Omega\,r^3.$$
Consequently, we apply \eqref{D51} to argue that
\begin{equation}\label{5522}
\begin{split}
|\rho(y)-\bar{\rho}|&\leq
\frac{C_\Omega}{|B_{\varepsilon d}|}\int_{B_{\varepsilon d}(y)\cap\Omega}\big(|\rho(y)-\rho(x)|+|\rho(x)-\bar{\rho}|\big)\,dx\\
&\leq\frac{C\,\|\rho\|_{C^\alpha(\Gamma,\Gamma')}}{|B_{\varepsilon d}|}\,\bigg(
\int_{B_{\epsilon d}}\left|x\right|^\alpha dx  \bigg)+\bigg(\frac{C}{|B_{\epsilon d}|}\int_{\Omega}|\rho(x)-\bar{\rho}|^2dx\bigg)^\frac{1}{2}\\
&\leq C \varepsilon^\alpha\|\rho\|_{C^\alpha(\Gamma,\Gamma')}+C\varepsilon^{-\frac{3}{2}}\|\rho-\bar{\rho}\|_{L^2(\Omega)}.
\end{split}
\end{equation}
In addition, since $G(x,z)=0$ for any $x\in\Omega, z\in\partial\Omega$ (see \cite[Section 2.4]{GT}), we argue that
\begin{equation*}
\int_{\Omega}\partial_zG(x,z)\,dz=0~~ \forall x\in\Omega,
\end{equation*}
which along with \eqref{C603} also yields that
\begin{equation}\label{6188}
\begin{split}
&\frac{1}{|x-y|}\left|\int_{B_{10\varepsilon d}(y)\cap\Omega}
\big(\partial_zG(x,z)-\partial_zG(y,z)\big)\,dz\right|\\
&=\frac{1}{|x-y|}\left|\int_{B^c_{10\varepsilon d}(y)\cap\Omega}\big(\partial_zG(x,z)-\partial_zG(y,z)\big)\,dz\right|\\
&\leq\frac{1}{|x-y|}\int_{B^c_{10\varepsilon d}(y)\cap\Omega}\big|\partial_zG(x,z)-\partial_zG(y,z)\big|\,dz\\
&\leq\frac{C}{|x-y|}\int_{\{\varepsilon d\leq|z|\leq\mathrm{diam}(\Omega)\}}\frac{|x-y|}{|z|^3}dz\\
&\leq C\big(1-\log\varepsilon\big).
\end{split}
\end{equation}
Thus combining \eqref{5522} and \eqref{6188} gives
\begin{equation}\label{C524}
J_1\leq C\,\varepsilon^{\alpha/2}\|\rho\|_{C^\alpha(\Gamma,\Gamma')}
+C\,\varepsilon^{-2}\|\rho-\bar{\rho}\|_{L^2(\Omega)}.
\end{equation}

The last term $J_2$ is treated by the similar method, and we declare that
\begin{equation}\label{C525}
\begin{split}
J_2&=\frac{C}{|x-y|}\left|\int_{B_{10\varepsilon d}(y)}\big(\partial_zG(x,z)-\partial_zG(y,z)\big)\cdot\big
(\rho(z)-\rho(y)\big)dz\right|\\
&\leq\frac{C}{|x-y|}\bigg|\int_{B_{5|x-y|}(y)}\big(\partial_zG(x,z)-\partial_zG(y,z)\big)
\cdot\big(\rho(z)-\rho(y)\big)\,dz\bigg|\\
&\quad+\frac{C}{|x-y|}\bigg|\int_{B_{10\varepsilon d}(y)\setminus
B_{5|x-y|}(y)}\big(\partial_zG(x,z)-\partial_zG(y,z)\big)
\cdot\big(\rho(z)-\rho(y)\big)\,dz\bigg|\\
&\leq \frac{C}{|x-y|}\int_{B_{5|x-y|}(y)}(\big|\partial_zG(x,z)\big|+\big|\partial_zG(y,z)\big|)\cdot\big
|\rho(z)-\rho(y)\big|\,dz\\
&\quad +\frac{C}{|x-y|}\int_{B_{10\varepsilon d}(y)\setminus
B_{5|x-y|}(y)}\big|\partial_zG(x,z)-\partial_zG(y,z)\big|\cdot
\big|\rho(z)-\rho(y)\big|\,dz.\\
\end{split}
\end{equation}
In view of \eqref{D51} and \eqref{C603}, the first term is controlled by
\begin{equation}\label{C609}
\begin{split}
&\frac{C}{|x-y|}\int_{B_{5|x-y|}(y)}(\big|\partial_zG(x,z)\big|+\big|\partial_zG(y,z)\big|)\cdot\big
|\rho(z)-\rho(y)\big|\,dz\\
&\leq C\|\rho\|_{C^\alpha(\Gamma,\Gamma')}\cdot|x-y|^{\alpha-1}\cdot\int_{|z|\leq 10|x-y|} |z|^{-2} \,dz\\
&\leq C\,\varepsilon^\alpha\|\rho\|_{C^\alpha(\Gamma,\Gamma')},
\end{split}
\end{equation}
while the second term is bounded by
\begin{equation}\label{CC610}
\begin{split}
&\frac{C}{|x-y|}\int_{B_{10\varepsilon d}(y)\setminus
B_{5|x-y|}(y)}\big|\partial_zG(x,z)-\partial_zG(y,z)\big|\cdot
\big|\rho(z)-\rho(y)\big|\,dz\\
&\leq C\|\rho\|_{C^\alpha(\Gamma,\Gamma')}\cdot \int_{|x-y|<|z|<10\varepsilon d}|z|^{-3+\alpha}\,dz\\
&\leq C\,\varepsilon^\alpha\|\rho\|_{C^\alpha(\Gamma,\Gamma')}.
\end{split}
\end{equation}
Substituting \eqref{C609} and \eqref{CC610} into \eqref{C525} which together with \eqref{61888}, \eqref{C520} and \eqref{C524} leads to
\begin{equation}\label{C0526}
\frac{|V_1(x)-V_1(y)|}{|x-y|}\leq C\,\varepsilon^{\alpha/2}\|\rho\|_{C^\alpha(\Gamma,\Gamma')}+
C\,\varepsilon^{-2} \|\rho-\bar{\rho}\|_{L^4(\Omega)}.
\end{equation}

\textit{Step 2. Estimates of $V_2(x)$.}  

In fact $V_2(x)$ is given by the elliptic system
\begin{equation}\label{C605}
\begin{cases} 
\Delta V_2=-\nabla F_{v_1} &\mathrm{in}\ \Omega, \\
V_2=0 &\mathrm{on}\ \partial\Omega.
\end{cases}
\end{equation}
Let us replace $\alpha$ by $\alpha/2$ in \eqref{C92} and utilize the standard $C^\alpha$ estimates of \eqref{C605} (see \cite[Theorems 4.15 \& 4.16]{GT}) to deduce that
\begin{equation}\label{C48}
\|\nabla V_2(x)\|_{L^\infty(\Omega)}
\leq C\,\|F_{v_1}\|_{C^{\alpha/2}(\Omega)} 
\leq C\big(\|\rho\|_{C^{\alpha/2}(\Gamma')}
+\|\rho-\bar{\rho}\|_{L^2(\Omega)}\big).
\end{equation}
In particular, $\forall x,y\in\Gamma'$, if $|x-y|>\varepsilon d$, we make use of \eqref{5522} to derive that
\begin{equation}\label{C50}
\begin{split}
\frac{|\rho(x)-\rho(y)|}{|x-y|^{\alpha/2}}
&\leq C\varepsilon^{-\alpha/2}\|\rho-\bar{\rho}\|_{L^\infty(\Gamma')}\leq C\big(\varepsilon^{\alpha/2}\|\rho\|_{C^\alpha(\Gamma,\Gamma')}
+\,\varepsilon^{-2}\|\rho-\bar{\rho}\|_{L^2(\Omega)}\big).
\end{split}
\end{equation}
While if $|x-y|\leq\varepsilon d$, we directly check that
\begin{equation}\label{C51}
\begin{split}
\frac{|\rho(x)-\rho(y)|}{|x-y|^{\alpha/2}}&=
\frac{|\rho(x)-\rho(y)|}{|x-y|^\alpha}\cdot|x-y|^{\alpha/2}\leq C\,\varepsilon^{\alpha/2}\|\rho\|_{C^\alpha(\Gamma,\Gamma')}.
\end{split}
\end{equation}
Thus, combining \eqref{C48}--\eqref{C51} leads to
\begin{equation*} 
\|\nabla V_2(x)\|_{L^\infty(\Omega)}\leq C\,\varepsilon^{\alpha/2}\|\rho\|_{C^\alpha(\Gamma,\Gamma')}
+C\,\varepsilon^{-2}\|\rho-\bar{\rho}\|_{L^2(\Omega)},
\end{equation*}
which along with \eqref{C0526} yields
\begin{equation*}
\frac{|v(x)-v(y)|}{|x-y|}\leq C\,\varepsilon^{\alpha/2}\|\rho\|_{C^\alpha(\Gamma,\Gamma')}+
C\,\varepsilon^{-2}\|\rho-\bar{\rho}\|_{L^4(\Omega)}.
\end{equation*}
The proof of Lemma \ref{le52} is therefore completed.
\end{proof}

\section*{Acknowledgements}    The research  is
partially supported by the National Center for Mathematics and Interdisciplinary Sciences, CAS,
NSFC Grant Nos.  11688101 and 12071200,  and Double-Thousand Plan of Jiangxi Province (No. jxsq2019101008), and the Basic Science Center Program (No: 12288201) of the National Natural Science Foundation of China.

\section*{Statements and Declarations}

~~~~\,\textbf{Conflict of interest statement} On behalf of all authors, the corresponding author states that there is no conflict of interest.

\textbf{Data availability statement} No new data were created or analysed in this study.

\end{document}